%
%
%
%
%
%
%
%
%
%
\scrollmode
\magnification=\magstep1
\parskip=\smallskipamount
\def\demo#1:{\par\noindent\it{#1}. \rm}
\def\noi{\noindent}               
\def\ll{\leftline}
\def\cl{\centerline}

\def\begin{\ll{}\vskip 10mm \nopagenumbers}  
\def\pn{\footline={\hss\tenrm\folio\hss}}   
\def\ii#1{\itemitem{#1}}

%
%
\outer\def\beginsection#1\par{\bigskip
  \message{#1}\leftline{\bf\&#1}
  \nobreak\smallskip\vskip-\parskip\noindent}
%
%
\outer\def\proclaim#1:#2\par{\medbreak\vskip-\parskip
    \noindent{\bf#1.\enspace}{\sl#2}
  \ifdim\lastskip<\medskipamount \removelastskip\penalty55\medskip\fi}

\def\endpr{\hfill $\spadesuit$ \medskip}

%
%
%
%
\def\rH{{\rm H}}
\def\rH{{\rm H}}
%
%
%
%

\def\C{{\bf C}}

\def\N{{\bf N}}

\def\R{{\bf R}}

\def\Z{{\bf Z}}


%
%
%
%

\def\cC{{\cal C}}

\def\cE{{\cal E}}

\def\cN{{\cal N}}
\def\cO{{\cal O}}

\def\cU{{\cal U}}

\def\cW{{\cal W}}

%
%
%
\def\a{\alpha}
\def\b{\beta}
\def\g{\gamma}
\def\d{\delta}
\def\e{\epsilon}
\def\z{\zeta}

\def\l{\lambda}
\def\r{\rho}

\def\u{\upsilon}
\def\c{\chi}

\def\L{\Lambda}
%
%
%
%
\def\bar{\overline}              
\def\bs{\backslash}              
\def\di{\partial}                
\def\dibar{\bar\partial}         
\def\hra{\hookrightarrow}

%
%

\def\cn{{\bf C}^n}

\def\cN{{\bf C}^N}

\def\cthree{{\bf C}^3}
\def\c*{{\bf C}^*}
%
%

\def\dim{{\rm dim}\,}                    
\def\holo{holomorphic}                   
\def\nbd{neighborhood}                   
\def\psc{pseudoconvex}                   
\def\spsc{strongly\ pseudoconvex}        
\def\ra{real-analytic}                   
\def\spsh{strongly\ plurisubharmonic}
\def\tr{totally real}                    
\def\pc{polynomially convex}             
\def\ss{\subset\!\subset}                

\def\supp{{\rm supp}\,}                  

\def\iff{if and only if}
\def\hvf{holomorphic vector field}

\def\Aut{{\rm Aut}}                         
\def\istar{i^*}

%
%
\def\cm{{\bf C}^m}
\def\Rn{{\bf R}^n}
\def\Rm{{\bf R}^m}
\def\RN{{\bf R}^N}
\def\zn{{\bf Z}^n}
\def\wh{\widehat}
\def\wt{\widetilde}
\def\End{{\rm End}}
\def\pr{{\prime}}             
\def\dpr{{\prime\prime}}      

%
%
%
%
\begin
\cl{SOLVING THE $d$ AND $\dibar$-EQUATIONS IN THIN TUBES AND}
\cl{APPLICATIONS TO MAPPINGS}
\bigskip\medskip
\cl{Franc Forstneri\v c, Erik L\o w and Nils \O vrelid}

\bigskip\bigskip
\noi \bf Abstract. \rm We construct a family of integral kernels for
solving the $\dibar$-equation with $\cC^k$ and H\"older estimates
in thin tubes around totally real submanifolds in $\C^n$ (theorems
1.1 and 3.1). Combining this with the proof of a theorem of Serre
we solve the $d$-equation with estimates for holomorphic forms
in such tubes (theorem 5.1). We apply these techniques and a
method of Moser to approximate $\cC^k$-diffeomorphisms between
totally real submanifolds in $\C^n$ in the $\cC^k$-topology by biholomorphic
mappings in tubes, by unimodular and symplectic biholomorphic mappings,
and by automorphisms of $\C^n$.

\beginsection 1. The results.

Let $\C^n$ denote the complex $n$-dimensional Euclidean space with
complex coordinates $z=(z_1,\ldots,z_n)$. A compact 
$\cC^k$-submanifold $M\subset \C^n$ ($k\ge 1$), with or without boundary,
is {\it totally real\/} if for each $z\in M$ the tangent space 
$T_z M$ (which is a real subspace of $T_z\C^n$) contains no complex line;
equivalently, the complex subspace $T_z^C M=T_z M +i T_z M$ of
$T_z\C^n$ has complex dimension $m=\dim_\R M$ for each $z\in M$.
We denote by ${\cal T}_\d M =\{z\in\C^n\colon d_M(z)<\d \}$
the tube of radius $\d>0$ around $M$; here $|z|$ is the Euclidean norm 
of $z\in\C^n$ and $d_M(z)=\inf\{|z-w|\colon w\in M\}$.

For any open set $U\subset \C^n$ and integers
$p,q\in\Z_+$ we denote by $\cC^l_{p,q}(U)$ the space
of differential forms of class $\cC^l$ and of bidegree
$(p,q)$ on $U$. For each multiindex $\a\in Z_+^{2n}$ we
denote by $\di^{\a}$ the corresponding partial derivative
of order $|\a|$ with respect to the underlying real
coordinates on $\C^n$.

The following is one of the main results of the paper;
for additional estimates see theorem 3.1 in sect.\ 3 below.

%
%
%
\proclaim 1.1 Theorem: Let $M \subset \C^n$ be a closed, totally
real, $\cC^1$-submanifold and let $0<c<1$. Denote by ${\cal T}_\d$
the tube of radius $\d>0$ around $M$. There is a $\d_0 > 0$ 
and for each integer $l\ge 0$ a constant $C_l>0$ such that the 
following hold for all $0<\d\le \d_0$, $p\ge0$, $q\ge 1$ and $l\ge 1$. 
For any $u\in \cC^l_{p,q}({\cal T}_\d)$ with $\dibar u=0$ there is a
$v\in \cC^l_{p,q-1}({\cal T}_{\d})$ satisfying $\dibar v=u$ in
${\cal T}_{c\d}$ and the estimates
$$
   \eqalignno{ \|v\|_{L^\infty({\cal T}_{c\d})} &\le
   C_0\, \d \|u\|_{L^\infty({\cal T}_\d)}; \cr
   ||\di ^{\a}v||_{L^\infty ({\cal T}_{c\d})} &\le
   C_{l}\! \left( \d ||\di^{\a}u||_{L^\infty ({\cal T}_{\d})} +
     \d^{1-|\a|} ||u||_{L^\infty ({\cal T}_{\d})} \right);
     \qquad |\a|\le l. &(1.1) \cr}
$$
If $q=1$ and the equation $\dibar v=u$ has a solution
$v_0\in \cC^{l+1}_{(p,0)}({\cal T}_\d)$, there is a solution
$v\in \cC^{l+1}_{(p,0)}({\cal T}_\d)$ of $\dibar v=u$ 
on ${\cal T}_{c\d}$ satisfying
for $1\le j\le n$ and $|\a|=l$:
$$ ||\di_j\di ^{\a}v||_{L^\infty ({\cal T}_{c\d})} \le
   C_{l+1} \left( \omega(\di_j\di ^{\a}v_0,\d) +
     \d^{-l} ||u||_{L^\infty ({\cal T}_{\d})} \right).
$$

In the last estimate $\omega(f,\d)=\sup\{|f(x)-f(y)| \colon |x-y|\leq \d\}$ 
is the {\it modulus of continuity} of a function; when $f$ is a differential
form on $\C^n$, $\omega(f,t)$ is defined as  the sum of the moduli of
continuity of its components (in the standard basis).
The constants $C_k$ appearing in the estimates are independent 
of $u$ and $\d$ (they depend only on $M$ and $c$).

The solution in theorem 1.1 is obtained by a family of integral kernels,
depending on $\d>0$ and constructed specifically for thin tubes
(and hence is given by a linear solution operator on each
tube ${\cal T}_\d$). Immediate examples show that the gain
of $\d$ in the estimate for $v$ is the best possible.
When $u$ is a $(0,1)$-form (or a $(p,1)$-form), the estimates for
the derivatives of $v$ in (1.1) follow from the sup-norm estimate by
shrinking the tube and applying the interior regularity for the
$\dibar$-operator (lemma 3.2). This is not the case in bidegrees
$(p,q)$ for $q>1$. We refer to section 3 below for further
details.

\pn

Another major result of the paper is theorem 5.1 in section 5
on solving the equation $dv=u$ for holomorphic forms in tubes
${\cal T}_\d$ with precise estimates. Theorem 5.1 is obtained
by using the solutions of the $\dibar$-equation,
provided by theorem 3.1, in the proof of Serre's theorem to
the effect that, on pseudoconvex domains, the de Rham cohomology
groups are given by holomorphic forms.

We now apply these results to the problem of approximating smooth
diffeomorphisms between \tr\ submanifolds in $\C^n$ by biholomorphic
maps in tubes ${\cal T}_\d$ and by holomorphic automorphism of $\C^n$.
The tools developed here give optimal results without any 
loss of derivatives in these approximation problems.

The {\it complex normal bundle\/}\ $\nu_M \to M$ of a
\tr\ submanifold $M\subset\C^n$ is defined as the quotient 
bundle $\nu_M=T\C^n|_{M}/T^C M$. It can be realized as a 
complex subbundle of $T\C^n|_M$ such that $T\C^n|_M =T^C M\oplus \nu_M$.
Given a diffeomorphisms $f\colon M_0 \to M_1$ between \tr\ submanifolds
$M_0,M_1\subset \C^n$, we say that the complex normal bundles
$\pi_i\colon \nu_i\to M_i$ are isomorphic over $f$ if there exists
an isomorphism of $\C$-vector bundles $\phi\colon \nu_0 \to\nu_1$
satisfying $\pi_1\circ \phi=f\circ\pi_0$.

\medskip\noi \bf 1.2 Theorem. \sl
Let $f\colon M_0\to M_1$ be a diffeomorphism of class $\cC^k$ between
compact \tr\ submanifolds $M_0,M_1\subset\C^n$, with or without
boundary ($n\ge 1$, $k\ge 2$). Assume that the complex
normal bundles to $M_0$ and $M_1$ are isomorphic over $f$.
Then there are numbers $\d_0>0$ and $a>0$ such that for
each $\d \in (0,\d_0)$ there exists an injective holomorphic
map $F_\d \colon {\cal T}_\d M_0 \to\C^n$ such that
$F_\d({\cal T}_\d M_0) \supset {\cal T}_{a\d}M_1$ and the
following estimates hold for $0\le r\le k$ as $\d\to 0$:
$$
    ||F_\d|_{M_0} -f||_{\cC^r(M_0)} = o(\d^{k-r}),\quad
    ||F_\d^{-1}|_{M_1} - f^{-1}||_{\cC^r(M_1)} = o(\d^{k-r}).  \eqno(1.2)
$$
\rm

The $\cC^r(M)$-norm is defined as usual by using a finite open
covering of $M$ by coordinate charts and a corresponding partition
of unity. An important aspect of theorem 1.2 is the precise relationship
between the rate of approximation on $M_0$ resp.\ on $M_1$
and the radius $\d$ of the tube on which the approximating
biholomorphic map $F_\d$ is defined. The condition that the complex 
normal bundles are isomorphic over $f$ is necessary since
the derivative of any biholomorphic map, defined near $M_0$ 
and sufficiently close to $f$ in the $\cC^1(M_0)$-norm,
induces such an isomorphism.
If $M_0$ and $M_1$ are contractible (such as arcs or \tr\ discs)
or if they are of maximal real dimension $n$, theorem 1.2 applies 
to any $\cC^k$-diffeomorphism $f\colon M_0\to M_1$.

When all data in theorem 1.2 are real-analytic, $f$ extends
to a biholomorphic map $F$ from a \nbd\ of $M_0$ onto a \nbd\
of $M_1$ (see remark 1 after the proof of theorem 1.2 in sect.\ 4).
In such case we say that $M_0$ and $M_1$ are
{\it biholomorphically equivalent;} such pairs of
submanifolds have identical local analytic properties
in $\C^n$. This is not so if $f$ is smooth but non \ra, for
there exist smooth arcs in $\C^n$ which are complete pluripolar
as well as arcs which are not pluripolar [DF], yet any diffeomorphism
between smooth arcs can be approximated as in theorem 1.2.

We don't know whether in general there exist biholomorphic maps
$F_\d$ in a {\it fixed} open \nbd\ of $M_0$ and satisfying (1.2)
as $\d\to 0$. However, in certain situations we can approximate
diffeomorphisms by global holomorphic automorphisms of $\C^n$.
Recall that a compact set $K\subset \C^n$ is {\it polynomially convex}
if for each $z\in\C^n\bs K$ there is a \holo\ polynomial $P$ on $\C^n$
such that $|P(z)|>\sup\{ |P(x)| \colon x\in K\}$. We denote by
$\Aut\C^n$ the group of all holomorphic automorphisms of $\C^n$.

\medskip\noi\bf Definition 1. \sl
\item{(a)}
A {\it $\cC^k$-isotopy} (or a {\it $\cC^k$-flow})
in $\C^n$ is a family of $\cC^k$-diffeomorphisms
$f_t\colon M_0\to M_t$ $(t\in [0,1])$ between $\cC^k$-submanifolds
$M_t \subset \C^n$ such that $f_0$ is the identity on $M_0$,
and both $f_t(z)$ and ${\di\over \di t} f_t(z)$ are continuous with
respect to $(t,z)\in [0,1]\times M_0$ and of class $\cC^k(M_0)$
in the second variable for each fixed $t\in [0,1]$.

\item{(b)} The isotopy in (a) is said to be {\it \tr}
(resp.\ {\it \pc}) if the submanifold $M_t\subset \C^n$ is
\tr\ (resp.\ compact \pc) for each $t\in [0,1]$.

\item{(c)} The {\it infinitesimal generator} of $f_t$
as in (a) is the time-dependent vector field $X_t$
on $\C^n$ which is uniquely defined along $M_t$ by the equation
${\di\over \di t}f_t(z)=X_t(f_t(z))$ ($z\in M_0$, $t\in [0,1]$).

\item{(d)} A {\it holomorphic isotopy}
(or {\it \holo\ flow}) on a domain $D\subset \C^n$ is a family of
injective \holo\ maps $F_t\colon D\to \C^n$ such that
$F_0$ is the identity on $D$ and such that the maps $F_t(z)$ and
${\di\over\di t}F_t(z)$ are continuous with respect to
$(t,z)\in [0,1]\times D$. Its infinitesimal generator $X_t$,
defined as in (c), is a \hvf\ on the domain $D_t=F_t(D)$
for each $t\in [0,1]$.
\medskip\rm

%
%
%
%
\proclaim 1.3 Theorem:
Let $M_0\subset \C^n$ be a compact $\cC^k$-submanifold
of $\C^n$ ($n\ge 2$, $k\ge 2$). Assume that
$f_t\colon M_0\to M_t \subset \C^n$ $(t\in [0,1])$ is a $\cC^k$-isotopy
such that the submanifold $M_t=f_t(M_0) \subset\C^n$ is \tr\ and \pc\
for each $t\in [0,1]$. Set $f=f_1\colon M_0\to M_1$.
Then there exists a sequence 
$F_j\in \Aut\C^n$ $(j=1,2,3,\ldots)$ such that
$$ \lim_{j\to\infty} ||F_j|_{M_0} -f||_{\cC^k(M_0)}=0, \qquad
   \lim_{j\to\infty} ||F_j^{-1}|_{M_1} - f^{-1}||_{\cC^k(M_1)}=0. \eqno(1.3)
$$

Combining theorem 1.3 with corollary 4.2 from [FR] we obtain:
%
%
\proclaim 1.4 Corollary: Let $f\colon M_0\to M_1$ be
$\cC^k$-diffeomorphism ($k\ge 2$) between compact, \tr, \pc\
submanifolds of $\C^n$ of real dimension $m$. If $1\le m\le 2n/3$,
there exists a sequence $F_j\in \Aut\C^n$ $(j=1,2,3,\ldots)$
satisfying (1.3).

Theorems 1.2 and 1.3 are proved in sect.\ 4 below.
A weaker version of theorem 1.3 (with loss of derivatives)
was obtained in [FL] by applying H\"ormander's $L^2$-method
for solving the $\dibar$-equations in tubes. For a converse to
theorem 1.3 see remark 2 on p.\ 135 in [FL]. When $f$ is a \ra\
diffeomorphism as in theorem 1.3 the approximating sequence
$F_j\in \Aut\C^n$ can be chosen such that it converges to 
a biholomorphic map $F$ in an open \nbd\ of
$M_0$ in $\C^n$ satisfying $F|_{M_0}=f$ [FR].

We now consider the approximation problem for maps
preserving one of the forms
$$ \omega = dz_1\wedge dz_2\wedge\cdots\wedge dz_n,            \eqno(1.4) $$
$$ n=2n',\quad \omega=\sum_{j=1}^{n'} dz_{2j-1}\wedge dz_{2j}. \eqno(1.5) $$
A holomorphic map $F$ between domains in $\C^n$ satisfying
$F^*\omega=\omega$ will be called a {\it holomorphic\ \ $\omega$-map.}
(1.4) is the (standard) {\it complex volume form} on $\C^n$; in this case
$F^* \omega =JF\cdotp \omega$, where $JF$ is the complex Jacobian determinant
of $F$, and $\omega$-maps are called {\it unimodular}. (1.5) is the
{\it standard holomorphic symplectic form}, and holomorphic $\omega$-maps
are called {\it symplectic holomorphic}. We denote
the corresponding automorphism group by
$$ \Aut_\omega \C^n = \{F\in\Aut\C^n \colon F^*\omega=\omega\}. $$

For convenience we state the approximation results for $\omega$-maps
(theorems 1.5, 1.7 and corollary 1.6) only for closed submanifolds;
for an extension to manifolds with boundary see the remark
following theorem 1.7.

%
%
%
\proclaim 1.5 Theorem: Let $\omega$ be any of the forms (1.4), (1.5).
Let $f\colon M_0\to M_1$ be a $\cC^k$-diffeomorphism between closed
\tr\ submanifolds in $\C^n$ $(k,n\ge2)$. Assume that there is
a $\cC^{k-1}$-map $L\colon M_0\to GL(n,\C)$ satisfying
$$ L_z|_{T_z M_0} = df_z,\qquad L_z^*\omega=\omega \qquad (z\in M_0).
                                                            \eqno(1.6) $$
Then for each sufficiently small $\d>0$ there is an injective
holomorphic map $F_\d\colon {\cal T}_\d M_0 \to\C^n$ such that
$F_\d^*\omega=\omega$ and (1.2) holds as $\d\to 0$.
If $M_0, M_1$ and $f$ are \ra\ and if there exists a continuous
$L$ satisfying (1.6), then $f$ extends to a biholomorphic map $F$
on a neighborhood of $M_0$ satisfying $F^*\omega = \omega$.

\noi
The notation $L_z^*\omega$ in (1.6) denotes the pull-back of
the multi-covector $\omega_{f(z)}$ by the $\C$-linear map
$L_z$ (which we may interpret as a map $T_z\C^n \to T_{f(z)} \C^n$).
Clearly (1.6) implies that the complex normal bundles
$\nu_j\to M_j$ are isomorphic over $f$. Denoting by  
$SL(n,\C)$ the special linear group on $\C^n$ 
and by $Sp(n,\C)$) the linear symplectic group on $\C^{2n}$,
we can express the condition in theorem 1.5 as follows:

\item{(*)} {\sl There exists a $\cC^{k-1}$-map
$L\colon M_0\to SL(n,\C)$ (resp.\ $L\colon M_0\to Sp(n,\C)$)
such that $L_z=df_z$ on $T_z M_0$ for each $z\in M_0$.}

\medskip 
The only obvious necessary condition for the approximation
of a $\cC^k$-diffeomorphism $f\colon M_0\to M_1$ by holomorphic
$\omega$-maps is that the complex normal bundles
$\nu_j\to M_j$ are isomorphic over $f$ and that
$f^*(i^*_1\omega) = i_0^*\omega$, where $i_j\colon M_j\hra\C^n$
is the inclusion. Theorem 1.5 reduces this analytic approximation
problem to the geometric problem of finding an extension $L$ of $df$
satisfying (1.6). The regularity of $L$ is not the key point;
it would suffice to assume the existence of a {\it continuous}
$L$ satisfying (1.6) since an argument similar to the one in the
proof of theorem 1.5 for the \ra\ case then allows us to approximate
$L$ by a $\cC^{k-1}$-map satisfying (1.6). We expect that such
extension does not always exist, although we do not have specific
examples. Here are some positive results.

%
%
\medskip\noi\bf 1.6 Corollary. \sl
Let $\omega$ be one of the forms (1.4), (1.5), and let $k,n\ge2$.
Let $f\colon M_0\to M_1$ be a $\cC^k$-diffeomorphism between
closed totally real submanifolds such that the complex normal
bundles to $M_0$ resp.\ $M_1$ in $\C^n$ are isomorphic over $f$
and $f^*\omega=i_0^*\omega$. Then the conclusion of theorem 1.5
holds in each of the following cases:
\item{(i)}   $\dim M_0=\dim M_1=n$,
\item{(ii)}  $\omega=dz_1\wedge\cdots\wedge dz_n$ and
$M_0$ is simply connected,
\item{(iii)} $\omega=dz_1\wedge\cdots\wedge dz_n$ and $\nu_0$
admits a complex line subbundle.
\medskip\rm

In cases (ii) and (iii) we have $f^*\omega=i_0^*\omega=0$
when $m<n$. Finally we present approximation
results for $\omega$-flows. We first introduce
convenient terminology.

%
%
\medskip\noi\bf Definition 2. \sl
Let $\omega$ be a differential form on $\C^n$ and let
$f_t\colon M_0\to M_t \subset\C^n$ $(t\in [0,1])$ be a
$\cC^k$-isotopy with the infinitesimal generator $X_t$
(definition 1).

\item{(a)} $f_t$ is an {\it $\omega$-flow} if the form $f_t^*\omega$
on $M_0$ is independent of $t\in [0,1]$.

\item{(b)} An $\omega$-flow $f_t$ is {\it closed} (resp.\ {\it exact})
if for each $t\in [0,1]$ the pull-back to $M_t$ of the form
$\alpha_t=X_t\rfloor \omega$ (the contraction of $\omega$ by $X_t$)
is closed (resp.\ exact).

\item{(c)} Let $U\subset \C^n$ be an open set and $\omega$
a \holo\ form on $\C^n$. A holomorphic flow
$F_t\colon U\to \C^n$ $(t\in [0,1])$ satisfying
$F^*_t\omega =\omega$ for all $t$ is called a
{\it holomorphic $\omega$-flow.}
\medskip \rm

\demo Remark:
If $d\omega=0$ (this holds for the forms (1.4), (1.5))
then a flow $f_t\colon M_0\to M_t$ is an $\omega$-flow
\iff\ the pull-back of $\alpha_t=X_t\rfloor \omega$ to $M_t$ is
a closed form on $M_t$ for each $t\in [0,1]$. This can be seen
from the following formula for the Lie derivative $L_{X_t}\omega$
([AMR], p.\ 370, Theorem 5.4.1.\ and p.429, Theorem 6.4.8.\ (iv)):
$$ {d\over dt}\bigl(f_t^*\omega \bigr)
    = f_t^*\bigl( L_{X_t} \omega \bigr)
    = f_t^*\bigl( d(X_t\rfloor \omega)
         + X_t\rfloor d\omega \bigr)
    = f_t^*(d\alpha_t).
$$
Hence $f_t^*\omega$ is independent of $t$ \iff\
$d(i^*_t \alpha_t)=0$ on $M_t$ for each $t\in [0,1]$.

%
%
%
%
\proclaim 1.7 Theorem: Let $\omega$ be any of the forms (1.4), (1.5).
Assume that $M_0 \subset \C^n$ is a closed \tr\ submanifold and
$f_t\colon M_0\to M_t \subset \C^n$ $(t\in[0,1)$ is a \tr\
$\omega$-flow of class $\cC^k$ for some $k\ge 2$. Then for each
sufficiently small $\d>0$ there is a holomorphic
$\omega$-flow $F_t^\d \colon {\cal T}_\d M_0 \to \C^n$ $(t\in [0,1])$
such that for $0\le r\le k$ we have the following estimates
as $\d>0$ (uniformly with respect to $t\in [0,1]$):
$$ ||F^\d_t - f_t||_{\cC^r(M_0)} = o(\d^{k-r}), \quad
   ||(F^\d_t)^{-1} - f_t^{-1}||_{\cC^r(M_t)} = o(\d^{k-r}).
$$
If in addition $n\ge 2$ and $f_t$ is an exact $\omega$-flow
which is \tr\ and \pc, there is for each $\e>0$ a \holo\
$\omega$-flow $F_t \in \Aut_\omega\C^n$ such that
for all $t\in[0,1]$
$$ ||F_t - f_t||_{\cC^k(M_0)} <\e, \quad
   ||F_t^{-1} - f_t^{-1}||_{\cC^k(M_t)} <\e.
$$

Theorems 1.5, 1.7 and corollary 1.6 extend to the
following situation. Let $M_0$ be a compact domain in a
\tr\ submanifold $M'_0\subset \C^n$, not necessarily closed
or compact. In particular, $M_0$ may be a \tr\ submanifold
with boundary $\di M_0$ and $M'_0$ a larger submanifold
containing $M_0$. In the context of theorem 1.5 or corollary 1.6
assume that $f\colon M'_0\to M'_1$ is a $\cC^k$-diffeomorphism
between \tr\ submanifolds in $\C^n$ ($k\ge2$) and
$L\colon M'_0\to GL(n,\C)$ is a $\cC^{k-1}$-map
satisfying (1.6) on $M'_0$. Then the conclusion of theorem 1.5
holds for $M_0$: There exist holomorphic $\omega$-maps
$F_\d \colon {\cal T}_\d M_0 \to \C^n$ for all sufficiently small
$\d>0$ satisfying (1.2) as $\d\to 0$. Likewise, if the flow
$f_t$ as in theorem 1.7 is defined on $M'_0$, the conclusion
of that theorem applies on the compact subdomain $M_0 \subset M'_0$.

In our last result we consider the problem of approximating a diffeomorphism
$f\colon M_0\to M_1$ by \holo\ $\omega$-automorphisms of $\C^n$.
Assuming that $M_0$ and $M_1$ are \pc\ we have two necessary conditions
for such approximation:
\item{--} $f^*\omega=i_0^* \omega$, and
\item{--} there is a \tr, \pc\ flow $f_t\colon M_0\to M_t \subset \C^n$
($t\in [0,1]$) with $f_0=Id_{M_0}$ and $f_1=f$.

The second condition is necessary since the group $\Aut_\omega\C^n$
is connected (see [FR]).
When $\dim M_0$ is smaller than the degree of $\omega$ the
first condition is trivial (both sides are zero).
We summarize some of the situations when such an approximation
is possible. Let $\beta$ be a \holo\ form on $\C^n$ satisfying
$d\b=\omega$; when $\omega$ is given by (1.4) we may take
$\beta = {1\over n} \sum_{j=1}^n (-1)^{j-1}
     dz_1\wedge\cdots \widehat{dz_j}\cdots \wedge dz_n,
$
and when $\omega$ is the form (1.5) we may take
$\beta=\sum_{j=1}^{n'} z_{2j-1}dz_{2j}$.

%
%
%
%
\medskip\noi\bf 1.8 Theorem. \sl
Let $n,k\ge 2$. Let $M_0\subset \C^n$ be a compact connected
$\cC^k$-submanifold of dimension $m$ and let $f_t\colon M_0\to M_t$
$(t\in[0,1])$ be a \tr, \pc\ $\cC^k$-flow. Assume either that $\omega$
is the volume form (1.4), $d\b=\omega$, and at least one of the
following four conditions holds:
\item{(i)}    $m\le n-2$;
\item{(ii)}   $m=n-1$ and $\rH^{n-1}(M_0;\R)=0$;
\item{(iii)}  $m=n-1$, $M_0$ is closed and orientable, and
$\int_{M_0} \b =\int_{M_0} f_1^* \b \ne0$;
\item{(iv)}   $m=n$, $M_0$ is closed and
satisfies $\rH^{n-1}(M_0;\R)=0$, and $f_t^*\omega$ is
independent of $t$,

\noi or that $n=2n^\prime$ ($n^\prime\geq 2$), $\omega$ is the form
(1.5), $d\b=\omega$, and at least one of the following three
conditions holds:
\item{(v)}    $M_0$ is an arc;
\item{(vi)}   $M_0$ is a circle and
              $\int_{M_0} \b =\int_{M_0} f_1^* \b$;
\item{(vii)}  $m=2$, $M_0$ is closed and satisfies
    $\rH^{1}(M_0;\R)=0$, and $f_t^*\omega$ is independent of
    $t\in [0,1]$.

\noi Set $f=f_1\colon M_0\to M_1$. Then there is a sequence
$F_j\in \Aut_\omega \C^n$ satisfying (1.3).
\medskip\rm

For \ra\ data theorem 1.7 was proved in [F2] in the symplectic case  
and in [F3] in the unimodular case. In that situation the 
sequence $F_j\in \Aut_\omega \C^n$ can be chosen 
such that it converges to a \holo\ $\omega$-map
$F$ in a \nbd\ of $M_0$.

The paper is organized as follows. In sect.\ 2 we collect
some preliminary material, mostly extensions of certain well
known results. In sect.\ 3 we construct a family of integral
kernels for solving the $\dibar$-equation in tubes and we prove
the stated estimates; we conclude the section by historical
remarks concerning such kernels. In sect.\ 4 we apply theorem 3.1
to prove theorems 1.2 and 1.3. In sect.\ 5 we solve the equation $dv=u$
in tubes, where $u$ is an exact \holo\ form and we find a \holo\
solution $v$ satisfying good estimates. In sections 6
and 7 we prove the results on approximating $\omega$-diffeomorphisms
by \holo\ $\omega$-maps and $\omega$-automorphisms.
At the end of section 4 we also include a correction to [FL].

\demo Acknowledgements:
The first author acknowledges partial support by an NSF grant,
a Vilas Fellowship at the University of Wisconsin--Madison,
and by the Ministry of Science of the Republic of Slovenia.
A part of this work was done during his visit to the
mathematics department of the University of Oslo, and he thanks
this institution for its hospitality. The second and third author
did initial work on this project while visiting the mathematics
departments of the University of Ljubljana and University of Michigan
at Ann Arbor, respectively, and they would also like to thank these
institutions for their hospitality.

\beginsection 2. Geometric preliminaries.

We denote by $X\rfloor v$ the contraction of a form $v$ by a
vector field $X$. We shall use the following version of the
Poincar\'e's lemma ([AMR], p.437, Deformation Lemma 6.4.17.):

%
%
\proclaim 2.1 Lemma: Let $M$ be a $\cC^2$-manifold and $w$ a closed
$\cC^1$ $p$-form on $I\times M$; $I=[0,1]$, $p>0$.
For $t\in I$ let $i_t\colon M\to I\times M$ be the injection
$x\to(t,x)$. Then the $(p-1)$-form
$v=\int_0^1 i_t^*({\di\over \di t} \rfloor w) dt$ on $M$
satisfies $dv=i_1^*w - i_0^* w$. In particular, let
$F\colon I\times M\to N$ be a $\cC^2$-map and $u$ a closed
$\cC^1$ $p$-form on $N$; $p>0$. Setting $f_t=F\circ i_t \colon M\to N$
and $w=F^* u$ we get $dv=f^*_1 u - f^*_0 u$.

We shall apply this to the case when $F$ is a deformation retraction
of a tubular \nbd\ ${\cal T}_\d= {\cal T}_\d M$ of
a submanifold $M\subset \C^n$ onto $M$.
This means that $f_1$ is the identity on ${\cal T}_\d$,
$f_t|_M$ is the identity
for all $t$, and $f_0({\cal T}_\d)=M$. Set $\pi=f_0$.
With $u$ a closed $\cC^1$ $p$-form on ${\cal T}_\d$ and
$v$ as above we get $dv=u-\pi^* u$ in ${\cal T}_\d$.

In the situation that we shall consider we have the following local
description of the retraction $F$. Let $M$ be a $\cC^k$-submanifold
in $\C^n$. For $U$ a small open \nbd\ in $M$ of a point
$z_0\in M$ there is a $\cC^k$-diffeomorphism
$\phi\colon O\to \pi^{-1}(U)$, where $O$ is open in
$\R^m\times\R^{2n-m}$, such that

\item{(i)} $F^{-1}(U) = O\cap (\R^m\times \{0\}^{2n-m})
= O'\times \{0\}^{2n-m}$,

\item{(ii)} for $x'\in O'$, the set
$O_{x'}=\{y'\in \R^{2n-m} \colon (x',y')\in O\}$ is starshaped
with respect to $0$,

\item{(iii)} the map $f_t=F\circ i_t$ is $\phi$-conjugate to
$(x',y')\to (x',ty')$ for each $t\in I=[0,1]$.

Let $u=\sum'_{|I|+|J|=p} u_{I,J}(x',y') {dx'}^I \wedge {dy'}^J$
in these coordinates. Then
$$ w=\sum_{|I|+|J|=p}{}{\!\!\!\!\!\! '}\,\, u_{I,J}(x',ty')
  {dx'}^I \wedge d(ty')^J
$$
%
%
and it is easy to check that
$$ v= \sum_{|I|+|K|=p}{\!\!\!\!\!\!\! '}\,\,   (-1)^{|K|} \sum_{j=1}^n
   \sum_{|J|=|K|+1} \e^{jK}_J \, y'_j
   \left( \int_0^1 u_{I,J}(x',ty')t^{|K|} dt\right)
   {dx'}^I \wedge {dy'}^K,
$$
where $\e^{jK}_J$ equals, if $jK$ is a permutation of $J$,
the signature of that permutation, and equals zero otherwise.

$F$ is constructed by retracting to $M$ along the fibers of a
vector bundle supplementary to the tangent bundle $TM$. The normal bundle
to $M$ in $\C^n$ is an obvious choice, but is only of class $\cC^{k-1}$
when $M$ is a $\cC^k$-submanifold. We shall show below that there are
$\cC^k$-subbundles $E$ of $M\times \C^n$ that are arbitrarily close to
the normal bundle. When $k>1$, it is easy to see that
$(z+E_z)\cap {\cal T}_\d$ is starshaped with respect to $z$ for all $z\in M$
when $\d>0$ is small enough and $E$ is sufficiently close to the
normal bundle. The map $G\colon E\to \C^n$, $G(z,v)= z+v$,
maps the zero section $0_E$ diffeomorphically onto $M$, and its derivative
$dG$ is an isomorphism at each point of $0_E$; hence $G$ is a
$\cC^k$-diffeomorphism of a \nbd\ $U_\d \subset E$ of $0_E$ onto ${\cal T}_\d$
for $\d>0$ small. We may assume that $U_\d\cap E_z$ is starshaped with
respect to $(z,0)$ for each $z\in M$. When $f_t$ is $G$-conjugate to the
map $(z,v)\to (z,tv)$ in $U_\d$ for $t\in I$, the map $F$ has the properties
listed above.

The local coordinates $(x',y')$ are constructed as follows.
Let $\varphi\colon O'\to U \subset M$ be a local $\cC^k$ parametrization,
and $s_1,\ldots,s_{2n-m}$ sections of $E\to M$ over $U$ which
form a $\cC^k$-trivialization of $E|_U$. We set
$\phi(x',y')=\varphi(x')+\sum_{j=1}^{2n-m} y'_j\,
   s_j(\varphi(x'))
$ for $x'\in O'$, $y'\in \R^{2n-m}$,
and restrict it to $O=\phi^{-1}({\cal T}_\d)$. Then the fiber $O_{x'}$
is starshaped for all $x'\in O'$ when $\d>0$ is small enough.
\endpr

\proclaim 2.2 Lemma: {\rm (Approximation of subbundles)} Let $M$
be a $\cC^k$-submanifold of $\C^n$ and $E\to M$ a $\cC^l$-subbundle
(real or complex) of $M\times \C^n$ for some $0\le l < k$.
Then there is a $\cC^k$-subbundle $E'$ of $M\times\C^n$
arbitrarily close to $E$ in the $\cC^l$ topology.
Moreover, if $M$ is \tr\ in $\C^n$ and the bundle $E$ is complex,
$E'$ may be taken as the restriction to $M$ of
a holomorphic subbundle of $U\times \C^n$
for some open \nbd\ $U$ of $M$ in $\C^n$.

\demo Proof: A proof may be based on the following standard
result. If $L\colon M\to {\rm Lin}_\C(\C^n,\C^n)$ is a $\cC^l$ map
such that $L_z$ has constant rank $r$ independent of $z\in M$
(abusing the language we shall say that $L$ has rank $r$), then
$$ E_L=\{(z,v)\in M\times \C^n \colon v\in L_z(\C^n) \}     \eqno(2.1) $$
is a complex $\cC^l$-subbundle of rank $r$ of the trivial bundle
$M\times \C^n$, and every subbundle $E$ of $M\times \C^n$ appears in
this manner, for instance by setting $L_z$ to be the orthogonal projection
of $\C^n$ onto the fiber $E_z$ for $z\in M$. The analogous result
holds for real vector bundles.

A more regular approximation to a subbundle $E$ may then be obtained
by approximating the corresponding map $L$ defining $E$ by a more regular
map of rank $r$. The problem is that the rank of a generic perturbation
of $L$ may increase. To overcome this we use the following result:

\sl Let $C$ be a positively oriented simple closed curve in $\C$,
and let $L\in {\rm Lin}_\C(\C^n,\C^n)$ be a linear map with no
eigenvalues on $C$. Then $\C^n=V_+\oplus V_-$, where $V_+$ resp.\
$V_-$ are $L$-invariant subspaces of $\C^n$ spanned by the generalized
eigenvectors of $L$ inside resp.\ outside of $C$. The map
$$ P(L)= {1\over 2\pi i} \int_C \left(\zeta I-L\right)^{-1}\,
                                                   d\z  \eqno(2.2) $$
is the projection onto $V_+$ with kernel $V_-$ \rm (see [GLR]).
Note that $P(L)$ depends holomorphically on $L$; thus, if $L$
depends $\cC^k$ or holomorphically on a parameter, so does $P(L)$.

We now take $C$ to be a curve which encircles $1$ but not zero;
for instance
$$ C=\{\z\in \C \colon |\z-1|=1/2\}.                      \eqno(2.3)  $$
Let $P$ be the associated projection operator (2.2). If $L$ is a projection
then $P(L)=L$. Moreover, for each $L'$ sufficiently near a projection
$L$, each eigenvalue of $L'$ is either near $0$ or near $1$ and
hence $P(L')$ is a projection with the same rank as $L$.

Thus, to smoothen $E$, let $L_z$ be the orthogonal projection onto
$E_z$ for $z\in M$; we approximate $L$ by a $\cC^k$-map
$L'\colon M\to {\rm Lin}_\C(\C^n,\C^n)$ and let $E'$ be the
bundle (2.1) associated to $P(L')$. By (2.2) the difference equals
$$ P(L')-L= {1\over 2\pi i} \int_C \left(
          (\zeta I-L')^{-1} - (\zeta I-L)^{-1}\right) \, d\z
$$
and is $\cC^l$-small when $L'-L$ is.

In the real case we extend $L\colon \R^n\to\R^n$ to a complex linear
map $L\colon \C^n\to \C^n$ and observe that $P(L)$ is also real
(i.e., it maps $\R^n$ to itself) when $C$ is the curve (2.3).
Hence the restriction of $P(L)$ to $\R^n$ solves the problem.

Let now $M$ be a \tr\ submanifold of $\C^n$ and $E\to M$ a
$\cC^l$ rank $r$ complex subbundle of $M\times \C^n$.
For each $z\in M$ let $L_z\colon \C^n\to E_z$ be the orthogonal
projection onto $E_z$. By [RS] we can approximate the $\cC^l$-map
$L\colon M\to {\rm Lin}_\C(\C^n,\C^n)$ as well as we like
in the $\cC^l$-topology on $M$ by the restriction to $M$
of a holomorphic map $L'\colon U\to {\rm Lin}_\C(\C^n,\C^n)$
defined on an open \nbd\ $U\subset \C^n$ of $M$. Shrinking $U$
we may assume that $L'_z$ has exactly $r$ eigenvalues inside $C$
(2.3) for each $z\in U$, so $P(L'_z)$ is a rank $r$ projection.
The map $z\to P(L'_z)$ is holomorphic in $U$ and determines a holomorphic
rank $r$ vector bundle $E'$ over $U$, with $E'|_M$ close to $M$.
\endpr

Let $d=\di+\dibar$ be the splitting of the exterior derivative
on a complex manifold.

\proclaim Definition 3: {\rm ($\dibar$-flat functions)}
If $M$ is a 
closed subset in a complex manifold $X$ and $u$
is a $\cC^k$-function ($k\ge 1$) defined in a
\nbd\ of $M$ in $X$, we say that $u$ is {\it $\dibar$-flat
(to order $k$)} on $M$ if $\di^\alpha(\dibar u)(z)=0$ for each
$z\in M$ and each derivative $\di^\alpha$ of total order
$|\alpha|\le k-1$ with respect to the underlying real local
coordinates on $X$.

We shall commonly use the phrase `$u$ is a $\dibar$-flat
$\cC^k$-function' when it is clear from the context
which subset $M\subset X$ is meant.

\proclaim 2.3 Lemma: {\rm ($\dibar$-flat partitions of unity)}
Let $M$ be a \tr\ $\cC^k$-submanifold of a complex manifold $X$;
$k\ge 1$. For every open covering ${\cal U}$ of $M$ in $X$ there
exists a $\cC^k$ partition of unity on a \nbd\ of $M$ in $X$, subordinate
to the covering ${\cal U}$ and consisting of functions that are
$\dibar$-flat to order $k$ on $M$.

\demo Proof: We may assume that ${\cal U}$ consists of coordinate
\nbd s. Let $\phi^0_\nu$ be a $\cC^k$ partition of unity subordinate
to ${\cal U}|_M =\{U\cap M\colon U\in {\cal U}\}$. We may assume that
the index sets agree, so $\supp\phi^0_\nu \subset U_\nu$ for each
$\nu$. By passing to local coordinates we may find a $\dibar$-flat
$\cC^k$ extension $\wt{\phi}_\nu$ of $\phi^0_\nu$ with
$\supp \wt{\phi}_\nu \subset U_\nu$. Since
$\rho= \sum_\nu \wt{\phi}_\nu =1$ on $M$, $\rho\ne 0$ in
a \nbd\ $V$ of $M$ in $X$. It is immediate that
$\phi_\nu = \wt{\phi}_\nu/\rho$ is a $\cC^k$ partition of
unity on $V$ which is $\dibar$-flat (to order $k$) on $M$.
\endpr

As a consequence of lemma 2.3 we see that the usual results
about $\dibar$-flat extensions of maps into $\C^N$ are also valid
for \tr\ submanifolds in arbitrary complex manifolds.

\proclaim 2.4 Lemma: {\rm (Asymptotic complexifications)}
Let $M$ be a \tr\ $\cC^k$-submanifold of $\C^n$ of real dimension
$m\le n$; $k\ge 1$. Then there exists a $\cC^k$-submanifold
$\wt{M} \supset M$ in $\C^n$, of real dimension $2m$, with the
following property: $\wt{M}$ may be covered by $\cC^k$ local
parametrizations $Z\colon U\to Z(U)\subset \wt{M}$,
with $U\subset \C^m$ open subsets, such that
$Z^{-1}(M)=U\cap \R^m$ and $Z$ is $\dibar$-flat on
$U\cap \R^m$. Moreover there is a $\cC^k$-retraction of a \nbd\
of $\wt{M}$ in $\C^n$ onto $\wt{M}$ which is $\dibar$-flat on $M$.

\demo Proof: By a theorem of Whitney ([Wh2], Theorem 1, p.\ 654)
there exists a $\cC^\omega$-manifold $M_0$ and a
$\cC^k$-diffeomorphism $G^0\colon M_0\to M$.
The manifold $M_0$ has a complexification $\wt{M}_0$ which is a
complex manifold containing $M_0$ as a maximal real submanifold.
The map $G^0$ has a $\dibar$-flat extension $G\colon \wt{M}_0\to \C^n$
which is an injective immersion at $M$. (To obtain $G$ it suffices to
patch local $\dibar$-flat extensions of $G_0$ by a $\dibar$-flat
partition of unity provided by lemma 2.3). Hence $G$ maps a \nbd\ of
$M_0$ in $\wt{M}_0$ diffeomorphically onto its image
$\wt{M}\subset \C^n$. When $Z^0\colon U\to \wt{M}_0$
($U$ open in $\C^m$) is a local holomorphic parametrization with
$(Z^0)^{-1}(M_0)=U\cap \R^m$, the map $Z=G\circ Z^0 \colon U\to \wt{M}$
is a local parametrization of the type described in lemma 2.4.
Note that $T_z\wt{M}=T^C_z M$ for each $z\in M$.

Next we prove the existence of a retraction onto $\wt{M}$
which is $\dibar$-flat on $M$.
Let $\nu\to M$ be the complex normal bundle of $M$ in $\C^n$.
By lemma 2.2 there is an open \nbd\ $O$ of $M$ in $\C^n$ and
a holomorphic rank $(n-m)$ subbundle $N \subset O\times \C^n$ such
that $N|_M$ approximates $\nu$ well. Shrinking $O$ we may assume
that $N$ is transversal to $\wt{M}$ in $O$. This means that
the map $\phi\colon N|_{\wt{M}} \to \C^n$,
$\phi(z,v) = z+v$, is a $\cC^k$-diffeomorphism from a \nbd\ $W$
of the zero section in $N|_{\wt{M}}$ onto its image
$O_0\subset O \subset \C^n$. We may assume that $W\cap N_z$ is starshaped
with respect to $0_z \in N_z$ for each $z\in \wt{M}$. Now the
deformation retraction $(z,v)\to (z,tv)$ ($t\in [0,1]$) of $W$
onto the zero section in $N|_{\wt{M}}$ may be transported by $\phi$
to a retraction $F\colon [0,1]\times O_0\to O_0$ of $O_0$ onto
the submanifold $\wt{M}\cap O_0$.
Set $\pi=F_0 \colon O_0\to \wt{M} \cap O_0$. Let
$U\subset \C^m$ and let $Z\colon U\to \wt{M}$ be a local
$\cC^k$-parametrization such that $Z(U\cap \R^m) \subset M$ and
$Z$ is $\dibar$-flat on $U\cap \R^m$. Choose \holo\ sections
$s_1,\ldots,s_{n-m}$ of $N$ which provide a trivialization
of $N$ near $Z(U)$. Then
$$ (z',w')\to Z(z')+ \sum_{j=1}^{n-m} w'_j s_j(Z(z')) $$
is a $\cC^k$-diffeomorphism of a \nbd\ $W$ of $U\times \{0\}^{n-m}$
in $\C^n$ onto $\pi^{-1}(Z(U))$, and it is $\dibar$-flat on
$(U \cap \R^m) \times \{0\}^{n-m}$.
In these coordinates the maps $F_t$ are given by
$(z',w')\to (z',tw')$, hence $F_t$ is $\dibar$-flat
on $\pi^{-1}(M)$.
\endpr

%
%
%
\medskip \noi\bf 2.5 Lemma. \sl
{\rm (The rough multiplication)}
Let $U$ be an open set in $\R^N$, $f\in \cC^k(U)$ and
$g\in \cC^{k-1}(U)$, where $k\ge 1$. Let $E$ be a closed
subset of $U$ such that $f(x)=0$ for all $x\in E$.
Then there exists a function $h\in \cC^k(U)$ such that

\item{(i)} $|\di^\a(h-fg)|=o(d_E^{k-|\a|})$ for $|\a|<k$,
uniformly on compacts in $U$,

\item{(ii)} at points of $E$ we have
$\di^\a h=\sum_{0\ne \beta\le \alpha} \left({\a\atop\b}\right)
      \di^\b\! f \, \di^{\a-\b}\! g$
for $|\a|\le k$, and

\item{(iii)} if $U\subset \C^N$ and if $f$ and $g$ as above
are $\dibar$-flat on $E \subset f^{-1}(0)$, then so is $h$.
\medskip\rm

The proof is similar to the better known `Glaeser-Kneser rough
composition theorem'; the main point is to verify that the
collection of functions $(\di^\a h)_{|\a|\le k}$ on $E$, defined by (ii),
are a Whitney system, i.e., they satisfy the assumptions of the
Whitney's extension theorem (see [Wh1] or [T]).
We shall leave out the details of this verification.
Let $h$ be a $\cC^k$-function provided by Whitney's theorem, having
the partial derivatives given by (ii) at points of $E$.
Then (i) follows easily by comparing the Taylor expansions
of $\di^\a h$, $\di^\b f$ and $\di^{\a-\b} g$ about the nearest
point in $E$. The case (iii) follows from (ii) which is seen as follows.
From (ii) we get at points of $E$ and for $|\a|\le k-1$
$$ \di^\a\dibar h= \di^\a(\dibar f\cdotp g) + \sum_{0\ne\b\le \a}
   \left({\a\atop\b}\right) \di^\b\! f\, \di^{\a-\b}\dibar g.
$$
If $f$ and $g$ are $\dibar$-flat on $E$, this expression
vanishes when $|\a|\le k-1$, so we get (iii).
\endpr

The following lemma is needed in the proof of theorem 1.4 and
its corollaries.

\proclaim 2.6 Lemma: Let $M$ be a \tr, $m$-dimensional
$\cC^k$-submanifold of $\C^n$, $f \colon M\to \C^p$ a $\cC^k$-map
and $l\colon M\to {\rm Lin}_\C(\C^n,\C^p)$ a $\cC^{k-1}$-map
such that for each $z\in M$, $l_z$ agrees with $df_z$
on $T_z M$. Then there is a \nbd\ $U\subset \C^n$ of $M$
and a $\cC^k$-map $F\colon U\to\C^p$ which is $\dibar$-flat
on $M$ and satisfies $F(z)=f(z)$ and $dF_z=l_z$ for all $z\in M$.

\demo Proof: It suffices to prove the result for functions
$(p=1)$; the general case then follows by applying it
componentwise. So we shall assume $p=1$.

We first consider the local case. Fix a point $z_0\in M$.
Choose $e_1,\ldots, e_{n-m} \in \C^n$ such that these vectors,
together with the tangent space $T_{z_0}M$, span a \tr\ subspace
of $T_{z_0}\C^n$ of maximal dimension $n$. When $\kappa\colon U\to M$
is a $\cC^k$-parametrization of a small \nbd\ of $z_0$ in $M$,
with $\kappa(0)=z_0$, and $V$ is a sufficiently small \nbd\ of $0$
in $\R^{n-m}$, the map
$\phi(x,y)=\kappa(x)+\sum_{j=1}^{n-m} y_j e_j$ ($x\in U,\ y\in V$)
is a $\cC^k$-diffeomorphism onto an $n$-dimensional \tr\ submanifold
in $\C^n$. Observe that for $x\in U$ and $(u,v)\in \R^m\times \R^{n-m}$
we have
$$ l_{\kappa(x)} \circ d\phi_{(x,0)}(u,v) =
   df_{\kappa(x)} \circ d\kappa_x(u) +
   \sum_{j=1}^{n-m} v_j\, l_{\kappa(x)}(e_j).
$$
Since $l_{\kappa(x)}(e_j)$ is only of class $\cC^{k-1}$ in $x$,
we apply the rough multiplication lemma to the pairs
$y_j$, $l_{\kappa(x)}(e_j)$ for $1\le j\le n-m$ to get a
$\cC^k$ function $h$ on $U\times V$ satisfying
${\di h \over \di x_i}(x,0)=0$,
${\di h \over \di y_j}(x,0)=l_{\kappa(x)}(e_j)$ for
$1\le i\le m$ and $1\le j\le n-m$. With
$F^0(x,y)= f(\kappa(x)) + h(x,y)$ it follows that
$dF^0_{(x,0)}=l_{\kappa(x)}\circ d\phi_{(x,0)}$. When $\wt{F}^0$
resp.\ $\wt{\phi}$ are $\cC^k$-extensions of $F^0$ resp.\ $\phi$
which are $\dibar$-flat on $\R^n$, we see that $\wt{\phi}$
is a $\cC^k$-diffeomorphism of a \nbd\ of $0\in \C^n$ onto a \nbd\
of $z_0\in \C^n$. Thus, near $z_0$, $F=\wt{F}^0 \circ \wt{\phi}^{-1}$
is a $\cC^k$ $\dibar$-flat extension of $f$. When $z\in M$ we have
$dF_z=l_z$ on a maximal \tr\ subspace, so these two linear maps are
equal on $T_z\C^n$. This establishes the local case.

For the global case let $\cU= \{U_i\}$ be an open covering of $M$
and $F^{(i)}$ a $\dibar$-flat extension of $f$ in $U_i$, with
$dF^{(i)}_z=l_z$ for $z\in U_i\cap M$. By lemma 2.3 there is a partition
of unity $\{\phi_i\}$ by $\dibar$-flat $\cC^k$-functions on a \nbd\
of $M$ subordinate to $\cU$. We set $F=\sum_i \phi_i F^{(i)}$,
where the term with index $i$ is zero outside $U_i$. When $z\in M$,
$dF_z= \sum_i \phi_i(z)dF^{(i)}_z + \sum_i f(z)d(\phi_i)_z$.
Since $\sum_i\phi_i=1$, $\sum_i d\phi_i =0$ and we get $dF_z=l_z$.
\endpr

\beginsection 3.
Solving the $\dibar$-equation in tubes around \tr\ manifolds.

In this section we construct a family of integral kernels, depending
on a parameter $\d>0$, for solving the $\dibar$-equation in
tubes ${\cal T}_\d M$ around compact \tr\ submanifolds
$M\subset \C^n$ of class $\cC^1$. The main result is
theorem 3.1 which is identical with theorem 1.1
except that it contains additional H\"older
estimates (3.3) and (3.4).

We denote by $d_M$ the Euclidean distance to $M$. If $M$
is of class $\cC^k$, it is well known that $\r = d_M^2$ is a
$\cC^k$ strictly plurisubharmonic function in a neighborhood of $M$
when $k > 1$, and when $k=1$ there is a strictly plurisubharmonic
$\cC^2$-function $\r$ such that $\r = d_M^2 + o(d_M^2)$.
As in sect.\ 1 let ${\cal T}_\d$ denote the tubular \nbd\
of $M$ of radius $\d$, i.e., the set of points whose distance
to $M$ is less than $\d$.

For a domain $D$ in $\Rn$ (or in $\C^n$), a bounded function $u$
in $D$ belongs to the H\"older class $\L^s(D)$ for some $0<s<1$
if $|u|_{s,D}\colon= {\rm sup} \{ |u(z+h)-u(z)| |h|^{-s}
         \colon h\neq 0,\ z,z+h\in D\} < \infty$;
in this case the H\"older $s$-norm of $u$ is defined by
$\|u\|_{\L^s(D)}=\|u\|_{L^\infty(D)}+|u|_{s,D}$.
When $s=1$ we set
$|u|_{1,D}\colon= {\rm sup}
  \{|u(z+h)+u(z-h)-2u(z)||h|^{-1} \colon h\neq 0,\ z,z-h,z+h\in D\}$;
$\L^1(D)$ is called the {\it Zygmund class} on $D$. When $D$ is a
tubular neighborhood  ${\cal T}_\d M$ of a submanifold $M$, we write
$|u|_{s,\d}$ for $|u|_{s,{\cal T}_\d M}$. When $s=k+\a$, $k\in \Z_+$ and
$0<\a\leq 1$, we take
$\|u\|_{\L^s(D)}=\|u\|_{\cC^k(D)}+|D^ku|_{\a,D}$. We sometimes
write $\cC^{k+\a}(D)$ for $\L^{k+\a}(D)$ when $0<\a<1$.

We extend function space norms to vector fields or differential forms on
open sets in $\Rn$ as the sum of the norms of the components. When $M$ is
a compact $\cC^k$-manifold, we define the norms on functions or forms on $M$
as follows: Let $\Phi_j \colon U_j \to V_j \subset M$, $j=1,\cdots,p$, be a
covering of $M$ by local parametrizations, and $\{\phi_1,\cdots,\phi_p\}$
a $\cC^k$ partition of unity subordinate to the covering $\{V_1,\cdots,V_p\}$
of $M$. Then we set
$\|u\| = \sum_{j=1}^p\|\Phi_j^*(\phi_ju)\|$, where $\|\cdot\|$ is a H\"older
or some other function space norm. Different choices of $\{\Phi_j\}$
and $\{\phi_j\}$ give rise to equivalent norms on the same space.

Let $z=(z_1,\ldots,z_n)$ be the complex coordinates and
$(x_1,y_1,\ldots,x_n,y_n)$ ($z_j=x_j+iy_j$) the underlying real
coordinates on $\C^n=\R^{2n}$. For $1\le j\le 2n$, $\di_j$ denotes the
partial derivative with respect to the $j$-th variable. If
$\a=(\a_1,\ldots,\a_{2n})$ is a multiindex of length $2n$ then
$\di^{\a}$ denotes the corresponding partial derivative of order
$|\a|=\a_1+\cdots+\a_{2n}$ with respect to the real variables
on $\C^n=\R^{2n}$.

If $f$ is a function or a form near $M$, we shall say that $f$
{\it vanishes to order $l$ on $M$} if $|f(z)| = o(d_M(z)^l)$ and
$\di^{\a}f = 0$ on $M$ when $|\a| \le l$. Recall that any
$\cC^k$-function $f$ on $M$ can be extended to a $\cC^k$-function
on $\C^n$ such that $\dibar f$ vanishes to order $k-1$ on $M$
(Lemma 4.3 in [H\"oW]).

We call a continuous function $\omega \colon \R_+\to\R_+$
a {\it modulus of continuity} if it is non-decreasing, sub-additive,
and $\omega (0)=0$. If $f \colon A \to \C$, $A\subset\Rn$, is
uniformly continuous, we define the modulus of continuity
of $f$ by $\omega(f,t)=\sup\{|f(x)-f(y)| \colon |x-y|\leq t\}$,
$t\geq 0$. $\omega(f,\cdot)$ is clearly a modulus of continuity as
defined above. We say that a modulus of continuity $\omega$ is a modulus
of continuity for a function $f$ if $\omega(f,t)\leq\omega(t)$ for
all $t\geq 0$. If $f$ is a form on $A$, $\omega(f,t)$ is defined as
the sum of the moduli of continuity of its components.

We denote by $\cC^l_{p,q}(U)$ the space of $(p,q)$-forms
of class $\cC^l$ on an open set $U\subset \C^n$.

%
%
%
\proclaim 3.1 Theorem:
Let $M \subset \C^n$ be a closed totally real $\cC^1$-submanifold
and let $0<c<1$. Denote by ${\cal T}_\d$ the tube of radius $\d>0$ 
around $M$. Then there is a $\d_0 > 0$ and for each 
integer $l\ge 1$ a constant $C_l>0$ such that 
the following hold for $0<\d\le \d_0$, $p\ge0$, $q\ge 1$:
For each $u\in \cC^l_{p,q}({\cal T}_\d)$ with $\dibar u=0$ there is a
$v\in \cC^l_{(p,q-1)}({\cal T}_{\d})$ satisfying $\dibar v=u$
in ${\cal T}_{c\d}$ and
$$ ||\di ^{\a}v||_{L^\infty ({\cal T}_{c\d})} \le
   C_{l} \left( \d ||\di^{\a}u||_{L^\infty ({\cal T}_{\d})} +
   \d^{1-|\a|} ||u||_{L^\infty ({\cal T}_{\d})} \right),
   \qquad |\a|\le l. \eqno(3.1)
$$
In particular we have
$\|v\|_{L^\infty({\cal T}_{c\d})} \le C\d \|u\|_{L^\infty({\cal T}_\d)}$.
If $q=1$ and the equation $\dibar v=u$ has a solution
$v_0\in \cC^{l+1}_{(p,0)}({\cal T}_\d)$, there is a solution
$v\in \cC^{l+1}_{(p,0)}({\cal T}_\d)$ of $\dibar v=u$ satisfying
for $1\le j\le n$
$$ ||\di_j\di ^{\a}v||_{L^\infty ({\cal T}_{c\d})} \le
   C_{l+1} \left( \omega(\di_j\di ^{\a}v_0,\d) +
     \d^{-l} ||u||_{L^\infty ({\cal T}_{\d})} \right),
     \qquad |\a|=l. \eqno(3.2)
$$
If we assume in addition that $\di^\a u\in \L^s({\cal T}_\d)$ for some
$|\a|\le l$ and $0<s \le 1$, we may choose $v$ as above satisfying
also the following estimates (with constants $C_{l,s}$ independent
of $u$ and $\d$):
$$ \eqalignno{ ||\di_j\di ^{\a}v||_{L^\infty ({\cal T}_{c\d})} &\le
   C_{l,s} \left( \d^s ||\di^{\a}u||_{\L^s ({\cal T}_{\d})} +
              \d^{-|\a|} ||u||_{L^\infty({\cal T}_\d)} \right) & (3.3)\cr
              ||\di_j\di ^{\a}v||_{\L^s ({\cal T}_{c\d})} &\le
   C_{l,s} \left( ||\di^{\a}u||_{\L^s ({\cal T}_{\d})} +
              \d^{-|\a|-s} ||u||_{L^\infty({\cal T}_\d)} \right) & (3.4)\cr}
$$

\demo Remarks: 1.\ If $u$ is of class $\cC^l$, there is in general 
no $\cC^{l+1}$ solution $v$ to $\dibar v=u$.

\noi 2.\ In (3.3) one may be tempted to delete $\d^s$ and use instead
the $L^\infty({\cal T}_\d)$ norm of $\di^\a u$ in the first term on the right
hand side. This however is false even when $n=1$ and is a well known
phenomenon. Since the Bochner-Martinelli operator used in the
proof is a homogeneous convolution operator, it gains one derivative
in norms such as H\"older, Zygmund, Sobolev, but not in the sup-norm
or the $\cC^l$-norm.

\noi 3.\ Theorem 3.1 has the following extension to non-closed
\tr\ $\cC^1$-submanifolds $M'$ in $\C^n$. Let $K$ be a compact subset
of $M'$ and let $K' \subset M'$ be a compact \nbd\ of $K$ in
$M'$. For $\d>0$ we set
$$ U_\d =\{z\in \C^n \colon d_K(z)<\d\},\qquad
   U'_\d=\{z\in \C^n\colon d_{K'}(z)<\d\}.
$$
Choose $c\in (0,1)$. Given a form $u\in \cC^l_{p,q}(U'_\d)$
with $\dibar u=0$, we can solve $\dibar v=u$ in $U_{c\d}$,
and the estimates in theorem 3.1 remain valid when the tube
${\cal T}_{c\d}$ is replaced by $U_{c\d}$ on the left hand side,
and ${\cal T}_\d$ is replaced by $U'_\d$ on the right hand side
of each estimate. The proof can be obtained by simple
modifications of the kernel construction below.
This applies to compact \tr\ submanifolds with boundary in $\C^n$
since any such is a compact domain in a larger \tr\ submanifold.
\endpr

In this section $C$ denotes some constant whose value
may change every time it occurs, but which does not depend
on quantities such as $u$, $\d$ etc.

For $(0,1)$-forms $u$ (and hence for $(p,1)$-forms, $0\le p\le n$)
a large part of the result comes from the interior elliptic
regularity of the $\dibar$-operator and has nothing to
do with the particular solution $v$:

%
%
%
%
\proclaim 3.2 Lemma: Let $0 < c < 1$. There exist constants
$C_l>0$ ($l\in \N$) satisfying the following. 
If $K \subset \cn$ is a compact subset,
${\cal T}_{\d}=\{z\colon d(z,K)<\d\}$,
and $v$ a continuous function in ${\cal T}_{\d}$
such that $\dibar v \in \cC^l_{(0,1)}({\cal T}_\d)$
for some $l\in \N$ then $v \in \cC^l({\cal T}_{\d })$ and
$$ ||\di^{\a}v||_{L^\infty ({\cal T}_{c\d})} \le
   C_{l} \bigl(\d ||\di^{\a}\dibar v||_{L^\infty ({\cal T}_{\d})} +
   \d^{-|\a|}||v||_{L^\infty ({\cal T}_{\d})} \bigr);\quad |\a|\le l.
$$
If $\dibar v = \dibar f$ for some $f \in \cC^{l+1}({\cal T}_{\d })$,
then $v$ is also of class $\cC^{l+1}$ and satisfies
$$ ||\di_j\di^{\a}v||_{L^\infty ({\cal T}_{c\d})} \le
   C_{l+1} \bigl( \omega (\di_j\di^{\a}f,\d) +
   \d^{-l-1}||v||_{L^\infty ({\cal T}_{\d})} \bigr); \quad |\a|=l.
$$

\demo Proof:
We apply the Bochner-Martinelli formula
$$ g(z) =  \int_{\di {\cal T}_{\d}} g(\zeta)B(\zeta,z) -
         \int_{{\cal T}_{\d}} \dibar g (\zeta)\wedge B(\zeta,z),
$$
valid for $g \in \cC^1(\bar {\cal T}_{\d})$, where $B(\zeta,z)$ is the
Bochner-Martinelli kernel
$$B(\zeta,z) = c_n \sum_{j=1}^n (-1)^{j-1}{\bar{\zeta _j - z_j}
 \over ||\zeta - z||^{2n}} d\bar{\zeta} [j] \wedge d\zeta
$$
which is a closed integrable $(n,n-1)$-form. Let $z \in {\cal T}_{c\d}$
and let $\chi\colon \R\to[0,1]$ be a cut-off function with
$\chi (t) = 1$ when $|t| \le {1 \over 2}$ and $\chi (t) = 0$ when
$|t| \ge 1$. For $w \in \cn$ set
$\chi _{\d}(w) = \chi ({2|w| \over (1-c)\d})$.
It follows that the partial derivatives satisfy
$|\di ^{\a}\chi _{\d}| \le C_{\a}\d ^{-|\a|}$
for all $\a$. Applying the Bochner-Martinelli formula to
$g(\zeta) = \chi _{\d}(\zeta -z)v(\zeta)$ we obtain
$$ \eqalignno{ v(z) &= -\int_{{\cal T}_{\d}} \dibar_{\zeta}
        \bigl(\chi _{\d} (\zeta -z) v(\zeta) \bigr) \wedge B(\zeta,z) \cr
    &= -\int_{{\cal T}_{\d}} \dibar v(\zeta)
         \wedge \chi_{\d}(\zeta -z)B(\zeta,z)
       -\int_{{\cal T}_{\d}} v(\zeta)\, \dibar_{\zeta}\chi_{\d}(\zeta -z)
        \wedge B(\zeta,z) \cr
    &= I_1(z) + I_2(z).  & (3.5) \cr}
$$
These are convolution operators and we may differentiate on either
integrand. This gives for $|\a|\le l$ :
$$ \eqalignno{ \di^{\a}v(z) &= \di ^{\a}I_1(z) + \di ^{\a}I_2(z) \cr
               &= -\int_{{\cal T}_{\d}} \di^{\a} \dibar v(\zeta)\wedge
                  \chi_{\d}(\zeta -z) B(\zeta,z)
                 -\int_{{\cal T}_{\d}} v(\zeta)\, \di_z^{\a}
     \bigl( \dibar \chi _{\d}(\zeta -z)\wedge B(\zeta,z) \bigr). \cr}
$$
Setting $c^{\prime} = {1-c \over 2}$, $c'' = {1 \over 2}c'$ and using
$|B(\zeta,z)| \le C|\zeta -z|^{1-2n}$ we can estimate the integrals
for $|\a|\le l$ as follows:
$$ \eqalignno{ |\di ^{\a}I_1(z)|
          &\le C\int_{|\zeta -z| \le c^{\prime }\d }
          |\di ^{\a}\dibar v(\zeta)|\cdotp |\zeta -z|^{1-2n} dV \cr
  & \le C||\di ^{\a}\dibar v||_{L^\infty ({\cal T}_{\d})}
    \int_{|\zeta -z| \le c^{\prime }\d }|\zeta -z|^{1-2n}dV \cr
  & \le C ||\di ^{\a}\dibar v||_{L^\infty ({\cal T}_{\d})}
          \int_0^{c^{\prime }\d}{r^{2n-1}dr \over r^{2n-1}}
          \le C \d ||\di ^{\a}\dibar v||_{L^\infty ({\cal T}_{\d})}, \cr}
$$
$$ |\di ^{\a}I_2(z)| \le C \int_{c^{\prime \prime }\d \le |\zeta -z|
    \le c^{\prime }\d} |v(\zeta)|\cdotp \d^ {-2n-|\a |} dV
    \le C ||v||_{L^\infty ({\cal T}_{\d})}\d^{-|\a |}.
$$
This proves the first estimate in lemma 3.2. The estimate
for $|\di^\a I_2|$ also holds for derivatives of order
$|\a|=l+1$.

We now assume that $\dibar v = \dibar f$ for some
$f\in \cC^{l+1}({\cal T}_\d)$; then $v-f$ is holomorphic
and hence $v$ is also $\cC^{l+1}$. We wish to estimate the
derivatives of order $l+1$ of $I_1(z)$. For $|\a|=l$ we have
$$ \di_j\di^{\a}I_1(z) =
   - \di_j \int_{{\cal T}_{\d}}
    \di^{\a}\dibar v(\zeta)\wedge \chi_{\d}(\zeta -z)B(\zeta,z).
$$
We now apply (3.5) to $f$, replacing $\dibar f$ by
$\dibar v(=\dibar f)$ in the first term on the right hand side
and differentiating under the integral, to get
$$\eqalign{ \di_j\di^{\a}f(z) &=
        -\di_j\int_{{\cal T}_{\d}} \di^{\a}\dibar v(\zeta)\wedge
                        \chi_{\d}(\zeta -z) B(\zeta,z) \cr
       &\qquad - \int_{{\cal T}_{\d}}\di_j\di^{\a} f(\zeta)\wedge
       \dibar \chi_{\d}(\zeta -z)\wedge B(\zeta,z). \cr}
$$
Observe that the first term on the right hand side
equals $\di_j\di^{\a}I_1(z)$ from the previous display.
For a fixed $z\in \C^n$ we also apply (3.5) to the constant
function $\di_j\di^\a f(z)$:
$$ \di_j\di^{\a}f(z) = - \int_{{\cal T}_{\d}}\di_j\di^{\a}f(z)\,
    \dibar \chi _{\d}(\zeta -z)\wedge B(\zeta,z).
$$
Combining the above three formulas we get
$$ \di_j\di^{\a}I_1(z) = \int_{{\cal T}_{\d}}
   \bigl( \di_j\di^{\a} f(\zeta )- \di_j\di^{\a}f(z) \bigr)\,
   \dibar \chi _{\d}(\zeta -z)\wedge B(\zeta,z)
$$
and hence $|\di_j\di^{\a}I_1(z)| \le C\omega (\di_j\di^\a f,\d)$.
\endpr

From lemma 3.2 it follows that the estimates (3.1) and
(3.2) in theorem 3.1 will be proved for $(p,1)$-forms $u$
if we can find a solution $v$ which satisfies a sup-norm estimate
$ ||v||_{L^\infty ({\cal T}_{c\d})} \le
   C \d||u||_{L^\infty ({\cal T}_{\d})}$.
Such a solution is obtained by a linear operator given
by an integral kernel that we now construct.

\demo Construction of the kernel for $(0,1)$-forms:
We shall use the Koppelman's formula which we now recall.
For $v,w \in \cn$ let $<v,w> = \sum_{i=1}^{n} v_i w_i$.
Let $V\subset \C^{2n}$ and $\Omega'\subset \Omega\ss\C^n$ be
open subsets such that $\Omega$ has piecewise $\cC^1$-boundary
and $\bar \Omega \times \Omega \subset V$. Let
$P= P(\zeta,z)=(P_1,\ldots,P_n) \colon V\to \C^n$ be a
$\cC^1$-map satisfying

\item{(i)}  $P(\zeta,z) = \bar{\zeta -z}$ in a neighborhood
of the diagonal of $\Omega^{\prime} \times \Omega^{\prime}$, and

\item{(ii)} the function $\Phi\colon V\to\C$,
$\Phi(\zeta,z) = <P(\zeta,z), \zeta - z>$,
satisfies $\Phi(\zeta,z) \ne 0$ when $z \in \Omega^{\prime}$
and $\zeta \in \bar\Omega\bs \{z\}$.

\medskip Any such map $P$ is called a {\it Leray map}
for the pair $\Omega'\subset \Omega$ and $\Phi$ is the
corresponding {\it support function}.
We shall use the notation
$$ \eqalign{ d\zeta &= d\zeta_1\wedge\cdots\wedge d\zeta_n, \cr
   \dibar_{\zeta}P[j] &= \dibar_{\zeta} P_1\wedge\cdots
   \widehat{\dibar_{\zeta}P_j}\cdots \wedge \dibar_{\zeta}P_n, \cr
   \dibar_{\zeta}P[i,j] &= \dibar_{\zeta}P_1\wedge\cdots
   \widehat{\dibar_{\zeta}P_i}\cdots
   \widehat{\dibar_{\zeta}P_j}\cdots\wedge \dibar_{\zeta}P_n. \cr}
$$
Define the integral kernels
$$ \eqalign{ K(\zeta,z) &= c_n  \Phi(\zeta,z)^{-n}
                 \sum_{j=1}^n (-1)^{j-1} P_j\,
                 \dibar_{\zeta}P[j]\wedge d\zeta, \cr
             L(\zeta,z) &= c_n  \Phi(\zeta,z)^{-n}
                 \sum_{i \ne j} (-1)^{i+j} P_j\, \dibar_z P_i \wedge
                 \dibar_{\zeta}P[i,j]\wedge d\zeta. \cr}
$$
Note that the kernel $K(\zeta,z)$ is locally integrable when
$z\in \Omega'$. It is also important to observe that,
if $a(\zeta,z)$ is a $\cC^1$ function, the kernels generated
by $P$ resp.\ $aP$ are identical outside the zero set $a=0$.
For a suitable choice of the constant $c_n\in\R$ we then
have the following Koppelman-Leray representation formula
for $\dibar$-closed $(0,1)$-forms $u \in \cC^1_{0,1} (\bar\Omega)$:
$$ u(z)= \int_{\di \Omega} L(\zeta,z) \wedge u(\zeta)
         + \dibar_z \int_{\Omega} K(\zeta,z) \wedge u(\zeta),
         \qquad z\in \Omega'.                  \eqno(3.6)
$$
This follows by applying the Stokes formula to the first integral
on the right hand side to transfer the integration to an
$\e$-sphere around $z$ and using $\dibar_z K= -\dibar_\zeta L$;
in the limit as $\e\to 0$ we obtain (3.6) by a usual residue
calculation. For $\zeta$ near $z$, the kernel $L$ coincides with
the B-M kernel for $(0,1)$-forms; in fact, for the Leray map
$P(\zeta,z)=\bar{\zeta-z}$, (3.6) is the classical
Bochner-Martinelli-Koppelman formula.

We have a lot of freedom in the choice of the map $P$ which determines
$K$ and $L$. If we choose it such that $P(\zeta,\cdotp)$ is
holomorphic in $\Omega'$ when $\zeta \in \di \Omega$, then $L(\zeta,z)=0$
for such $\zeta$ and $z$ (since each term in $L$ contains a derivative
$\dibar_z P_i$), and hence the function
$$    v(z) = \int_{\Omega} K(\zeta,z) \wedge u(\zeta) \eqno(3.7) $$
solves the equation $\dibar v = u$ in $\Omega^{\prime}$.

We shall construct the integral kernel of our solution operator
on ${\cal T}_\d$ by combining the Bochner-Martinelli kernel
near the diagonal $\zeta=z$ of the smaller tube
${\cal T}_{c\d}$ with the Henkin kernel when $\z$ is near the
boundary of ${\cal T}_\d$ and $z\in {\cal T}_{c\d}$.
This will give a family of linear solution operators
of the form (3.7) depending on $\d$ for small $\d>0$.

Let $\rho$ be the \spsh\ function mentioned in the beginning of
this section. Since
$\{\rho  < (1-\e)\d^2 \} \subset {\cal T}_\d \subset \{\rho < (1+\e)\d^2 \}$
for sufficiently small $\d>0$, we may replace the tube ${\cal T}_\d$
with the sublevel sets $\{\rho < \d^2\}$ which we still denote
by ${\cal T}_\d$.

The construction of the kernel will proceed through several lemmas.
First we recall from [HL] the following well known result about
the existence of the Henkin support function $\Phi$
and the corresponding Leray map $P$ on a fixed \spsc\
domain which in our case is a tube ${\cal T}_{\d_0}$.

%
%
\medskip\noi \bf 3.3 Lemma. \sl
There exist constants $C , R > 0$ such that, for $\d_0 > 0$
sufficiently small, there are functions $\Phi(\zeta,z)$ and
$A(\zeta,z)$ in $\cC^1({\cal T}_{\d_0}\times {\cal T}_{\d_0})$,
with $\Phi$ holomorphic in $z$, and there is a $\cC^1$-function
$B(\zeta,z)$, defined for $\zeta,z \in {\cal T}_{\d_0}$
and $|\zeta-z| \le R$, satisfying the following:

\item{(i)}   $\Phi(\zeta,z) = A(\zeta,z)B(\zeta,z)$,

\item{(ii)}  $|B(\zeta,z)| \ge C$ and
${\rm Re} A(\zeta,z) \ge \rho(\zeta)-\rho(z)+C|\zeta-z|^2$ when
$|\zeta-z| \le R$,

\item{(iii)} $|\Phi(\zeta,z)| \ge C$ when
$|\zeta-z| \ge {R\over 2}$, and

\item{(iv)}  with $\Phi$ as above, there exists a map
$P=P(\zeta,z) = (P_1,\cdots,P_n)$ such that for all $j$,
$P_j  \in \cC^1({\cal T}_{\d_0}\times {\cal T}_{\d_0})$,
$P_j$ is holomorphic in $z$, and
$\Phi(\zeta,z) = <P(\zeta,z),\zeta - z>$.
\medskip\rm

\demo Proof: This follows from the proof of Theorems 2.4.3 and
2.5.5.\ in [HL]. $A(\zeta,z)$ is an approximate Levi polynomial
in $z\in \C^n$ of the form
$$ A(\zeta,z) = 2\sum_{j=1}^n
              {\di \rho (\zeta) \over \di \zeta_j}(\zeta_j - z_j)
             -\sum_{j,k=1}^n a_{jk}(\zeta)
              (\zeta_j - z_j)(\zeta_k -z_k),
$$
where $a_{jk}$ are $\cC^1$ functions which approximate the partial
derivatives $\di^2\rho/\di\z_j\di \z_k$ sufficiently well
on ${\cal T}_{\d_0}$ (Lemma 2.4.2 in [HL]). In fact, when $\rho$ is
of class $\cC^3$ or better, we might simply take
$a_{jk}=\di^2\rho/\di\z_j\di \z_k$.

The only small change from [HL] is that, in our situation, the
maps $\Phi$ and $P$ may be defined globally for
$\zeta\in {\cal T}_{\d_0}$, and not only for $\zeta$ near the
boundary of ${\cal T}_{\d_0}$, provided that $\d_0>0$ is sufficiently
small. This follows from the thinness of the tube ${\cal T}_{\d_0}$
and can be seen as follows. Observe that for $\z\in M$
the linear terms in $A(\zeta,\cdotp)$ vanish and we have
$\Re A(\zeta,z)<0$ for all points $z\in M\bs \{\z\}$
sufficiently close to $\z$. Hence for $\d_0>0$ small we
can choose $\e>0$ (depending on $\d_0$) such that
$\Re A(\zeta,z)<0$ whenever $z, \zeta \in {\cal T}_{\d_0}$ and
$\e \le |\z-z| \le 2\e$. The proof of Theorem 2.4.3 in [HL]
(which proceeds by cutting of $\log A$ on
$B(\z,2\e) \cap {\cal T}_{\d_0}$ and solving a $\dibar$-equation
on ${\cal T}_{\d_0}$) then gives a globally defined $\Phi$
(and hence $P$).
\endpr

Let $\Phi$, $P$, $A$ and $B$ be as in lemma 3.3, constructed
on a fixed tube ${\cal T}_{\d_0}$. $P$ is not quite
a Leray map since it does not equal $\bar{\zeta-z}$ near the diagonal,
and we shall now modify it suitably on tubes ${\cal T}_{\d}$
for $0< \d \le \d_0$. Let $0< c<c^{\prime}<1$.
Choose a cut-off function $\l_{\d}$ such that $\l_{\d} = 1$
in ${\cal T}_{c^{\prime}\d}$ and $\l_{\d} = 0$ near
$\di {\cal T}_{\d}$. We may assume that its (real) gradient satisfies
$||\nabla \l _{\d}|| \le C\d^{-1}$ for some $C>0$
independent of $\d$. We will show that for a
suitably chosen function $\phi (\zeta,z)$ on
$\bar {\cal T}_\d\times \bar {\cal T}_\d$,
the conditions in Koppelman's formula (3.6) are satisfied for the pair
of domains $\Omega = {\cal T}_{\d}$ and $\Omega^{\prime} = {\cal T}_{c\d}$
if we define the Leray map $\wt{P}$ by
$$\wt{P}(\zeta,z) = (1-\l_{\d}(\zeta))\phi (\zeta,z)
        P(\zeta,z) + \l_{\d}(\zeta)\bar{\zeta -z},
$$
with the corresponding support function given by
$$ \wt{\Phi} = <\wt{P},\zeta -z> =
     (1-\l_{\d})\phi \Phi + \l_{\d}|\zeta -z|^2.
$$
We need to find $\phi$ such that $\wt{\Phi} (\zeta,z) \ne 0$
when $z \in {\cal T}_{c\d}$ and $\zeta \in \bar {\cal T}_\d \bs\{z\}$.
When $\zeta\in \bar {\cal T}_{c\d}$, we have $\wt{\Phi}(\z,z)=|\z-z|^2$,
so the condition is satisfied for any choice of $\phi$.
Hence it suffices to consider the points $\zeta$ where
$\rho (\zeta) > \rho (z)$. Let $\psi\colon \R\to [0,1]$ be a
cut-off function such that $\psi (t) = 1$ for $|t| \le {1\over 2}R$
and $\psi (t) = 0$ for $|t| \ge {2\over 3}R$. Set
$$\phi(\zeta,z) = \psi (|\zeta - z|) B(\zeta,z)^{-1} +
                 \bigl(1-\psi (|\zeta - z|)\bigr)
                 \bar{\Phi(\zeta,z)},
$$
where $B$ is as in lemma 3.3. Then
$\phi \Phi = \psi A + (1-\psi)|\Phi |^2$
(since $B^{-1}\Phi=A$), and we have the following
estimates for the real part
$\theta(\zeta,z) := {\rm Re} \phi \Phi (\zeta,z)$
when $\rho(\z)>\rho(z)$:

\item{--} when $|\zeta - z| \le {1\over 2}R$,
$\theta ={\rm Re} A \ge C|\zeta -z|^2$,
\item{--} when ${1\over 2}R \le |\zeta - z| \le R$,
$\theta= \psi {\rm Re} A + (1-\psi )|\Phi |^2 \ge
   \psi C|\zeta - z|^2 + (1-\psi )C^2 > 0$,
\item{--} when $|\zeta - z| > R$, $\theta= |\Phi |^2 \ge C^2$ .

\medskip\noi
This verifies the required properties, and hence (3.6) is valid
when $K(\zeta ,z)$ and $L(\zeta ,z)$ are the kernels generated by the
Leray map $\wt{P}$. For $\zeta $ near $\di {\cal T}_{\d}$ we have
$\wt{P} = \phi P$; since $\phi \ne 0$ there, the kernel $L$ is
identical to the one generated by the holomorphic Leray map $P$,
and hence the first term in (3.6) is zero. This gives
us the solution formula (3.7) for the equation $\dibar v = u$
in ${\cal T}_{c\d}$. This completes the construction of the kernel
for $(0,1)$-forms.

\demo Proof of the sup-norm estimates:
It suffices to show that the sup-norm estimate holds in our situation
when $n\ge 3$. In case $n<3$ we simply identify $\C^n$ with
$\C^n \times \{0\} \subset \cthree$ and extend $f$ independently
of the additional variables; the solution to the extended problem will
satisfy the estimates, and its restriction to $\C^n$ will be a
solution to the original $\dibar$-problem.

\proclaim 3.4 Lemma: The solution
$v(z) = \int_{{\cal T}_{\d}} K(\zeta,z) \wedge u(\zeta) $
defined above satisfies the estimate
$ ||v||_{L^\infty ({\cal T}_{c\d})} \le C \d||u||_{L^\infty ({\cal T}_{\d})}$
when $n\ge 3$.

\demo Proof: Let $P^0 = \phi P$. $P^0$ is independent of $\d $
and $\wt{P} = (1-\l )P^0 + \l (\bar{\zeta -z})$.
This gives
$$ \eqalignno{ \dibar \wt{P}_j
       & = \dibar \l (\bar{\zeta _j - z_j} - P^0_j) +
          (1-\l )\dibar P^0_j + \l d\bar{\zeta}_j \cr
       & =: \dibar \l (\bar{\zeta _j - z_j} - P^0_j) + \eta _j & (3.8)\cr}
$$
The terms in $K(\zeta ,z)$ are of the form
$\wt{\Phi} ^{-n}\wt{P}_j \dibar \wt{P} [j]\wedge d\zeta$.
Since $\dibar \l \wedge \dibar \l = 0$, this is a sum of terms of the
following two types:
$$   \wt{\Phi} ^{-n}\wt{P}_j \eta [j]\wedge d\zeta \quad {\rm and}\ \
     \wt{\Phi} ^{-n}\wt{P}_j (\bar{\zeta _k - z_k} -
         P^0_k)\dibar \l \wedge \eta [j,k]\wedge d\zeta.
$$
We shall estimate the integrals of these over ${\cal T}_{\d}$ when
$z \in {\cal T}_{c\d}$. We have already shown that
$Re\wt{\Phi}(\zeta ,z) \ge C|\zeta -z|^2$.
For $|\zeta - z| \le {1\over 2}R$ we have
$$ \eqalignno{ \sum_{j=1}^n P_j^0(\zeta ,z)(\zeta _j - z_j)
       & = <P^0(\zeta ,z),\zeta -z> = A(\zeta ,z) \cr
       & = 2\sum_{j=1}^n {\di \rho (\zeta) \over \di \zeta_j}
          (\zeta_j - z_j) + \cO (|\zeta -z|^2). \cr}
$$
This implies
$P_j^0(\zeta ,z) = 2{\di \rho (\zeta) \over \di \zeta_j} +
                   \cO (|\zeta -z|) = \cO (\d + |\zeta -z|)
$.
By choice of $\l$ this gives
$(1-\l(\zeta ,z))P_j^0(\zeta ,z) = \cO (|\zeta -z|)$
and therefore $\wt{P}_j (\zeta ,z) = \cO (|\zeta -z|)$.
Since $|\eta _j| \le C$, we get
$$ \eqalign{ \wt{\Phi} ^{-n}\wt{P}_j \eta [j]\wedge d\zeta &=
   \cO (|\zeta -z|^{1-2n}), \cr
   \wt{\Phi} ^{-n}\wt{P}_j(\bar{\zeta _k - z_k} -P^0_k)\dibar \l
           \wedge \eta [j,k]\wedge d\zeta
       &= \cO (|\zeta -z|^{1-2n} + \d ^{-1}|\zeta -z|^{2-2n}), \cr}
$$
which shows that the kernel $K(\zeta ,z)$ has a singularity of the
same type as the Bochner-Martinelli kernel on the diagonal.

Locally we may straighten $M$, i.e., for each $p \in M$ there is a
neighborhood $V_p$ of $p$ and a $\cC^1$-diffeomorphism
$\Psi \colon U \rightarrow V_p$, where $U$ is a neighborhood of
the origin in ${\bf R}^{2n}$, such that $\Psi$ is nearly volume
and distance preserving, and
$\Psi (U \cap {\bf R}^m) =  V_p \cap M$, where $m$ is the dimension
of $M$. We denote the points in ${\bf R}^{2n}$ by
$(u^{\prime },u^{\prime \prime }) \in {\bf R}^m \times {\bf R}^{2n-m}$.
By compactness we may assume that ${\cal T}_{\d}$ is covered by a finite number
(independent of $\d$) of sets
$$K_j^{\d } = \Psi_j (\{ (u^{\prime },u^{\prime \prime });
       \ |u^{\prime }| \le a, |u^{\prime \prime }| \le \d \})
$$
for some constant $a$. We keep the notation $\zeta$ and $z$ for
the points in the new coordinates also.
We then have the estimate ($\zeta = (u^{\prime },u^{\prime \prime }))$:
$$ \eqalignno{ \big| \int_{K_j} K(\zeta ,z)\wedge u(\zeta ) \big|
       & \le C||u||_{L^\infty ({\cal T}_{\d})} \int_{|u^{\prime }|
         \le a, |u^{\prime \prime }|\le \d}
         (|\zeta -z|^{1-2n} + \d ^{-1}|\zeta -z|^{2-2n}) dV(\z) \cr
       & \le C||u||_{L^\infty ({\cal T}_{\d})} \int_{|u^{\prime }|
         \le a, |u^{\prime \prime }|\le \d}
         (|\zeta |^{1-2n} + \d ^{-1}|\zeta |^{2-2n}) dV(\z). \cr}
$$
For $m < t < 2n$ we estimate these integrals as follows:
$$ \int_{|u^{\prime }| \le a, |u^{\prime \prime }|\le \d}
                              {1 \over |\zeta |^t}
     \le  C\left( \int_0^{\sqrt2 \d} {r^{2n-1}dr \over r^t} +
          \int_{\d }^a {\d ^{2n-m}r^{m-1}dr \over r^t} \right)
     \le  C\d ^{2n-t}.                          \eqno(3.9)
$$
Hence
$$  \big| \int_{K_j} K(\zeta ,z)\wedge u(\zeta ) \big|
        \le C||u||_{L^\infty ({\cal T}_{\d})}
        \bigl(\d + \d^{-1}\d^2 \bigr)
        = 2C ||u||_{L^\infty ({\cal T}_{\d})}
$$
when $2n-2 > m$. Since $m\le n$, this holds for $n > 2$.
\endpr

\demo Construction of the kernel for forms of higher degree:
We consider the form
$$
 K(\zeta,z)= c_n\wt{\Phi}(\zeta,z)^{-n}
            \sum_{j=1}^n(-1)^{j-1}\wt{P}_j\dibar\wt{P}[j]
            \wedge d(\zeta -z)
$$
on ${\cal T}_\d \times {\cal T}_{c^\prime\d}$,
where $\dibar$ is now taken with
respect to both $\zeta$ and $z$. We decompose
$$ K(\zeta,z)=\sum_{p\leq n}\sum_{q\leq n-1} K_{p,q}(\zeta,z),$$
where $K_{p,q}$ has bidegree $(p,q)$ with respect to $z$ and
$(n-p,n-q-1)$ with respect to $\zeta$. When $q>0$, $K_{p,q}(\zeta,z)=0$
when $z\in {\cal T}_{c^\pr\d}$ and $\zeta$ is near $\di {\cal T}_\d$.
(Recall that $K(\zeta,z)=K^\d(\zeta,z)$ depends on $\d$ via
the cutoff function $\l_\d$.)
It follows that the $(p,q-1)$-form
$$
    v(z)=\int_{{\cal T}_\d} K_{p,q-1}(\zeta,z) \wedge
    u(\zeta)=(-1)^{p+q}\int_{{\cal T}_\d}
    u(\zeta) \wedge K_{p,q-1}(\zeta,z)
$$
solves $\dibar v = u$ in ${\cal T}_{c'\d}$ for each $\dibar$-closed
$(p,q)$-form $u$ in ${\cal T}_\d$, $q>0$. The precise meaning of the
integral is as follows. Write
$$
   K_{p,q-1}(\zeta,z)=\sum_{|I|=p}
   \sum_{|J|=q-1}k_{I,J}(\zeta,z)dz^I\wedge d\bar{z}^J,
$$
where $k_{I,J}(\zeta,z)$ is an $(n-p,n-q)$-form in $\zeta\in {\cal T}_\d$
depending smoothly on $z\in {\cal T}_{c'\d}$. Then
$$
   v(z)=\sum_{|I|=p}\sum_{|J|=q-1}(-1)^{p+q}
          \left( \int_{{\cal T}_\d}u(\zeta)\wedge k_{I,J}(\zeta,z)
          \right)
          dz^I\wedge d\bar{z}^J.
$$
This completes the construction of the kernel. The reader may find
some additional references and historical remarks about the solution
formula at the end of this section.
\medskip

Before proceeding we make the following elementary

\demo Geometric observations:
Let $M$ be a compact $m$-dimensional $\cC^1$-submanifold of $\RN$.
There exists a constant $B>0$ such that, if $z_0,z_1 \in {\cal T}_\d(M)$ for
sufficiently small $\d$, then $z_0$ and $z_1$ may be joined by a path
in ${\cal T}_\d(M)$ of length no more than $B|z_1-z_0|$. This is due to the
fact that the tubes may be locally straightened, in a uniform way, to tubes
around $\Rm \times \{0\}$ in $\RN$.

From this we get the following:
If $u\in\cC^1({\cal T}_\d M)$, $\|u\|_{L^\infty({\cal T}_\d)}\leq A$,
$\|u\|_{\cC^1({\cal T}_\d)}\leq At^{-1}$ for $t\leq 1$ and $0<s<1$, then
$|u|_{s,\d}\leq {\rm max}(2,B)At^{-s}$. We see this as follows:
If $|h|\leq t$, we can integrate $Du$ from $z$ to $z+h$ to get
$|u(z+h)-u(z)||h|^{-s}\leq BAt^{-1}|h|^{1-s}\leq BAt^{-s}$.
If $|h|\geq t$, the triangle inequality gives
$|u(z+h)-u(z)||h|^{-s}\leq 2At^{-s}$.

We also have a corresponding result for compact manifolds $M$: if
$\|u\|_{\cC^r(M)}\leq A$ and $\|u\|_{\cC^{r+1}(M)}\leq At^{-1}$ for
$t\geq 0$, then $\|u\|_{\cC^{r+s}(M)}\leq CAt^{-s}$,
where $C$ is a constant independent of $u$.

\demo Proof of the estimates for forms of higher degree:
The proof of the sup-norm estimate, which we gave for
$(0,1)$-forms, carries over almost verbatim to the general case.
However, lemma 3.2 fails and we must proceed differently
to estimate the derivatives.

With $c_0=c'-c$ we introduce smooth cut-off functions
$\chi_\d\in\cC^\infty_0(B(0,c_0\d))$ with $\chi_\d(w)=1$
when $|w|<c_0\d/2$ and $|\di^\a\chi_\d|\leq C_\a\d^{-|\a|}$.
Then we decompose $v$ as $v^\pr + v^\dpr$, with
$$ \eqalign{
  v^\pr(z) &= \int_{{\cal T}_\d}\chi_\d(\zeta-z)
                K_{p,q-1}(\zeta,z)\wedge u(\zeta), \cr
  v^\dpr(z) &= \int_{{\cal T}_\d}(1-\chi_\d(\zeta-z))
                K_{p,q-1}(\zeta,z)\wedge u(\zeta), \cr}
$$
and estimate each summand separately.

Recall that when $z\in {\cal T}_{c\d}$ and $|\zeta-z|\leq c_0\d$, $K(\zeta,z)$
equals the Bochner-Martinelli kernel. Thus $v'(z)$ is obtained
for $z\in {\cal T}_{c\d}$ by applying a convolution operator to $u$;
hence
$$
     \di^\a v^\pr(z) =
     \int_{{\cal T}_\d}\chi_\d(\zeta-z)K_{p,q-1}(\zeta,z)\wedge\di^\a u(\zeta).
$$
Thus the components of $\di^\a v^\pr(z)$ are linear combinations of terms
$h(z)=(k*g)(z)$, where $k(w)=\chi_\d(w)\bar{w}_j|w|^{-2n}$ and $g$ is a component
of $\di^\a u$. Since $|k(w)|\leq|w|^{1-2n}$ and $k$ is supported
by $B(0,c_0\d)$, an obvious estimate gives
$|h(z)|\leq C\d\|g\|_\infty$, so
$$
     \|\di^a v^\pr\|_{L^\infty({\cal T}_{c\d})}\leq
     C\d \|\di^\a u\|_{L^\infty({\cal T}_\d)}.
$$
To estimate the finer norms of $h$ we introduce the auxiliary kernels
$$
    k_t(w)=\chi_\d(w)\bar{w}_j(t^2+|w|^2)^{-n};\qquad t>0.
$$
This is a smooth function of $(t,z)$ satisfying
$|k_t(z)|\le |k(z)|$ and $\lim_{t\to 0} k_t(z)=k(z)$.
Since each $k_t$ has compact support, it follows that
$\int\di_jD^\b k_t(w)dV(w)=0$ for every $(t,z)$-derivative $D^\b$.
Thus, setting $h_t(z)=(k_t*g)(z)$, we see that
$$
   D^\b\di_jh_t(z) = \int_{{\cal T}_\d} \di_jD^\b k_t(w)
   \bigl( g(z-w)-g(z) \bigr)\, dV(w).
$$
Observing that $|\di^\g\chi_\d(w)|\leq C_\g|w|^{-|\g|}$ on
${\rm supp}\chi_\d$, a simple calculation gives
$$
     |D^\b\di_jk_t(w)| \leq  C_\b|w|^{-|\b|}(t+|w|)^{-2n}.
$$
Assume that $\di^\a u\in\L^s({\cal T}_\d)$ for some $s\in(0,1)$.
We have $g\in\L^s({\cal T}_\d)$, and for $t>0$ we can estimate
in polar coordinates:
$$   \eqalign{
     |\di_jh_t(z)| &\leq
     C|g|_{s,\d}\int_{|w|<c_0\d}(t+|w|)^{s-2n}dV(w) \cr
     &\leq
     C|g|_{s,\d}\int_0^{c_0\d}r^{s-1}dr=C_s\d^s|g|_{s,\d}.\cr}
$$
For the first order derivatives with respect to $(t,z)$
we get in the same way:
$$   \eqalignno{ |D\di_jh_t(z)| & \leq C|g|_{s,\d}
             \int_{|w|<c_0\d}|w|^{s-1}(t+|w|)^{-2n}dV(w) \cr
     & \leq C|g|_{s,\d}\int_0^{c_0\d}(t+r)^{s-2}dr
       \leq C_s^\pr t^{s-1}|g|_{s,\d}.              & (3.10)\cr}
$$
By the dominated convergence theorem we have
$h_t(z)\to h(z)$ and
$$ \di_jh_t(z)\to h_{(j)}(z) = \int_{{\cal T}_\d} \di_jk(w)
   \bigl( g(z-w)-g(z) \bigr) \,dV(w)
$$
as $t\to 0$. We also have
$$
   |\di_jh_t(z)-h_{(j)}(z)|\leq\int_0^t|
   {\di\over\di\tau}\di_jh_\tau(z)|d\tau\leq
   C|g|_{s,\d}t^s;
$$
hence the convergence of the derivatives is uniform and
therefore $h_{(j)}(z)=\di_jh(z)$.
Thus $|\di_jh_t(z)|\leq C_s\d^s|g|_{s,\d}$, and we conclude that
$$
\|\di_j\di^\a v^\pr\|_{L^\infty({\cal T}_{c\d})}\leq C_s\d^s|\di^\a u|_{s,\d}.
$$
We have also shown that
$$
        \di_j\di^\a v^\pr(z) =
        \int_{{\cal T}_\d}\di_j \bigl( \chi_\d(\zeta-z)K_{p,q-1}(\zeta,z) \bigr)
        \wedge \bigl( \di^\a u(\zeta)-(\di^a u)_z \bigr).
$$
In order to estimate the $\L^s$-norm of $\di_j\di^\a v^\pr$
we need the following standard

\proclaim Lemma:  Let $\phi\in\cC({\cal T}_\d)$ have an extension
$\wt{\phi}\in\cC^1(\R^+\times {\cal T}_\d)$ satisfying
$|D\wt{\phi}(t,z)|\leq At^{s-1}$ for some $0<s<1$.
Then $\phi\in\L^s({\cal T}_\d)$ and $|\phi|_{s,\d}\leq (B+2/s)A$.

This is just a slight modification of Proposition 2 in the Appendix 1
of [HL]. Applying this to $\phi=\di_jh$ and $\wt{\phi}(t,z)=\di_jh_t(z)$,
(3.10) gives $|\di_jh|_{s,c\d}\leq C_s^\pr(B+2/s)|g|_{s,\d}$. Thus
$$
|\di_j\di^\a v^\pr|_{s,c\d}\leq C_s|\di^\a u|_{s,\d}.
$$

In order to study $v^\dpr$ we set
$K^\dpr_{p,q-1}(\zeta,z)= (1-\chi_\d(\zeta-z))K_{p,q-1}(\zeta,z)$
on ${\cal T}_\d\times {\cal T}_{c\d}$. This kernel has continuous $z$-derivatives
of all orders, and it equals zero when $|\zeta-z|\leq c_0\d/2$.
It follows that $v^\dpr$ is a smooth form with
$$
    \di^\a v^\dpr(z)=\int_{{\cal T}_\d}
    \di_z^\a K^\dpr_{p,q-1}(\zeta,z)\wedge u(\zeta).
$$
We recall the formula (3.8) and point out that $|\dibar\l_\d|=\cO(\d^{-1})$,
$|\wt{\phi}(\zeta,z)|\geq C|\zeta-z|^2$ on ${\cal T}_\d\times {\cal T}_{c\d}$,
and the quantities $|D_z\wt{\phi}(\zeta,z)|$, $|\eta_j|$ and
$|P_j^0|$ are all bounded by $C|\zeta-z|$, while their derivatives
with respect to $z$ are bounded independently of $\d$ (since $\l_\d$
is independent of $z$).

An induction on $|\a|$ shows that the components of
$\di_z^\a K_{p,q-1}^\dpr(\zeta,z)$ are linear combinations of terms
of the type
$$
\wt{\phi}^{-n-k}\di_z^\b(1-\chi_d(\zeta-z))a_0(\zeta,z)\cdots a_t(\zeta,z)
$$
with $\b\leq\a$, $k\leq|\a-\b|$ and $t\geq 2k+1-|\a-\b|$, and of
terms of the type
$$
    \wt{\phi}^{-n-k} {\di\l_\d(\zeta)\over\di\bar{\zeta}_i}
    a_0(\zeta,z)\cdots a_t(\zeta,z)
$$
with $k\leq|\a|$ and $t\geq 2k+2-|\a|$, where the $a_j(\zeta,z)$ have continuous
$z$-derivatives of all orders that have upper bounds independent of $\d$, and
$|a_j(\zeta,z)|\leq C|\zeta-z|$ when $t>0$ and $1\leq j\leq t$.
Since $|\zeta-z|>c_0\d/2$ when $K_{p,q-1}^\dpr\neq 0$, it follows easily that
$$
|\di_z^\a K_{p,q-1}^\dpr(\zeta,z)|\leq C_\a\d^{-1}|\zeta-z|^{2-2n-|\a|}.
$$
Thus
$$  \eqalign{ |\di^\a v^\dpr(z)| &\leq
        C_\a\d^{-1} \|u\|_{L^\infty({\cal T}_\d)}
        \int_{{\cal T}_\d\setminus B(z,c_0\d/2)}
        |\zeta-z|^{2-2n-|\a|}\, dV(\zeta)  \cr
    &\leq C_\a\d^{1-|\a|}\|u\|_{L^\infty({\cal T}_\d)} \cr}
$$
for $z\in {\cal T}_{c\d}$, $\a\in\Z^n_+$. The last estimate follows
for $|\a|>2$, $|\a|=2$ and $|\a|=1$, respectively, from the following
three integral estimates:
$$ \eqalignno{ &\int_{|\zeta-z|>\d} |\zeta-z|^{-2n-t}dV(\zeta) =
                                 C_t\d^{-t},\quad  t>0,    & (3.11)\cr
   &\int_{{\cal T}_\d\bs B(z,c_1\d)} |\zeta-z|^{-2n}dV(\zeta)
                        \leq C(c_1),\quad  z\in {\cal T}_\d,      & (3.12)\cr
   &\int_{{\cal T}_\d}|\zeta-z|^{s-2n}dV(\zeta)
                        \leq C_s\d^s;\quad  0<s<2n-m.      &(3.13)\cr}
$$
(3.11) follows immediately by a change of variable. (3.12) is proved exactly
like (3.9); in the sum in the middle of (3.9) the first integral has lower
limit $c_1\d$ instead of $0$. Finally (3.13) follows by setting
$t=2n-s$ in (3.9).

Using the geometric observation following the construction
of the kernel we have
$$
      \|\di^\a v^\dpr\|_{\L^s({\cal T}_{c\d})}
      \leq C_{\a,s}\d^{1-|\a|-s}\|u\|_{L^\infty({\cal T}_{c\d})}.
$$
This completes the proof of the H\"older estimates in theorem 3.1
for the case $0<s<1$. The proof for $s=1$ goes along the same lines,
with certain small modifications; since that case will not be used in
the paper, we omit the details.
\endpr

\demo Remarks on constructions of kernels:
The first integral kernel operators with holomorphic kernels,
solving the $\dibar$-equation on \spsc\ domains in $\C^n$,
have been constructed by  Henkin and, independently,
R.\ de Arellano (see the references in [HL]).
Henkin's approach is to patch the Bochner-Martinelli
and Leray kernels on the boundary $\di\Omega$. Our patching of
the two kernels (by first multiplying by $\phi$)
is the same as in \O vrelid [\O1, \O2]. The whole
construction is similar to the one by Harvey and Wells [HaW].

It seems that the first really precise $L^\infty$ and
$\cC^k$-estimates for the $\dibar$-equation in thin
tubes around a \tr\ submanifold $M\subset \C^n$, proved
by means of integral solution operators, are due to
Harvey and Wells [HaW] in 1972. A little later Range and Siu
[RS] (1974) used a more refined kernel construction to prove
estimates for the highest order derivatives of their solution on $M$
and deduced $\cC^k$-approximation of $\cC^k$-functions on a
$\cC^k$-submanifold $M \subset\C^n$ by \holo\ functions, a case
left open in [HaW]. In fact this approximation problem has been
one of the original motivations in proving such estimates.
Later on this approximation has been accomplished more efficiently,
and in greater generality, by Baouendi and Treves [BT1, BT2] by
using the convolution with the complex Gaussian kernel. This
latter method does not seem to give the approximation of diffeomorphisms
obtained in this paper because we must work in tubular
\nbd s and not solely on the submanifold.

As said earlier our construction of the kernel in this paper
is close to [HaW], and our main contribution is the way we
estimate the solutions. We find it quite striking that this
simple and seemingly crude construction of the kernel
gives rise to results that are essentially optimal
for the applications to mappings presented in this paper.
For the benefit of the reader we have given a fairly self
contained presentation based on the text [HL]. Another closely
related paper is [BB] where Bruna and Burgu\'es
approximate $\dibar$-closed jets on a \tr\ set $X$
in H\"older norms by functions \holo\ in a \nbd\ of $X$.
It seems likely that their method, making use of weighted
integral kernels of Anderson and Berndtsson type [AB],
may also be used to prove our results.
However, we believe that our approach is simpler and more
elementary. Our results, suitably reformulated, may also
be proved for \nbd s of \tr\ sets.

\beginsection 4.  Proof of theorems 1.2 and 1.3.

\demo Proof of theorem 1.2: We consider first the case
$\dim M_0=\dim M_1=n$. Let $d(z)$ denote the Euclidean distance
of $z$ to $M_0$, and let ${\cal T}_\d$ (resp.\ ${\cal T}'_\d$) denote the open
tube of radius $\d$ around $M_0$ (resp.\ around $M_1$).
The $\cC^k$-diffeomorphism $f\colon M_0\to M_1$ can be extended to a
$\cC^k$-map on $\C^n$, still denoted $f$, which is $\dibar$-flat to
order $k$ at $M_0$:
$$ |\di^\a(\dibar f)(z)| = o\bigl( d(z)^{k-1-|\a|} \bigr);
        \qquad 0\le |\a| \le k-1.
$$
In particular, the derivative $Df(z)$ is a non-degenerate $\C$-linear
map at each point $z\in M_0$ (the complexification of
$df_z\colon T_z M_0\to T_{f(z)} M_1$), and hence $f$ is a
$\cC^k$ diffeomorphism in some \nbd\ of $M_0$ in $\C^n$.
The $(0,1)$-form $u=\dibar f$ of class $\cC^{k-1}$
satisfies $\dibar u=0$ and
$$
    ||\di^\a u||_{L^\infty ({\cal T}_{\d})} = o(\d^{k-1-|\a|});
                            \qquad 0\le |\a| \le k-1
$$
as $\d\to 0$. Applying theorem 3.1 (specifically the estimates
(3.1), with $l=k-1 \ge 0$ and a fixed constant $0<c<1$),
we get for each sufficiently small
$\d>0$ a solution $v_\d$ to $\dibar v_\d=u$ in
${\cal T}_\d$ satisfying the following estimates:
$$ \eqalign{ || \di^\a v_\d||_{L^\infty ({\cal T}_{c\d})} &\le
   C\left( \d ||\di^{\a}u||_{L^\infty ({\cal T}_{\d})} +
       \d^{1-|\a|} ||u||_{L^\infty ({\cal T}_{\d})} \right) \cr
       &\le C\left( \d o(\d^{k-1-|\a|}) + \d^{1-|\a|} o(\d^{k-1}) \right) \cr
       &= o(\d^{k-|\a|});\qquad\qquad |\a|\le k-1. \cr}
$$
Moreover, since $\dibar v=u$ has a solution of class $\cC^{l+1}=\cC^k$,
namely $f$, we can choose $v_\d$ which in addition satisfies
the estimates (3.2) for the derivatives of top order $k$:
$$ ||\di^\a v_\d||_{L^\infty ({\cal T}_{c\d})} \le C \left(
      \omega_k(f;\d)+ \d^{-k+1} ||\dibar f||_{L^\infty ({\cal T}_{c\d})} \right)
      = o(1); \qquad |\a|=k.
$$
Here $\omega_k(f;\d)$ denotes the modulus of continuity
of the $k$-th order derivatives of $f$. Set $F_\d=f-v_\d$ in ${\cal T}_\d$.
Then $\dibar F_\d =0$ and the estimates on $v_\d$ imply
$$ ||F_\d - f||_{\cC^r({\cal T}_{c\d})} = ||v_\d||_{\cC^r({\cal T}_{c\d})}
   = o(\d^{k-r});                     \qquad 0\le r\le k      $$
which gives the first estimate in (1.2).
It remains to prove that $F_\d$ is biholomorphic and
satisfies the inverse estimates in (1.2) for all sufficiently small
$\d>0$. To simplify the notation we replace $\d$ by $c\d/2$, so that
$F_\d$ is \holo\ in the tube ${\cal T}_{2\d}$ and it satisfies
$$
    ||F_\d - f||_{\cC^r({\cal T}_{2\d})} = o(\d^{k-r}); \qquad 0\le r\le k
                                                             \eqno(4.1)
$$
as $\d\to 0$. Since $f$ is a diffeomorphism near $M_0$, so is any sufficiently
close $\cC^1$ approximation of $f$; hence (4.1) with $r=1$ implies that for
$\d>0$ sufficiently small, say $0< \d\le \d_0\le 1$, the map $F_\d$ is
diffeomorphic (and hence biholomorphic) in ${\cal T}_{2\d}$. Decreasing $\d_0$
if necessary, there is a number $a>0$ such that
$$
    |f(z)-f(z')|\ge 2a|z-z'|; \qquad  z,z'\in {\cal T}_{\d_0}.
$$
Since $f(M_0)=M_1$, the above implies that $f({\cal T}_\d)$ contains the
tube ${\cal T}'_{2a\d}$.

Fix an $\e>0$. By (4.1), applied with $r=0$, we get a constant
$\d_1=\d_1(\e)$, with $0<\d_1\le \d_0$, such that
$||F_\d - f||_{L^\infty({\cal T}_{2\d})} < a\e\d^k$ for $0<\d\le \d_1$.
Fix a point $z\in {\cal T}_\d$ and let $w=f(z)$. For each $z'$ with
$|z'-z|=\e\d^k$ we have
$$
    \eqalign{ |F_\d(z')-w| &= |(F_\d(z')-f(z'))+ (f(z')-f(z))| \cr
        &\ge |f(z')-f(z)| - |F_\d(z')-f(z')| \cr
        &\ge 2a\e\d^k - a\e\d^k =a\e\d^k. \cr}
$$
This means that the image by $F_\d$ of the sphere
$S=\{z'\colon |z'-z|= \e\d^k\}$ is a hypersurface containing
the ball $B(w;a\e\d^k)= \{w'\colon |w'-w|<a\e\d^k\}$ in the
bounded component of its complement. By degree theory
the $F_\d$-image of the ball $B(z;\e\d^k)$ contains the ball
$B(w;a\e\d^k)$. Hence there is a point $\z \in B(z;\e\d^k)$
such that $F_\d(\z)=w=f(z)$, and we have
$|F^{-1}_\d (w) - f^{-1}(w)| = |\z-z| < \e\d^k$.
Since this applies to any point $w\in {\cal T}'_{2a\d}$,
we conclude that $F_\d({\cal T}_{2\d}) \supset {\cal T}'_{2a\d}$ and
$$
    ||F^{-1}_\d - f^{-1}||_{L^\infty({\cal T}'_{2a\d})} \le \e\d^k;
                               \qquad 0<\d\le \d_1(\e). \eqno(4.2)
$$
Since $\e>0$ was arbitrary, this gives the inverse estimate
in (1.2) for $r=0$.

We proceed to estimate the derivatives of the inverse maps.
Denote by $||A||$ the spectral norm of a linear map $A\in GL(\R,2n)$.
Note that $Df^{-1}(w)=Df(z)^{-1}$ where $w=f(z)$.
Fix a point $w\in {\cal T}'_{2a\d}$ and let $z=f^{-1}(w)$,
$z_\d=F_\d^{-1}(w)$ (these are points in ${\cal T}_{2\d}$).
By (4.2) we have $|z-z_\d|\le \e\d^k$.
Writing $A=Df(z)$, $B=DF_\d(z_\d)$, we get
$$ \eqalign{ ||DF_\d^{-1}(w)- Df^{-1}(w)|| &= ||A^{-1}-B^{-1}|| \cr
   &=    ||A^{-1}(B-A)B^{-1}|| \cr
   &\le  ||A^{-1}||\,\cdotp ||A-B|| \,\cdotp ||B^{-1}||.\cr}
$$
Since $f$ is a diffeomorphism and $F_\d$ is $\cC^1$-close to
$f$, the eigenvalues of $A$ and $B$ are uniformly bounded away
from zero, and this gives a uniform estimate on $||A^{-1}||$ and
$||B^{-1}||$ (independent of $\d$). The middle term is
$$ ||A-B||=||Df(z)-DF_\d(z_\d)|| \le
   ||Df(z)-Df(z_\d)|| + ||Df(z_\d)-DF_\d(z_\d)||.
$$
The second term on the right hand side is of size $o(\d^{k-1})$
according to (4.1). As $\d\to 0$, we have $z_\d\to z$, and
hence the first term on the right hand side goes to zero
(by continuity of $Df$).
Hence $\sup\{||DF_\d^{-1}(w)-Df^{-1}(w)|| \colon w\in {\cal T}'_{2a\d}\}$
goes to zero as $\d\to 0$. This completes the proof when $k=1$.
If $k>1$, we can further estimate
$
    ||Df(z)-Df(z_\d)|| \le C|z-z_\d| \le C\e\d^k,
$
where $C$ is an upper bound for the second derivatives
of $f$. This gives
$$
    \sup \{||DF_\d^{-1}(w)- Df^{-1}(w)||
    \colon w\in {\cal T}'_{2a\d} \} = o(\d^{k-1})
$$
as required by (1.2) for derivatives or order $r=1$.
To get the estimates (1.2) for the higher derivatives of
$F^{-1}_\d - f^{-1}$ we may apply the same method to
the tangent map, i.e., the induced map on tangent bundles
over the tubes which equals the derivative of the given map
on each tangent space. We leave out the details.
This proves theorem 1.2 when $\dim M_0=n$.

Suppose now that $m=\dim M_0< n$. We are assuming that there is
an isomorphism $\phi\colon \nu_0\to \nu_1$ of the complex normal
bundles $\nu_0\to M_0$ resp.\ $\nu_1\to M_1$ over $f$;
by approximation we may assume that $\phi$ is of class $\cC^{k-1}$.
For each $z\in M_0$ we have $T_z\C^n =T_z^\C M_0\oplus \nu_{0,z}$.
Let $l_z$ be the $\C$-linear map on $\C^n$ which is uniquely
defined by taking $l_z=df_z$ on $T_z^\C M_0$ and $l_z=\phi_z$
on $\nu_{0,z}$. Clearly $l_z\in GL(n,\C)$ for each $z\in M$.
Applying lemma 2.6 we obtain a $\cC^k$-extension $\wt{f}$
of $f$ which is $\dibar$-flat on $M_0$. Now the proof may
proceed exactly as before. This proves theorem 1.2.
\endpr

%
%
\demo Remarks: 1.\ If $f_0\colon M_0\to M_1$ is a
{\it \ra\ diffeomorphism} and if the complex normal bundles to $M_0$
resp.\ $M_1$ are isomorphic over $f$ then $f$ {\it extends} to a
biholomorphic map $F$ from \nbd\ of $M_0$ onto a \nbd\ of $M_1$.
We see this as follows. Let $\phi \colon \nu_0\to \nu_1$ be the continuous
isomorphism (over $f$) of the complex normal bundles to $M_0$ resp.\ $M_1$.
There exist complexifications $\wt{M}_i \subset \C^n$ of $M_i$
($i=0,1$) such that $f$ extends to a biholomorphic map
$\wt{f}\colon \wt{M}_0 \to \wt{M}_1$ and such that the complex normal
bundles $\nu_i \to M_i$ extend to holomorphic vector bundles
$\wt{\nu}_i \to \wt{M}_i$. We define a continuous
map $\psi \colon M_0 \to GL(n,\C)$ by $\psi(z) = d\wt{f}_z \oplus \phi_z$.
Since $M_0 \subset \wt{M}_0$ is totally real,
$\psi$ may be approximated by a holomorphic map
$\wt{\psi}\colon \wt{M}_0 \to GL(n,\C)$. We now define
$\wt{\phi} \colon \wt{\nu}_0\to \wt{M}_1 \times \C^n$ by
$\wt{\phi}_{z}(v)=(\wt{f}(z),\wt{\psi}_z(v))$. Clearly $\wt{\phi}$
is a holomorphic vector bundle isomorphism between $\wt{\nu}_0$ and a
holomorphic sub-bundle $\wt{\nu}_2 \subset \wt{M}_1 \times \C^n$ which is
an approximation of $\wt{\nu}_1$. In particular, $\wt{\phi}$ is a
biholomorphic map between neighborhoods $V_i$ of the zero sections
of $\wt{\nu}_i$. These neighborhoods map biholomorphically onto 
neighborhoods of $M_0$ resp.\ $M_1$ under the projection maps
(the Docquier--Grauert theorem). This gives the desired 
biholomorphic extension of $f$.

\medskip
\noi 2.\ If instead of theorem 3.1 we use H\"ormander's $L^2$-estimates
when solving $\dibar v_\d=u(=\dibar f)$ in ${\cal T}_\d$, the resulting
holomorphic maps $F_\d=f-v_\d$ can be shown to satisfy the
weaker estimate
$||F_\d|_{M_0}-f||_{\cC^r(M_0)} = o(\d^{k-r-l})$
for $0\le r\le l$, where $l$ is the smallest integer larger
than ${1\over 2}\dim M_0$. This approach had been used in [FL].

\medskip
\demo Proof of theorem 1.3:
The proof can be obtained by repeating the proof of theorem 1.1
in [FL] (or its more technical version, theorem 2.1 in [FL]),
except that one applies theorem 3.1 above whenever solving
a $\dibar$-equation. This gives the improved estimates in (1.3)
with no loss of derivatives. We leave out the details.
\endpr

%
%
\noi\it A correction to {\rm [FL]}. \rm
We take this opportunity to correct an error in the proof
of Lemma 4.1 in [FL]. The equation numbers below refer to that paper.
The lemma is correct as stated,  but the proof of the estimate
(4.5) is not correct. Using the notation of that proof,
we have the higher variational equations
$$ {\partial \over \partial t} D^p\phi_t(x) =
   DX_t(\phi_t(x))\circ D^p\phi_t(x) + H^p_X(t,x)
$$
for $p \leq k$, where $D^pf$ denotes the $p$-th order derivative of a
map $f \colon \Omega\subset\Rn\to\Rn$, so $D^pf \in L^p(\Rn,\Rn)$.
$H^p_X(t,x)$ is a sum of terms involving derivatives of the vector
field $X_t$ and derivatives of order less than $p$ of the flow $\phi_t$,
and $H^1_X=0$. We use the same notation for $Y_t^\e$ and its flow $\psi_t^\e$.

Choose unit vectors $v_1,\cdots,v_p \in \Rn$ and set
$y(t)=[D^p\phi_t(x)-D^p\psi_t^\e(x)](v_1,\cdots,v_p)$. It will be
sufficient to show that $\|y(t)\| = o(\e^{k-p})$, uniformly for
$0\leq t \leq t_0$,  $x \in K(\e)$ and unit vectors
$v_1,\cdots,v_p$. $y(t)$ satisfies the differential equation
$$ \eqalign{y^\prime(t) &= DY_t^\e(\psi_t^\e(x))\cdotp y(t)
  +(DX_t(\phi_t(x))-DY_t^\e(\psi_t^\e(x))\circ D^p\phi_t(x)(v_1,\cdots,v_p)\cr
  &\quad +(H^p_X(t,x)-H^p_{Y^\e}(t,x))(v_1,\cdots,v_p). \cr}
$$
This is a linear system $y^\prime=A(t)\circ y + b(t)$, $y \in \Rn$.
Suppose the matrix norms satisfy $\|A(t)\| \leq A$ and $\|b(t)\| \leq b$
for $t \in [0,t_0]$. The function $u(t)=\|y(t)\|$ is differentiable
outside the zeroes of $u$, with
$u^\prime(t)= y^\prime(t)\cdotp y(t)/ \|y(t)\| \leq \|y^\prime(t)\|$,
so $u^\prime(t)\leq Au(t)+b$ outside the zeroes of $u$. Since
$\phi_0=\psi_0^\e=Id$, we have $y(0)=0$. We shall first show that
$u(t)\leq {b\over A}(e^{At}-1)$ for $t \in [0,t_0]$. If $u(t)=0$,
there is nothing to prove. If not, let $t_1$ be the largest
zero of $u$ on $[0,t]$. Thus $u^\prime(s)\leq Au(s)+b$ for
$s \in (t_1,t]$. Setting $v(s)=u(s)e^{-As}$ we get
$v^\pr(s)\leq be^{-As}$ for $s \in (t_1,t]$. Integration from
$t_1$ to $t$ gives $v(s)\leq {b\over A}(e^{-At_1}-e^{-At})$.
Thus $u(t)\leq {b\over A}(e^{A(t-t_1)}-1)\leq {b\over A}(e^{At}-1)$.

In our situation, the matrix norm of $A(t)=DY_t^\e(\psi_t^\e(x))$
is bounded independently of $\e>0$, $x\in K(\e)$ and $t$, by
(4.4). It is therefore sufficient to prove that $b=o(\e^{k-p})$,
uniformly in $x$, $t$ and unit vectors $v_1,\cdots,v_n$. It is
shown in [FL] that the matrix norm
$\|DX_t(\phi_t(x))-DY_t^\e(\psi_t^\e(x)\|_{L^\infty(K(\e))}=o(\e^{k-1})$.
Since the flow $\phi_t(x)$ is of class $\cC^k$, it follows that
the matrix norm $\|D^p\phi_t(x)\|$ is uniformly bounded for
$x \in K(\e)$ and $t \in [0,t_0]$. Applying (4.4) and (4.5)
inductively as in [FL] we obtain
$\|H^p_X-H^p_{Y^\e}\|_{L^\infty(K(\e))}=o(\e^{k-p})$,
uniformly in $t$, which proves the claim.
\endpr

\beginsection 5. Solving the equation $dv=u$ for holomorphic forms in tubes.

Let $d$ denote the exterior derivative. In this section we solve
the equation $dv=u$ with sup-norm estimates for holomorphic forms
in tubes ${\cal T}_\d={\cal T}_\d M$ around \tr\ submanifolds $M\subset \C^n$.
We denote by $\L^s$ the H\"older spaces as in sect.\ 3 above.
We first state our main result for closed submanifolds; for an extension
to compact submanifolds with boundary see remark 3 following
theorem 5.1.

%
%
%
%
\noi\bf 5.1 Theorem. \sl Let $i\colon M\hra\C^n$ denote the
inclusion of a closed, $m$-dimensional, \tr\ submanifold
of class $\cC^2$ in $\C^n$. Let a positive constant
$c<1$ be given. Then there exist positive constants
$C$, $\d_0$ and $C_s$ for all $s\in (0,1)$ such that,
if $u$ is a $d$-closed \holo\ $p$-form in the tube
${\cal T}_\d={\cal T}_\d M$ for some $0<\d\le \d_0$ and $1\le p\le n$, then:
\item{(a)} If $p>m$, the equation $dv=u$ has a \holo\ solution $v$
in ${\cal T}_\d$ satisfying
$$
     ||v||_{L^\infty({\cal T}_{c\d})} \le
     C\d ||u||_{L^\infty({\cal T}_{\d})}.
                                                    \eqno(5.1)
$$
\item{(b)} If $p\le m$ and the form $i^*u$ is exact on $M$,
then for any solution of $dv_0=i^*u$ of class $\L^s(M)$ ($0<s<1$)
there is a \holo\ solution $v$ of $dv=u$ in ${\cal T}_\d$ satisfying
$$
     ||v||_{L^\infty({\cal T}_{c\d})} \le C_s\left(
      \d ||u||_{L^\infty({\cal T}_{\d})} + ||v_0||_{L^\infty(M)}
     + \d^s ||v_0||_{\L^s(M)} \right).               \eqno(5.2)
$$
\item{(c)} If $p\le m$ and $i^*u$ is exact on $M$,
there is a \holo\ solution of $dv=u$ with
$$
    ||v||_{L^\infty({\cal T}_{c\d})} \le
    C ||u||_{L^\infty({\cal T}_{\d})}.
                                                    \eqno(5.3)
$$
\rm

\medskip
\demo Remarks: 1. If $\Omega$ is a Stein manifold, the Rham cohomology
groups $H^p(\Omega;\C)$ can be calculated by holomorphic forms
in the following sense: Each closed form is cohomologous to a
closed holomorphic $p$-form, and if a holomorphic form $u$ is
exact (i.e., $u=dv_0$ for some, not necessarily holomorphic,
$(p-1)$-form $v_0$), then also $u=dv$ for a holomorphic $(p-1)$-form
$v$ on $\Omega$. (See [H\"o], Theorem 2.7.10.)

\noi 2. On a $\cC^2$-manifold $M$, $\cC^1$-forms and the
$d$-operator $d\colon \cC^1_{p-1}(M) \to \cC^0_{p}(M)$ are
intrinsically defined. By duality, the notion $dv=u$ (weakly)
is well defined on $M$. The condition in theorem 5.1 that
$i^*u$ be exact on $M$ need only hold in the weak sense.

\noi 3. Theorem 5.1 has an extension to non-closed
\tr\ $\cC^1$-submanifolds $M'$ in $\C^n$. Let $K$ be a compact subset
of $M'$ and let $K' \subset M'$ be a compact \nbd\ of $K$ in
$M'$. (For instance, $K=M$ may be a compact \tr\
submanifold with boundary in $\C^n$.) For $\d>0$ we set
$$ U_\d =\{z\in \C^n \colon d_K(z)<\d\},\qquad
   U'_\d=\{z\in \C^n\colon d_{K'}(z)<\d\}.
$$
Choose $c\in (0,1)$. Assume that $u$ is a $d$-closed \holo\
$p$-form in $U'_\d$, with $i^*u$ exact on $U'_\d \cap M'$
(where $i\colon M'\hra \C^n$ is the inclusion map).
Then there is a \holo\ solution of $dv=u$ in $U_{c\d}$ such that
the estimates (5.1), (5.2) and (5.3) are valid when
${\cal T}_{c\d}$ is replaced by $U_{c\d}$ and ${\cal T}_\d$ is replaced by
$U'_\d$.
\endpr

\demo Proof of theorem 5.1:
We give the details in the case when $M$ is closed
(compact and without boundary); for the non-closed case
see remark 4 following the proof.

Since $M$ is a strong deformation retraction of the tube ${\cal T}_\d$,
the equation $dv=u$ has a differentiable solution on ${\cal T}_\d$ under
the assumptions above. The strategy is to first find a good
differentiable solution $v_1$ and then successively get rid of
its $(p-q-1,q)$-components for $q>0$. The second part, lemma 5.2
below, follows the proof of Serre's theorem
(Theorem 2.7.10 in [H\"o]) which amounts to solving a $\dibar$-equation
at each step. We use the solution provided by theorem 3.1; it is here
that we need the sharp estimates (3.3) and (3.4) for the H\"older
norms.

\proclaim 5.2 Lemma: Let $0<c<c_1<1$. Let $u$ be a closed
holomorphic $p$-form on ${\cal T}_\d$ for $0<\d\le \d_0$ as in theorem 5.1.
Suppose that there exists a differentiable
$(p-1)$-form $v_1$ on ${\cal T}_{c_1\d}$ satisfying $dv_1=u$ and
$$ ||v_1||_{L^\infty({\cal T}_{c_1\d})} \le A_\d, \qquad
   ||v_1||_{\Lambda^s ({\cal T}_{c_1\d})} \le A_\d\, \d^{-s},    \eqno(5.4)
$$
where $A_\d$ depends on $\d$ and $u$. Then there exists a
holomorphic $(p-1)$-form $v$ in ${\cal T}_{c\d}$ satisfying
$dv=u$ and $||v||_{L^\infty({\cal T}_{c\d})} \le C_0A_\d$ for
$0<\d\le \d_0$, where $C_0$ is an absolute constant.

\demo Proof of lemma 5.2: Let $v_1 = \sum_{q\le q_0} v_{(q)}$
where $v_{(q)}$ is of bidegree $(p-1-q,q)$. By comparing the
terms of bidegree $(p-1-q_0,q_0+1)$ in the equation
$dv_1=\di v_1+\dibar v_1=u$ and taking into account that
$u$ is holomorphic, we see that $\dibar v_{(q_0)}=0$.
If $q_0>0$, we get by theorem 3.1 a form $w$ on ${\cal T}_\d$ solving
$\dibar w= v_{(q_0)}$ and satisfying the following estimates
for some fixed $c<c_2<1$ and for all $1\le j\le 2n$:
$$\eqalign{
   ||\di_j w||_{L^\infty({\cal T}_{c_2\d})} &\le
   C_1\left( ||v_{(q_0)}||_{L^\infty({\cal T}_{c_1\d})} +
             \d^{s} ||v_{(q_0)}||_{\L^s({\cal T}_{c_1\d})} \right)
   \le 2C_1 A_\d,  \cr
   ||\di_j w||_{\L^s({\cal T}_{c_2\d})} \,\, &\le
   C_1\left( ||v_{(q_0)}||_{\L^s({\cal T}_{c_1\d})} +
             \d^{-s} ||v_{(q_0)}||_{L^\infty({\cal T}_{c_1\d})} \right)
   \le 2C_1 A_\d\, \d^{-s}. \cr}
$$
Thus the form $v_2=v_1-w$ solves $dv_2=dv_1=u$, it only
has components of bidegree $(p-q-1,q)$ for $q<q_0$, and it
satisfies
$$ ||v_2||_{L^\infty({\cal T}_{c_2\d})}   \le C' A_\d, \qquad
   ||v_2||_{\Lambda^s ({\cal T}_{c_2\d})} \le C' A_\d\, \d^{-s}.
$$
Repeated use of this argument gives a holomorphic solution
of the equation $dv=u$ satisfying
$||v||_{L^\infty({\cal T}_{c\d})} \le C_0A_\d$.
\endpr

To prove theorem 5.1 it thus suffices to construct a good
differentiable solution satisfying lemma 5.2, with $A_\d$ as small
as possible. Let
$F=\{F_t\}\colon [0,1]\times {\cal T}_{\d_0}\to {\cal T}_{\d_0}$ be a
$\cC^2$ deformation retraction of a tube ${\cal T}_{\d_0}$
onto $M$, with $F_1$ the identity map and
$\pi = F_0 \colon {\cal T}_{\d_0} \to M$ a retraction onto $M$. Let
$i_t\colon M\to [0,1]\times M$ be the map $x\to (t,x)$. By lemma
2.1 the form
$\tilde v = \int_0^1 i_t^*({\di\over \di t} \rfloor F^*u) dt$
solves $d\tilde v=u - \pi^* u$ in ${\cal T}_{\d_0}$. In the
special local coordinates provided by lemma 2.4, if
$u=\sum_{|I|+|J|=p} u_{I,J}dx^I\wedge dy^J$, the components of
$\tilde v$ are linear combinations of terms
$y_j\int_0^1 t^{|J|-1} u_{I,J}(x,ty)dt$ for $j\in J$.
Since the variables $y_j$ are transverse to $M$, we have
$|y_j|=O(\d)$ on  ${\cal T}_\d$ and hence
$||\tilde v||_{L^\infty({\cal T}_{\d})} \le C \d
 ||u||_{L^\infty({\cal T}_{\d})}$,
with $C$ independent of $\d$.
Replacing $\d$ by $c_1\d$ and changing $C$ in each step below
if necessary (but keeping it independent of $\d$),
it follows from Cauchy's inequalities that
$||Du||_{L^\infty({\cal T}_{c_1\d})} \le C \d^{-1}
||u||_{L^\infty({\cal T}_{\d})}$ and thus
$$
   ||D\tilde v||_{L^\infty({\cal T}_{c_1\d})} \le
                   C||u||_{L^\infty({\cal T}_{\d})}, \qquad
   ||\tilde v||_{\L^s({\cal T}_{c_1\d})} \le
     C \d^{1-s} ||u||_{L^\infty({\cal T}_{\d})}.
$$

In part (a) of theorem 5.1 we have $p>m$, and hence
$\pi^*u =\pi^*(i^*u)=0$ by degree reasons; so the form
$v_1=\tilde v$ satisfies $dv_1=u$ and the estimate (5.4)
with $A_\d=C\d||u||_{L^\infty({\cal T}_{\d})}$. Lemma 5.2
now completes the proof in this case.

To prove case (b) we set $v_1=\tilde v+\pi^* v_0$, where
$v_0 \in \L^s(M)$ solves $dv_0=i^*u$. We get
$$ \eqalign{ ||v_1||_{L^\infty({\cal T}_{c_1\d})} &\le
             C\left( \d||u||_{L^\infty({\cal T}_{\d})}
               + ||v_0||_{L^\infty(M)} \right), \cr
             ||v_1||_{\L^s({\cal T}_{c_1\d})} \,\, &\le
              C\left( \d^{1-s} ||u||_{L^\infty({\cal T}_{\d})}
               + ||v_0||_{\L^s(M)} \right). \cr}
$$
Lemma 5.2 then provides a holomorphic solution of $dv=u$
satisfying the estimates (5.2).

Finally, to prove part (c) in theorem 5.1, we shall construct
a good solution of $dv_0=i^*u$ on $M$ belonging to $\L^s(M)$
and apply (b). In order to circumvent problems caused
by low differentiability of $M$ we use the following
result of Whitney [Wh2]:
\sl If $M$ is a compact $\cC^k$ manifold, $k\ge 1$,
possibly with boundary, the underlying topological manifold
may be given a structure of a $\cC^\infty$ manifold, denoted $M_0$,
such that the set-theoretical identity map
$i_0\colon M_0\to M$ is a $\cC^k$-diffeomorphism. \rm

Let $i_0\colon M_0\to M$ be as above. We choose a smooth
Riemann metric on $M_0$. We refer to Wells [We] for what
follows. Let $d^*$ denote the Hilbert space adjoint of the exterior
derivative $d$ with respect to the corresponding inner product
on forms. The {\it Laplace operator} $\triangle=d^*d+dd^*$ has a
corresponding Green's operator
$G\colon L^2_{(p)}(M_0)\to H^2_{(p)}(M_0))$ with the property
that $\b=d^*G(\a)$ is the solution of $d\b=\a$ with minimal
$L^2$-norm (orthogonal to the null-space of $\triangle$),
provided that the equation $d\b=\a$ is (weakly) solvable.
For further details see section 4.5 in [We].

The Green operator is a classical pseudodifferential
operator of order $-2$, so it induces bounded operators
$L^\infty\to \L^2$ and $\L^s\to \L^{s+2}$ for $s>0$.
(See [S], section VI, 5.3.) Now $i_0^*u$ is a $\cC^1$-form
on $M_0$, and $v_0=(i_0^{-1})^*(d^* G i_0^* u)$ is a
$\cC^1$-form on $M$ with $dv_0=i^*u$, satisfying
$$ ||v_0||_{\L^s(M)} \le C'_s ||i^* u||_{L^\infty(M)}
                     \le C_s ||u||_{L^\infty({\cal T}_\d)}.
$$
Substituting this into (5.2) gives (5.3).
\endpr

\demo Remarks: 1.\ In general the constant $C=C_\d$ in the
estimate (5.3) cannot be chosen so that $\lim_{\d\to 0}C_\d=0$.
To see this, let $i^*u\ne0$ and choose a form $\phi\in \cC^1_{(m-p)}(M)$
with $\int_M u\wedge \phi \ne 0$. If $v_\d$ solves $d v_\d =u$
in ${\cal T}_{c\d}$ and satisfies
$\lim_{\d\to 0} ||v_\d||_{L^\infty({\cal T}_\d)} =0$, we get
$$ \int_M u\wedge \phi = \int_M dv_\d \wedge \phi =
   \pm \int_M v_\d \wedge d\phi \to 0
$$
as $\d\to 0$, a contradiction.

\noi 2. If $M$ is only of class $\cC^1$, the operator $d$ is not
well defined on $M$. Instead, we call a $p$-form $\a$ on $M$ exact
if there exists an integrable $(p-1)$-form $\b$ on $M$ such that
for each smooth $(m-p)$-form $\phi$ on a \nbd\ of $M$ we have
$\int_M \b\wedge i^*(d\phi) = (-1)^p \int_M \a\wedge i^*\phi$.
Then it is not hard to verify that $d(i^*_0 \b)=i_0^* \a$ (weakly)
on $M_0$, and also $d(\pi^*\b)=\pi^* \a$ on ${\cal T}_{\d_0}$.
Using this, the proof carries over with only minor changes
to the case when $M$ is of class $\cC^{1+\e}$ for some $\e>0$,
when $dv_0=i^*u$ is interpreted as above.

\noi 3. If $M$ is of class $\cC^{2+\e}$ for some $\e>0$, a more
refined argument gives a holomorphic solution of $dv=u$
that also satisfies
$||v||_{\cC^1({\cal T}_{c\d})} \le
C\log(1/\d) ||u||_{L^\infty({\cal T}_\d)}$
whenever $i^* u$ is exact. This reflects the fact that one
expects to `gain almost a derivative' in the interior estimates
for the $d$-equation. We cannot establish such estimates with
a constant independent of $\d$. In fact, when
$M = \{z\in\C^n\colon |z_j|=1,\ 1\le j\le n\}$,
this would lead to the estimate
$||\b||_{\cC^1(M)} \le const ||\a||_{L^\infty(M)}$
for a solution of $d\b=\a$, a contradiction.

\noi 4. Small changes are needed to prove theorem 5.1 when
$M=K$ is a compact subset of a larger \tr\ $\cC^2$-submanifold
$M'\subset\C^n$ (see remark 3 following the statement of theorem 5.1).
We follow the same proof as above, using the appropriate version
of $\dibar$-results given by remark 3 following theorem 3.1.
During the proof we shrink $K' \supset K$ and $\d>0$ several times.
In the proof of (5.2), we observe that the $L^2$-minimal solution
of $dv_0=i^*u$ in $U'_\d \cap M'$ also satisfies $d^*v_0=0$ when $p>1$,
and we may apply the interior elliptic estimates to obtain H\"older
estimates for $v_0$ in a \nbd\ of $K$. There are also arguments
to get the necessary control of $||v_0||_{L^2}$, for instance
the Hodge decomposition in a manifold with boundary.

\beginsection 6. Proof of theorem 1.5.

In this section we prove theorem 1.5. We shall adapt a method
of J.\ Moser [M] to the \holo\ setting.

Let $\omega$ be either the holomorphic volume form
$dz_1\wedge \cdots \wedge dz_n$ or the holomorphic
symplectic form $\sum_{j=1}^{n^\prime} dz_{2j-1}\wedge dz_{2j}$,
$n = 2n^\prime$. Write $M=M_0$ and let
$f\colon M=M_0\to M_1$ be a $\cC^k$-diffeomorphism as in
theorem 1.5 $(k\ge 2)$, satisfying the condition (1.6)
for some $\cC^{k-1}$-map $L\colon M\to GL(n,\C)$. Let
$i \colon M \hra \cn$ denote the inclusion. We assume in the proofs
that $M$ is compact and without boundary. As usual we denote by
${\cal T}_\d={\cal T}_\d M$ the tube of radius $\d$ around $M$.

By lemma 2.6 there is a \nbd\ $U\subset \C^n$ of $M$ and
a $\cC^k$-diffeomorphism $\wt{f} \colon U\to \wt{f}(U) \subset \C^n$
extending $f$ such that $\wt{f}$ is $\dibar$-flat on $M$ and satisfies
$(\wt{f}^*\omega)_z = \omega_z$ at all points $z\in M$.
The proof of theorem 1.2 then gives for each small $\d> 0$
a holomorphic map $F_\d^\prime \colon {\cal T}_\d \to \cn$ of
the form
$$ F_\d^\prime = \wt{f} + R_\d, \quad
   \|R_\d \|_{\cC^j({\cal T}_\d M)} = o(\d^{k-j});
   \qquad 0\leq j\leq k.                             \eqno(6.1)
$$
To prove theorem 1.5 we must construct biholomorphic maps
$F_\d$ as above which in addition satisfy $F^*_\d\omega=\omega$.
We need the following two lemmas.

%
%
\proclaim 6.1 Lemma:
{\rm (Existence of a good $\dibar$-flat extension.)}
If $\wt{f}$ is any  $\dibar$-flat $\cC^k$-extension of $f$
satisfying $(\wt{f}^*\omega)_z = \omega_z$ for all $z\in M$,
there exists another $\dibar$-flat $\cC^k$-extension $\wh{f}$ of $f$
satisfying $|\wh{f}^*\omega - \omega| = o(d_M^{k-1})$ near $M$
and $d\wt{f}_z = d{\wh{f}}_z$ for all $z\in M$.

%
%
\proclaim 6.2 Lemma:
{\rm (Approximation of a good $\dibar$-flat extension.)}
Assume that $\wt{f}$ is any $\dibar$-flat
$\cC^k$-extension of $f$ satisfying
$|\wt{f}^*\omega - \omega| = o(d_M^{k-1})$.
Then for all sufficiently small $\d >0$ there exist biholomorphic
maps $F_\d \colon {\cal T}_\d \to \cn$ with
$F_\d^*\omega = \omega$ and
$\|F_\d - \wt{f}\|_{\cC^j({\cal T}_\d M)} = o(\d^{k-j})$
for $0 \leq j \leq k$.

We postpone the proof of lemmas 6.1 and 6.2 for a moment.

\demo Proof of theorem 1.5 in the smooth case:
Let $f\colon M=M_0\to M_1$ be a $\cC^k$-diffeomorphism
as in theorem 1.5. By lemma 2.6 there is a $\dibar$-flat
extension $\wt{f}$ of $f$ satisfying $\wt{f}^*\omega = \omega$
at points of $M$. By lemma 6.1 we can modify this extension,
still denoting it $\wt{f}$, such that
$|\wt{f}^*\omega - \omega| = o(d_M^{k-1})$. Finally
we apply lemma 6.2 to get biholomorphic maps $F_\d$
in tubes ${\cal T}_\d$ around $M$ satisfying $F^*_\d\omega=\omega$
and the estimates (1.2). This proves theorem 1.5 in the
smooth case, granted that lemmas 6.1 and 6.2 hold.
We postpone the proof in the \ra\ case to the end of
this section.
\endpr

\demo Proof of Lemma 6.2:
Let $\wt{f}$ be as in lemma 6.2 and let
$F'_\d \colon {\cal T}_\d \to \cn$ (for small $\d>0$) be
holomorphic maps of the form (6.1) obtained as in the
proof of theorem 1.2. From the estimates on $R_\d$ in (6.1)
and the assumption $|\wt{f}^*\omega - \omega| = o(d_M^{k-1})$
it follows that
$$ \|(F'_\d)^* \omega - \omega\|_{\cC^j({\cal T}_\d)} =
   o(\d^{k-j-1}), \qquad 0 \leq j \leq k-1.
$$
Set $\omega^\d = (F'_\d)^*\omega$; this is a \holo\
$p$-form on ${\cal T}_\d$ which is close to $\omega$.
Choose constants $0<a<c<1$. Using Moser's method [M]
we shall construct a holomorphic map
$G_\d \colon {\cal T}_{a\d} \to {\cal T}_\d$ which is very close
to the identity map and satisfies $G_\d^* \omega^\d = \omega$
on ${\cal T}_{a\d}$. The holomorphic map
$F_\d=F'_\d\circ G_\d \colon {\cal T}_{a\d} \to \C^n$
is then close to $F'_\d$ (and hence to $\tilde f$),
and it satisfies $F_\d^*\omega=G^*_\d(\omega^\d)=\omega$.

We first outline the Moser's method, postponing the estimates
for a moment. Set $\omega^\d_1= (F'_\d)^*\omega$ and
$\omega^\d_t = (1-t)\omega + t\omega^\d_1$ for $t\in [0,1]$.
Then $d\omega^\d_t=0$, and $\omega^\d_t$ is close to $\omega$
for each $t$ and $\d$. Our goal is to construct a $\cC^1$-family
of holomorphic maps
$G_t=G_{\d,t} \colon {\cal T}_{a\d} \to {\cal T}_\d$ satisfying
$G_0 = Id$ and $G_t^*\omega^\d_t = \omega$ for all $t\in [0,1]$;
the time-one map $G_\d=G_{\d,1}$ will then solve the problem.

To simplify the notation we suppress $\d$ for the moment,
writing $\omega^\d_t=\omega_t$ and $G_{\d,t}=G_t$.
Suppose that such a flow $G_t$ exists. Denote by $Z_t$
its infinitesimal generator; this is a holomorphic
time-dependent vector field on the image of $G_t$, satisfying
${d \over dt}G_t(z)= Z_t(G_t(z))$ for each $t\in [0,1]$
and each $z$ in the domain of $G_t$.
Differentiating the equation $G_t^*\omega_t = \omega$
on $t$ and applying the time-dependent Lie derivative theorem
([AMR], Theorem 5.4.5., p.\ 372), we have
$$ 0= {d\over dt}(G_t^*\omega_t) =
   G_t^*\bigl( L_{Z_t}\omega_t + {d\over dt}\omega_t \bigr)
   = G^*_t\bigl( d(Z_t\rfloor\omega_t) + \omega_1-\omega \bigr).
                                                     \eqno(6.2)
$$
We have also used the Cartan formula for the Lie
derivative $L_{Z_t}\omega_t$, as well as $d\omega_t=0$.
This shows that $G^*_t\omega_t=\omega$ holds for all $t\in [0,1]$
if and only if the generator $Z_t$ satisfies the equation
$d(Z_t \rfloor \omega_t) + \omega_1 - \omega = 0$ for
all $t \in [0,1]$.

At this point we observe that $\omega$ is exact holomorphic
on $\C^n$, $\omega=d\beta$; in fact when $\omega$ is the volume
form (1.4) we may take $\beta={1\over n} \sum_{j=1}^n (-1)^{j+1}dz[j]$,
and when $\omega$ is the symplectic form (1.5) we may take
$\beta=\sum_{j=1}^{n'} z_{2j-1}dz_{2j}$. Hence the difference
$\omega_1-\omega=F_\d^{\prime*} d\beta - d\beta
 = d\bigl( F_\d^{\prime*}\beta -\beta \bigr)$
is exact \holo\ on ${\cal T}_\d$. By theorem 5.1 we can solve the
equation $dv = \omega_1 - \omega$ to get a small \holo\
$(p-1)$-form $v=v_\d$ in ${\cal T}_{c\d}$. Let $Z_t$ be the unique \hvf\ on
${\cal T}_{c\d}$ solving the (algebraic !) equation
$Z_t \rfloor \omega_t + v=0$. Integrating $Z_t$ we get
a flow $G_t$ satisfying $G^*_t\omega_t=\omega$ on its domain
of definition.

For this approach to work we must choose $v_\d$ on ${\cal T}_{c\d}$
to have as small sup-norm as possible; this will imply that
$|Z_t|$ is small, and hence its flow $G_t(z)$ will not
escape the tube ${\cal T}_{c\d}$ (on which $Z_t$ is defined)
before time $t=1$, provided that the initial point
$G_0(z)=z$ belongs to the smaller tube ${\cal T}_{a\d}$.
(In particular, the solution $v_\d=(F'_\d)^*\beta -\beta$
may not work since $F'_\d$ is not close to the identity map.)

In order to apply theorem 5.1 efficiently we must first
show that $dv_0 = i^*(\omega_1 - \omega)$ has a solution
on $M$ with small norm. Consider the map
$h \colon [0,1] \times M \to \cn$, $h(t,z) =\wt{f}(z) + tR_\d(z)$,
and set $w=h^*\omega$. Also let $i_t \colon M\to [0,1] \times M$
denote the injection $i_t(z)=(t,z)$ ($z\in M$, $t \in [0,1]$).
It follows from lemma 2.1 that
$v_0 = \int_0^1 i_t^*({\di \over \di t}\rfloor w) dt$
solves $dv_0 = i_1^*w - i_0^*w$. We have $i_1^*w = i^*\omega_1$
and $i_0^*w = i^*\wt{f}^*\omega = i^*\omega$, so
$dv_0 = i^*(\omega_1 - \omega)$. It follows from the formula above that
$v_0 = \sum_{j=1}^n r_j^\d v_j$, where $r_1^\d,\cdots,r_n^\d$ are
the components of $R_\d$ and $v_1,\cdots,v_n$ are $(p-1)$-forms on $M$
with $\|v_j\|_{\cC^{k-1}(M)}$ bounded independently of $\d$.
This gives $\|v_0\|_{\cC^l(M)} = o(\d^{k-l})$
for $0\leq l\leq k-1$. It follows that
$\|v_0\|_{\L^s(M)} = o(\d^{k-s})$ for a given $s\in (0,1)$.
Since $\|\omega_1 - \omega\|_{L^\infty({\cal T}_\d)} = o(\d^{k-1})$,
it follows from theorem 5.1 that for all sufficiently small
$\d>0$ we have a holomorphic solution of
$dv_\d = \omega_1 - \omega$ in ${\cal T}_{c\d}$,
satisfying $\|v_\d\|_{L^\infty({\cal T}_{c\d})} = o(\d^k)$.

Let $Z^\d_t$ be the \hvf\ in ${\cal T}_{c\d}$ satisfying
$Z^\d_t \rfloor \omega^\d_t = v_\d$. The above estimate on
$v_\d$ implies $\|Z^\d_t\|_{L^\infty({\cal T}_{c\d})} = o(\d^k)$,
uniformly in $t\in [0,1]$. The standard formula for the rate
of escape of the flow shows that we can choose
$\d_0>0$ sufficiently small such that for all $\d\in (0,\d_0)$
and all initial points $z\in {\cal T}_{a\d}$, the flow
$G_{\d,t}(z)$ of $Z^\d_t$ remains in ${\cal T}_{c\d}$ for all
$t\in [0,1]$. At $t=1$ we get a map
$G_\d= G_{\d,1} \colon {\cal T}_{a\d}\to {\cal T}_{c\d}$
satisfying $G_\d^*\omega^\d_1=\omega$ and
$|G_\d(z) - z| = o(\d^k)$ for $z\in {\cal T}_{a\d}$.

Set $F_\d = G_\d \circ F'_\d$. Since the maps $F'_\d$ have
uniformly bounded $\cC^1$-norms on ${\cal T}_\d $, we see that
$\|F_\d - F'_\d\|_{L^\infty({\cal T}_{a\d})} = o(\d^k)$.
Replacing $a$ by a smaller constant and applying
the Cauchy inequalities we also get
$$
\|F_\d-\wt{f}\|_{\cC^j({\cal T}_{a\d}}) \le
   \|F_\d - F_\d^\prime\|_{\cC^j({\cal T}_{a\d})}
   + \|F_\d^\prime - \wt{f}\|_{\cC^j({\cal T}_{\d})} = o(\d^{k-j}),
   \qquad j \leq k.
$$
By construction we have $F_\d^*\omega=\omega$,
so $F_\d$ solves the problem.
\endpr

\demo Remark: This method applies on any domain $D\ss \C^n$ on which
we can solve the $\dibar$-equations with estimates (e.g., on \psc\
domains); it shows that for any holomorphic map
$F'\colon D\to \C^n$ for which $|F^{\prime*}\omega -\omega|$
is sufficiently uniformly small on $D$ there exists a \holo\
map $F\colon D'\to \C^n$ on a slightly smaller domain $D'\ss D$
such that $F^*\omega=\omega$ and $F$ is uniformly close to $F'$
on $D'$. We obtain $F$ in the form $F=F'\circ G$, where
$G\colon D'\to D$ is a holomorphic map close to the identity,
chosen such that $G^*(F^{\prime *}\omega)=\omega$. The precise amount of
shrinking of the domain depends on
$\|F^{\prime *}\omega-\omega\|_{L^\infty(D)}$
and on the constants in the solutions of the $\dibar$-equations; we
do not know if there is a solution to this problem on all of $D$.
\medskip

We now turn to the proof of lemma 6.1. We shall need
the following:

\proclaim 6.3 Lemma: Let $u$ be a $d$-closed $p$-form of class $\cC^{k-1}$
in a neighborhood of $M$, with $p\geq 1$, such that the $(p,0)$-component
$u^\prime$ of $u$ is $\dibar$-flat on $M$, and $u^\dpr=u-u^\pr$ is
$(k-1)$-flat on $M$. Assume $i^*u=0$, where $i\colon M\hra\C^n$
is the inclusion. Then there exists a $(p-1,0)$-form
$v$ in a neighborhood of $M$ such that $v=\sum_{j=1}^N\z_jv_j$, where
each $\z_j$ is a $\dibar$-flat $\cC^k$-function vanishing on $M$, each
$v_j$ is a $\dibar$-flat $\cC^{k-1}$-form, and $|u-dv|=o(d_M^{k-1})$.
If $u=0$ on $M$, we may take $v_j=0$ on $M$ for all $j$.

\demo Remark:  Using the rough multiplication (lemma 2.5) we see that
there is a $\dibar$-flat $(p,0)$-form $v$ of class $\cC^k$ that also satisfies
$|dv-u|=o(d_M^{k-1})$. However, the version stated above is often technically
more convenient since we may wish to postpone the use of rough multiplication.

\demo Proof of Lemma 6.3:
In the case $m=n$ we may take $v=0$ which can be seen as follows.
We have $u^\prime=\sum_{|I|=p}u_Idz^I$, where the
coefficients $u_I$ are $\cC^k$-functions that are $\dibar$-flat on
$M$; hence $i^*u=0$ means that $u_I=0$ on $M$ for all $I$
(since the coefficients of $u^\dpr$ vanish on $M$).
It follows from the Cauchy-Riemann equations that each $u_I$
is flat on $M$, so we way choose $v=0$.

When $m<n$, we use the asymptotically holomorphic extension $\wt{M}$ of $M$
(lemma 2.4) and the $\dibar$-flat retraction $F$ to $\wt{M}$. Recall
that a neighborhood of $M$ may be covered by $\cC^k$-charts
$G_i \colon U_i \to V_i$,
$G_i(z)= (z_{(i)}^\pr(z),w_{(i)}^\pr(z)) \in \C^m \times
\C^{n-m}$, $1\le i\le r$, satisfying

\item{--} $G_i$ is $\dibar$-flat on $M$,
$G_i(M\cap U_i) = V_i\cap(\Rm\times\{0\})$,
$G_i(\wt{M}\cap U_i) = V_i\cap(\cm\times\{0\})$;

\item{--} the retraction $F$ is given in these local coordinates
by $(t,(z^\pr,w^\pr)) \to (z^\pr,tw^\pr)$.

\medskip
Let $\wt{i} \colon \wt{M} \hra \cn$ be the inclusion. Arguing as in the case
$m=n$ and making use of the $\dibar$-flat local parametrizations of $\wt{M}$,
we see that $\wt{i}^*u$ is flat on $M$, and so is $\wt{\pi}^*u=\wt{\pi}^*\wt{i}^*u$,
where $\wt{\pi}=F_0$. When $F\colon [0,1]\times W \to W$ is the retraction to
$\wt{M}$, the form
$$
   \wh{v}=\int_0^1i_t^*({\partial\over \partial t}\rfloor F^*u)dt \eqno(6.3)
$$
solves $d\wh{v}=u-\wt{\pi}^*u$ on a neighborhood of $M$ according to
lemma 2.1. Expressing $u$ in the $G_i$-coordinates
$(z_{(i)}^\pr(z),w_{(i)}^\pr(z))$ (which are $\dibar$-flat
on $M$) we get
$$
     u=\sum_{|I|+|J|=p} a_{I,J}(z_{(i)}^\pr, w_{(i)}^\pr)
     dz_{(i)}^{\prime I}\wedge dw_{(i)}^{\prime J}+ r'_{(i)}
$$
on $U_i$, where the $a_{I,J}$ are $\cC^{k-1}$-functions that are
$\dibar$-flat on $\Rm\times\{0\}$ and $r'_{(i)}$ is a $\cC^{k-1}$-form
that is flat on $M$. Using the formula following lemma 2.1 we see
that $\wh{v}$ (6.3) is a linear combination of terms
$$
   w_{(i),j}^\pr \left(
   \int_0^1a_{I,J}(z_{(i)}^\pr,tw_{(i)}^\pr)t^{|K|}dt \right)
   \, dz_{(i)}^{\prime I}\wedge dw_{(i)}^{\prime K},
$$
where $|I|+|K|=p-1$ and $1\le j\le n-m$, plus a remainder
term $r_{(i)}^\dpr$ satisfying $|\partial^\a r_{(i)}^\dpr|=o(d_M^{k-|\a|})$
on $U_i$ for $|\a|\leq k-1$. Here $w_{(i),j}^\pr$ denotes
the $j$-th component of $w'_{(i)}$. Since $G_{(i)}$ is $\dibar$-flat,
it follows that
$$ \wh{v}= \sum_{j=1}^{n-m} \sum_{|L|=p-1}
   w_{(i),j}^\pr g_{j,L}^{(i)}dz^L + r_{(i)}
$$
in $U_i$, where each $g_{j,L}^{(i)}$ is a $\dibar$-flat
$\cC^{k-1}$-function and $r_{(i)}$ behaves like $r_{(i)}^\dpr$.

Choose a $\dibar$-flat partition of unity $\{\psi_i\}_{i=1}^r$
subordinate to the covering $\{U_i\}_{i=1}^r$, and choose
$\dibar$-flat cut-off functions $\chi_i \in \cC_0^\infty(U_i)$,
with $\chi_i=1$ near ${\rm supp}\psi_i\cap M$ for $i=1,\cdots,r$.
Let $\z_1,\cdots,\z_N$ (with $N=r(n-m)$) be some enumeration of
the collection of functions
$\{\psi_i w_{(i),j}^\pr \colon i\leq r, j\leq n-m\}$.
Furthermore let $v_1,\cdots,v_N$ be the corresponding
enumeration of the forms $\chi_i\sum_{|L|=p-1}g_{j,L}^{(i)}dz^L$,
prolonged by zero outside $U_i$. Set $v=\sum_{l=1}^N \z_l v_l$.
Clearly $|dv-u|=o(d_M^{k-1})$. Furthermore, if $u=0$ on $M$, we also
see that $\int_0^1a_{I,J}(z^\pr,tw^\pr)t^{|K|}dt=0$ on
$V_i\cap(\Rm\times\{0\})$, and hence $v_1=\cdots=v_N=0$ on $M$.
\endpr

\demo Proof of Lemma 6.1:
In the unimodular case, $\omega=dz_1\wedge\cdots\wedge dz_n$,
we could successively increase the order of vanishing of
$\wt{f}\omega - \omega$ on $M$ by adding certain correction terms
to $\wt{f}$. This seems harder to do in the symplectic case, so we
shall instead present an argument that works uniformly in both cases.
It is a modification of Moser's method:
With $\omega_t = (1-t)\omega + t\wt{f}^*\omega$, we shall
construct a $\cC^1$-family of $\dibar$-flat $\cC^k$-maps $g_t$ on a
neighborhood of $M$, with $g_0 = id$ and
$|{d \over dt}g_t^*\omega_t|=o(d_M^{k-1})$ uniformly in $t$.
Given such a family, integration in $t$ gives
$\|g_1^*\omega_1 - \omega\| = o(d_M^{k-1})$.
We will also show that $g_1$ is $\dibar$-flat on $M$.
Hence the map $\wh{f} = \wt{f}\circ g_1$ will satisfy lemma 6.1.
Furthermore, we shall see that $|g_1(z) - z| = O(d_M(z)^2)$, so
$Dg_1 = Id$ on $M$ and hence $\wh{f}$ and $\wt{f}$ have the same
differential on $M$.

We shall obtain $g_t$ by integrating a certain real time-dependent
vector field $X_t$ of class $\cC^k$. Differentiating
${d \over dt}g_t^*\omega_t$ as in (6.2) we see that $X_t$ must satisfy
$|d(X_t\rfloor\omega_t)+\omega_1-\omega|=o(d_M^{k-1})$.
We shall now construct such a vector field. More precisely,
we shall construct a continuous family of $\cC^k$ real vector fields
$X_t$ on a tube ${\cal T}_0={\cal T}_{\d_0}$, satisfying
the following properties for each $t\in [0,1]$:

\item{(1)} $X_t$, considered as a map ${\cal T}_0 \to \cn$, is $\dibar$-flat
on $M$. (Here we identify a real tangent vector
$X=\sum_{j=1}^n a_j{\di/\di x_j}+b_j{\di/\di y_j} \in T_z\C^n$
with the corresponding complex vector $(a_1+ib_1,\ldots,a_n+ib_n)\in \C^n$.)

\item{(2)} $|X_t(z)| \leq C d_M(z)^2$ for some $C>0$ independent
of $t\in [0,1]$.

\item{(3)} $|d(X_t\rfloor\omega_t)+\omega_1 -\omega|=o(d_M^{k-1})$,
uniformly in $t \in [0,1]$.

\medskip
Let us first show that this solves the problem. We must show that
$X_t$ can be integrated from $t=0$ to $t=1$ for all initial values
in a smaller tube. Recall that, after shrinking $\d_0$ if necessary,
the function $d_M$ is differentiable in ${\cal T}_0\bs M$, with a gradient
of length one. Let $z(t)$ be an integral curve of $X_t$ in ${\cal T}_0\bs M$,
$t \in [0,t_0]$, and set $u(t) = d_M(z(t))$. Then
$$ u^\prime(t) = \nabla d_M(z(t))\cdotp X_t(z(t))
   \leq |X_t(z(t))| \leq C u(t)^2.
$$
Here we denoted by $v\cdotp w$ the real inner product of the
vectors $v,w \in \C^n$. Integrating the inequality
$u^\prime(t)/u(t)^2 \leq C$ from $0$ to $t$ gives
$1/u(0) - 1/u(t) \leq Ct$ and thus $u(t)(1-Ctu(0)) \leq u(0)$
for $0\leq t\leq t_0$. Let the initial value
$z(0) \in {\cal T}_{\d_1}\bs M$, where
$\d_1 \leq \min(\d_0/2, 1/2C)$. It follows that
$u(t) \leq u(0)/(1-Ctu(0)) \leq 2u(0)$, and hence the integral curve
extends to all values $t \in [0,1]$. Since $|X_t(z(t))| \leq Cu(t)^2$,
we see that $|z(t)-z(0)| \leq 4Cu(0)^2t$. In other words, the
time-$t$ diffeomorphisms $g_t$ are well defined on ${\cal T}_{\d_1}$
for all $t\in [0,1]$ and they satisfy
$|g_t(z)-z| \leq 4Ctd_M(z)^2$. In particular, $g_t(z) = z$ and
$Dg_t(z) = Id$ for $z \in M$ and $t\in [0,1]$.

To show that the $\cC^k$-maps $g_t$ are $\dibar$-flat on $M$,
we consider the variational equation
${\partial \over \partial t}D_zg_t(z)=D_zX_t(g_t(z))\circ D_zg_t(z)$
with the initial condition $D_z g_0=Id$. Decomposing the differential
$D\phi$ as the sum of a $\C$-linear part $D^\prime\phi$ and a
$\C$-conjugate part $D^{\prime\prime}\phi$, we get
$$ \eqalign{
  {\di \over \di t} D_z^{\prime\prime}g_t(z) &=
  D_z^{\prime\prime}\bigl( {\di \over \di t} g_t(z)\bigr)=
   D_z^{\prime\prime}\bigl( X_t(g_t(z))\bigr) \cr
  &= (D_z^\prime X_t)(g_t(z))\circ D_z^{\prime\prime}g_t(z)
  + (D_z^{\prime\prime}X_t)(g_t(z))\circ D_z^\prime g_t(z). \cr}
$$
We apply both sides to a unit vector $v\in \C^n$ and set
$y(t) = D_z^{\prime\prime}g_t(z)v \in \C^n$. We obtain a linear
differential equation $y'(t)=A(t)y(t)+b(t)$ with the initial
condition $y(0)=D_z^{\prime\prime}g_0(z) v=0$.
The function $u(t)=|y(t)|$ is differentiable when $u(t)\neq 0$
and $u^\prime(t)=y'(t)\cdotp y(t)/|y(t)| \leq |y'(t)|$.
Thus, if $|A(t)|\leq A$ and $|b(t)|\leq b$, we see that
$u'(t)\leq Au(t)+b$ where $u(t)\neq 0$. We shall prove that
$u(t)\leq {b \over A}(e^{At}-1)$, $t \in [0,1]$. If $u(t)=0$, there is
nothing to prove. If not, let $t_0$ be the largest zero of $u$
on the interval $[0,t]$. Then $v(s)=u(s)a^{-As}$ satisfies the differential
inequality $v^\prime(s)\leq be^{-As}$ for $s \in (t_0,t]$. Integration from
$t_0$ to $t$ gives $v(t)\leq {b \over A}(e^{-At_0}-e^{-At})$ and
$u(t)\leq {b \over A}(e^{A(t-t_0)}-1)\leq {b \over A}(e^{At}-1)$.

We know that $|D_zX_t(z)|$ and $|D_zg_t(z)|$ are bounded uniformly in
$z \in {\cal T}_{\d_1}$ and $t \in [0,1]$, while $|D_z^{\prime\prime}X_t(z)=
o(d_M(z)^{k-1})$. Thus we may choose the upper bound $A$ for $|A(t)|$
independently of $z \in {\cal T}_{\d_1}$ and the unit vector $v$,
and we may choose the upper bound $b$ of $|b(t)|$ to be of
size $b=o(d_M(z)^{k-1})$, uniformly in $v$.
Since $u(t)=|D_z^{\prime\prime}G_t(z)v|$, it follows
that $|D_z^{\prime\prime}g_t(z)|=o(d_M(z)^{k-1})$, so each $g_t$ is
$\dibar$-flat on $M$.

By assumption we have
$|d(X_t\rfloor\omega_t)_z+(\omega_1-\omega_0)_z|= o(d_M(z)^{k-1})$.
Since $d_M(g_t(z)) \leq 2d_M(z)$ and the norms $|D_zg_t(z)|$ are
bounded uniformly in $z \in {\cal T}_{\d_1}$ and $t \in [0,1]$,
we have
$|{\partial\over\partial t}(g_t^*\omega_t)_z|=o(d_M(z)^{k-1})$,
uniformly in $t$. By integration in $t$ we obtain
$|(g_1^*\omega_1-\omega)_z|=o(d_M(z)^{k-1})$.
Setting $\wh{f}=\wt{f}\circ g_1$, we see that $\wh{f}$ is a $\dibar$-flat
$\cC^k$-extension of $\wt{f}$, $D\wh{f}=D\wt{f}$ on $M$, and
$|(\wh{f}^*\omega-\omega)_z|=o(d_M(z)^{k-1})$.
Thus $\wh{f}$ satisfies lemma 6.1.

It remains to construct the vector field $X_t$. Applying lemma 6.3
to $\omega-\omega_1$ we obtain a $(p-1,0)$-form $v$ near $M$
with $|dv-(\omega-\omega_1)|=o(d_M^{k-1})$ and $v=0$ on $M$.
We decompose $\omega_t$ as $\omega_t^\pr +\omega_t^\dpr$, where
$\omega_t^\pr$ is the $(p,0)$-component of $\omega_t$.
Then $\omega_t^\pr=\omega+ t(\omega_1^\pr-\omega)$,
and $\omega'_t = \omega$ on $M$ for each $t$. Hence the map
$\phi \colon Z \to Z\rfloor\omega_t^\pr$, taking the $(1,0)$-vectors
$Z\in T^{(1,0)}_z\C^n$ to $\L^{(p-1,0)}T_z^*\cn$,
is an isomorphism for $z$ near $M$ and $t \in [0,1]$.
Hence the equation $Z_t^\pr\rfloor\omega_t=v$ uniquely
defines a time-dependent $(1,0)$ vector field $Z_t^\pr$
on $\C^n$ near $M$.

With respect to the basis
${\partial\over\partial z_1},\cdots,{\partial\over\partial z_n}$ for
$(1,0)$-vectors and the basis $dz[1],\cdots,dz[n]$
(respectively $dz_1,\cdots,dz_n$) for the $(p-1,0)$-covectors,
the map $\phi$ is represented by an $(n \times n)$-matrix
valued function $A(t,z) = A_0+tB(z)$, where $A_0$ is constant and
invertible, and the entries of $B(z)$ are $\dibar$-flat $\cC^{k-1}$-functions
that vanish on $M$. It follows that the entries of $A(t,z)^{-1}$ are rational
functions $b(t,z)$ in $t$, with coefficients that are $\dibar$-flat
$\cC^{k-1}$-functions. From the properties of $v$, as given by lemma 6.3,
it follows that
$Z_t^\pr=\sum_{j=1}^N\z_j\sum_{k=1}^nr_{jk}(t,z)\partial/\partial z_k$,
where $\z_1,\cdots,\z_N$ are $\cC^k$-functions that vanish on $M$
and are $\dibar$-flat on $M$, and each $r_{jk}$ is a rational function
in $t$ with coefficients that are $\dibar$-flat $\cC^{k-1}$-functions
with $r_{jk}(t,z)=0$ for $z \in M$.

We next apply the rough multiplication lemma to the pairs
$(\z_j(z),r_{jk}(t,z))$ with respect to the compact subset $M \times [0,1]$
in $\cn \times \R$ and obtain $\cC^k$-functions $a_l(z,t)$, $1\leq l\leq n$,
$\dibar$-flat on $M$ with respect to $z$, such that
$|\sum_{j=1}^N\z_jr_{jl}(t,\cdotp)-a_l(t,\cdotp)|=o(d_M^k)$, uniformly in $t$.
(Remark : Use of the parametrized version of rough multiplication gives a
smooth family of $\cC^k$-functions, but we do not need that.)

We set $Z_t=\sum_{l=1}^na_l(t,z)\partial/\partial z_l$ and
$X_t=Z_t+{\bar Z}_t$. Writing $a_l=u_l+iv_l$, with
$u_l$ and $v_l$ real, we have
$X_t=\sum_{l=1}^n u_l(z,t)\di/\di x_l + v_l(z,t)\di/\di y_l$.
If we consider $X_t$ as a map ${\cal T}_0\to \cn$, this means
that $X_t=(a_1(t,z),\cdots,a_n(t,z))$ and is $\dibar$-flat on $M$.
Furthermore, since $\z_j(z)$ and $r_{jl}(t,z)$ both vanish when $z\in M$,
we see that $X_t$ vanish to the second order on $M$.

Finally, we must show that (3) is satisfied. Writing
$X_t={\bar Z}_t+(Z_t-Z_t^\pr)+Z_t^\pr$, we see that
$$
   d(X_t\rfloor\omega_t)+\omega_1-\omega =
   d({\bar Z}_t\rfloor\omega_t^\dpr)+
   d((Z_t-Z_t^\pr)\rfloor\omega_t^\pr)+(dv+\omega_1-\omega).
$$
The first term on the right hand side is $o(d_M^k)$
since $\omega_t^\dpr$ vanishes to order $k-1$ and $Z_t$ vanishes
to the second order on $M$. Furthermore, $Z_t-Z_t^\pr$ vanishes to
the $k$-th order, so the second term is $o(d_M^{k-1})$, and the
third term is $|dv+\omega_1-\omega|=o(d_M^{k-1})$. Thus (3) holds,
uniformly in t, since the derivatives are continuous in $(z,t)$.
\endpr

\demo Proof of theorem 1.5 in the \ra\ case:
By assumption there is a continuous map
$\psi_0 \colon M_0 \to SL(n,\C)$
(resp.\ $\psi_0\colon M_0\to Sp(n,\C)$) such that
$\psi_{0,z}$ agrees with $d_z f$ on $T_zM_0$ for each
$z\in M_0$. By Remark 1 following the proof of theorem 1.2
(sect.\ 4) $\psi_0$ may be approximated by a holomorphic map
$\psi_1$ from a neighborhood of $M_0$ to $GL(n,\C)$
with $\psi_{1,z}=d_z\wt{f}$ on $T_z \wt{M}_0$ for each
$z\in \wt{M}_0$. Since $\psi_{0,z}^*\omega=\omega$
for $z\in M_0$ and since $\psi_1$ approximates $\psi_0$
on $M_0$, it follows that the form
$\psi_{1,z}^*\omega = (\det\psi_{1,z}) \omega$ is
close to $\omega$ for all $z\in \wt{M}_0$ sufficiently
near $M_0$.

We may think of $\psi_1$ as a \holo\ automorphism of the trivial
bundle $\wt{M}_0\times \cn \to \wt{M}_0$. We claim that there is
another holomorphic automorphism $g$ of $\wt{M}_0\times \cn$
such that $g|_{T\wt{M}_0}=Id$ and $g^*\psi_1^*\omega = \omega$.
In the unimodular case we let $g$ act as the identity
on $T \wt{M}_0$ and as multiplication by
$(\det \psi_1)^{-1 /(n-m)}$ on $\wt{\nu}_0$
(the holomorphic extension of the complex normal
bundle $\nu_0$ to $\wt{M}_0$); the root is well-defined
since the function $\det \psi_{1,z}$ is close to $1$.
In the symplectic case $g$ is a reduction to symplectic normal
form with holomorphic dependence on $z \in \wt{M}_0$.
In both cases the map $\psi = \psi_1 \circ g$ is an automorphism
of the trivial bundle $\wt{M}_0\times \cn$ satisfying
$\psi^*\omega=\omega$.

Let $F_1$ be a biholomorphic extension of $\wt{f}$
constructed from $\psi = \psi_1 \circ g$ as in Remark 1 (sec.\ 4),
satisfying $d_z F_1 = \psi_z$ at points $z\in \wt{M}_0$.
Thus $F_1^*\omega = \omega$ at points of $\wt{M}_0$. Applying Moser's
method as above we can construct a biholomorphism $G$ in a
tubular \nbd\ of $\wt{M}_0$ which equals the identity on
$\wt{M}_0$ and satisfies $G^*(F_1^*\omega )= \omega$.
Then $F=F_1\circ G$ is a biholomorphic map near $M_0$
which extends $f$ and satisfies $F^*\omega=\omega$.
\endpr

\beginsection 7. Proof of theorems 1.7 and 1.8.

We will have to consider maps which have different degree of
smoothness with respect to the time variable and the space variable,
and we shall use the following terminology.

%
%
\proclaim Definition 4:
Let $U$ be an open subset of $[0,1] \times \R^m$.
A mapping $f\colon U\to\R^n$ is called a {\it $\cC^l$-family
of $\cC^k$-maps} if $\di_t^j (\di_x^{\a}f)$ is continuous in $U$
for $0\le j \le l$ and $|\a | \le k$. There is an obvious extension
of this notion to maps $f\colon [0,1] \times M\to N$ where
$M$ and $N$ are $\cC^k$ manifolds. If in addition $f_t=f(t\cdotp)$
is a diffeomorphism (of its domain onto its image) for each
$t\in [0,1]$, we call $f=\{f_t\}$ a {\it $\cC^l$-family
of $\cC^k$-diffeomorphisms}.

Thus a $\cC^1$-family of $\cC^k$-diffeomorphisms is the same
as a $\cC^k$-isotopy (or a $\cC^k$-flow) in the sense of
definition 1 in sect.\ 1. We remark that if $f_t$ is a $\cC^l$-family
of diffeomorphisms on domains $U_t \subset \R^n$ for
$t\in [0,1]$, the family of inverses $f_t^{-1}$ are not
necessarily a $\cC^l$-family if $l>0$, the reason being that
the $t$-derivatives of the (derivatives of the) inverse map
will involve higher order $x$-derivatives of the original map.

In the situation in theorem 1.7 we shall say that a
time-dependent family of $\cC^k$-forms on submanifolds
$M_t\subset \C^n$, $\a_t = \sum_{|I|=p} \a_{I,t}dz^{I}$
with $\a_{I,t} \in \cC^k(M_t)$, is a {\it continuous family
of $\cC^k$-forms} if $\a_{I,t}\circ f_t$ is a continuous
family of $\cC^k$-functions on $M$ for all multiindices $I$.
Recall that ${\cal T}_\d={\cal T}_\d M$ is the open tube
of radius $\d$ around a submanifold $M\subset \C^n$.

The main step in the proof of Theorem 1.7 is the following result.

\proclaim 7.1 Theorem:
Let $f_t \colon M=M_0 \to M_t \subset \C^n$ $(t\in [0,1])$
be a $\cC^1$-family of $\cC^k$-diffeomorphisms between compact,
totally real, $\cC^k$-submanifolds of $\cn$, with $f_0$ the identity
on $M$. By $i_t \colon M_t \hra\C^n$ we denote the inclusion map.
Let $\a_t$ ($t\in [0,1]$) be a continuous family of $(p,0)$-forms
of class $\cC^k$ on $M_t$ such that $\istar_t \a_t$ is closed
on $M_t$ for each $t$. Then there exists an extension of $\a_t$ to
a continuous family $\wh{\a}_t$ of $(p,0)$-forms of class $\cC^k$ on
a neighborhood of $\wt{M} = \bigcup_{t \in [0,1]} \{t\} \times M_t$
in $[0,1]\times\C^n$ such that for all sufficiently small $\d > 0$
there exists a continuous family of closed holomorphic $p$-forms
$u_t^\d$ on
$U_{\d} = \bigcup_{t \in [0,1]} \{t\} \times {\cal T}_{\d}M_t$
satisfying
$$ \|u^{\d}_t - \wh{\a}_t \|_{\cC^r({\cal T}_\d M_t)} = o(\d^{k-r}),
   \qquad 0\le r\le k,
$$
uniformly in $t\in [0,1]$. If $i_t^*\a_t$ is exact on $M_t$
for each $t\in [0,1]$, we may choose $u_t^{\d}$ exact for every
$t$; in this case $u_t^\d$ can be chosen to be entire if
each $M_t$ is polynomially convex.

In the simplest case when $M_t = M$ and $\a_t = \a$ for all
$t\in [0,1]$, the main steps in the proof of theorem 7.1 are
as follows (we write ${\cal T}_\d={\cal T}_\d M$):

\item{(i)} We construct a $(p,0)$-form $\wh{\a}$ on a neighborhood
of $M$ such that $d\wh{\a}$ is flat on $M$. In particular, $\wh{\a}$
is $\dibar$-flat on $M$;

\item{(ii)} we approximate the coefficients of $\wh{\a}$ by
holomorphic functions to obtain a holomorphic $p$-form
$u^{\prime}$ in ${\cal T}_\d$ with
$\|du^{\prime}\|_{L^\infty({\cal T}_\d)} = o(\d^{k-1})$;

\item{(iii)} we solve $dv = du^{\prime}$, with $v$
holomorphic and $\|v\|_{L^\infty({\cal T}_\d} = o(\d^{k})$, and
set $u = u^{\prime}-v$;

\item{(iv)} if $i^*\a$ is exact, the norm of the de Rham cohomology
class of $i^*u$ is $o(\d^k)$, and this class may be represented by a
holomorphic $p$-form $u_0$ on ${\cal T}_\d$ of size $o(\d^k)$.
Then $u_1 = u-u_0$ is exact and it approximates $\a$ to the right
order on $M$.

\medskip
In the parametric case we perform these steps such that the
solutions are continuous with respect to the parameter $t$.
Before giving the proof of theorem 7.1 we summarize (slight extensions of)
certain well known results that we shall need.

We begin by considering the {\it parameter dependence in
Whitney's extension theorem.} Instead of a general compact subset
$K \subset \Rn$ (or $K\subset \C^n$ we consider the case when $K$
is a compact $\cC^1$-submanifold, with or without boundary.
This is a so-called 1-regular set, so we have the following more
precise results (see [T], chapter IV, sec.\ 1 and 2,
in particular p.\ 76):

\item{(i)} Let $A = \{ \a \in \zn_+ \colon |\a|\le k \}$. The
collections $F = (f_{\a})_{\a\in A} \in \cC(K)^A$,
satisfying the Whitney condition, form a closed subspace $\cE^k(K)$
of $\cC(K)^A$ with respect to the sup-norm; we shall call such
collections {\it Whitney functions.}

\item{(ii)} The Whitney extension operator
${\cW}\colon \cE^k(K) \to \cC^k_0(K^{\prime})$,
where $K^{\prime} \subset \R^n$ is a closed neighborhood of $K$,
is linear and norm-continuous. Thus $\di^{\a}{\cW}(F) = f_{\a}$ on $K$
for each $\a\in A$ and
$$
   \|{\cW}(F)\|_{\cC^k(K^{\prime})} \le C
   \sup \{\|f_{\a}\|_{L^\infty (K)}\colon |\a| \le k\}.
$$

\item{(iii)} There exists a constant $C>0$ such that $C\omega$
is a modulus of continuity for $\di^{\a}{\cW}(F)$, $|\a|=k$,
whenever $\omega$ is a modulus of continuity for
all $f_{\a}$, $|\a|=k$.

\medskip
From this it follows immediately that if $f_{\a,t}$, $\a \in A$,
are $\cC^l$-families of continuous functions on $K$ and if
$F_t=(f_{\a,t})_{\a \in A}$ is a Whitney function for each
$t \in [0,1]$, then their Whitney extensions ${\cW}(F_t)$ are a
$\cC^l$-family of $\cC^k$-functions, and we may bound the $t$
and $x$ derivatives of ${\cW}(F_t)$ in terms of $F_t$.

\medskip
Using the above results, the proof of lemma 2.5 (sec.\ 2)
gives the following:

\proclaim 7.2 Lemma:
{\rm (Parameter-dependent rough multiplication.)}
Let $K \subset \R^n$ be a compact $\cC^1$-submanifold, with or
without boundary. Let $f_t$ be a $\cC^l$-family of $\cC^{k}$-functions
and $g_t$ a $\cC^l$-family of $\cC^{k-1}$-functions on a
neighborhood of $K$ in $\R^n$ such that $f_t = 0$ on $K$ for
each $t\in [0,1]$. Then there exists a $\cC^l$-family of $\cC^k$-functions
$h_t$ on a \nbd\ of $K$ such that $|h_t - f_t g_t| = o(d_K^k)$,
uniformly in $t\in [0,1]$. If $K\subset \C^n$ and if $f_t$, $g_t$
are $\dibar$-flat on $K$, then so is $h_t$.

We next prove an extension lemma.

\proclaim 7.3 Lemma: Let $M \subset \C^n$ be a compact, totally
real, $\cC^k$-submanifold. For any $\cC^l$-family of
$\cC^k$-maps $f_t \colon M \to \cN$ $(t\in [0,1])$
there exist an open set $U\subset \C^n$ containing $M$
and a $\cC^l$-family of $\cC^k$-maps $\wt{f}_t \colon U \to \cN$
such that each $\wt{f}_t$ is $\dibar$-flat on $M$ and it restricts
to $f_t$ on $M$. If $N=n$ and $f_t \colon M\to M_t=f_t(M) \subset \C^n$
is a diffeomorphism for each $t\in [0,1]$, we can choose
$\wt{f}_t$ as above to be a $\cC^l$-family of
$\cC^k$-diffeomorphisms on $U$.


\demo Proof of Lemma 7.3:
Let $m=\dim_\R M\le n$. We consider first the case when $M=\bar V$
is a smoothly bounded compact domain in $\R^m \subset \C^m \subset \C^n$.
Write $z_j = x_j + iy_j$ with $x_j,y_j \in \R$.
Given $f \in \cC^k(\bar{V})$, we consider the following
Whitney function on $\bar V$ for the real coordinates
$x_1,\ldots,x_m,y_1,\ldots,y_m$ in $\C^m$:
$$  F \colon  f_{(\a',\a'')} = i^{|\a''|}\, \di_x^{\a' + \a''}\!(f),
    \qquad \a',\a''\in\Z_+^m,\ |\a'| + |\a''| \le k.
$$
From the Cauchy-Riemann equations
${\di g \over \di y_j} = i {\di g \over \di x_j}$ $(1\le j \le m)$
for a function $g$ in a \nbd\ of $\bar V$ in $\C^m$
it follows that the Whitney extension $\wt{f} = {\cW}(F)$ of $F$
to $\C^m$ is $\dibar$-flat on $\bar{V}$. If $m < n$, we extend
${\cW}(F)$ trivially in the variables $z_{m+1},\cdots,z_n$
to get a Whitney extension on $\C^n$. Moreover, if
$\{f_t \colon t\in [0,1]\}$ is a $\cC^l$-family of $\cC^k$-functions
on $\bar{V}$ and we define $F_t$ as above, the Whitney extensions
${\cW}(F_t)$ are a $\cC^l$-family of $\cC^k$-functions which are
$\dibar$-flat on $\bar V$.

Next we consider a local $\cC^k$-parametrization
$\phi \colon U \to M$ around a point $w_0\in M$,
where $U$ is an open set in $\Rm$. Let
$z_0 =\phi^{-1}(w_0)\in U$. Choose a smoothly bounded domain
$V\ss U$ containing $z_0$ and set $W=\phi(\bar V)\subset M$.
Let $\wt{\phi}$ be an extension of $\phi$ to $\C^n$ constructed
above which is $\dibar$-flat on $\bar V$. If $m<n$, we also choose
a basis $v_{1},\cdots,v_{n-m}$ of the complex normal space
$(T_{w_0}^{\C}M)^{\bot}$ to $M$ at $w_0$. The map
$\Phi(z)= \wt{\phi}(z)+ \sum_{j=1}^{n-m} z_{m+j} v_j$
is then a $\cC^k$-diffeomorphism in a \nbd\ of $z_0$
which is $\dibar$-flat on $\bar V$; hence its inverse
$\Phi^{-1}$ is well defined in a \nbd\ $\wt{W} \subset \C^n$
of $w_0$ and is $\dibar$-flat on $W\cap \wt{W} \subset M$.

The first part of the proof also provides an extension $\psi_t$
of the map $f_t \circ \phi \colon \bar V\to \C^n$ to a \nbd\
of $\bar V$ in $\C^n$ such that $\psi_t$ is $\dibar$-flat on $\bar V$.
The composition $\psi_t\circ \Phi^{-1} \colon \wt{W}\to\C^n$
is a $\cC^k$-extension of the map $f_t$ which is $\dibar$-flat
on $W\cap \wt{W} \subset M$.

This gives us a local $\dibar$-flat $\cC^k$-extension of $f_t$
in a \nbd\ of each point $w_0\in M$. We can patch these local extensions
by a $\dibar$-flat partition of unity along $M$ as in lemma 2.6
to obtain a desired $\cC^l$-family $\wt{f}_t$ satisfying lemma 7.3.

It remains to consider the case when $f_t \colon M\to M_t$
is a diffeomorphism for each $t\in [0,1]$. Let
$\wt{M} = \bigcup_{t \in [0,1]} \{t\}
  \times M_t \subset [0,1]\times \cn$,
and let $\wt{f} \colon [0,1] \times M \to \wt{M}$ be the map
$\wt{f}(t,z) \to (t,f_t(z))$. Let $\nu$ denote the complex
normal bundle of $M$ and $\nu^t$ the complex normal bundle of $M_t$
in $\C^n$. Then $\wt{\nu} = \bigcup_{t \in [0,1]} \{t\} \times \nu^t$
is, in an obvious way, a vector bundle over $\wt{M}$, and
$[0,1] \times \nu$ is a vector bundle over $[0,1] \times M$.
By standard bundle theory (see Lemma 1.4.5 of [Ati]) there exists
a bundle equivalence $\psi \colon [0,1] \times \nu \to \wt{\nu}$
over $\wt{f}$. Thus we have continuously varying isomorphisms
$\nu_z \to \wt{\nu}^t_{f_t(z)}$ $(z\in M,\ t\in [0,1])$
which we extend to a continuous map
$A' \colon [0,1] \times M \to \End_{\bf C}(\cn)$. Then we
approximate $A'$ by a $\cC^l$-family of $\cC^k$-maps
$A\colon [0,1]\times M \to \End_{\bf C}(\cn)$ so that
$A(t,z)\nu_z$ is a supplementary subspace to $(Df_t)_z(T_z^C M)$
for each $(t,z)\in [0,1]\times M$. Let $L(t,z)$ equal
$(Df_t)_z^C$ on $T_z^C M$ and $A(t,z)$ on $\nu_z$.
Since $T_z\C^n= T_z^C M \oplus \nu_z$, $L(t,z)$ belongs
to $GL(n,\C)$, and it is not hard to check that
$L_t=L(t,\cdotp) \colon M\to GL(n,\C)$ is a
$\cC^l$-family of $\cC^{k-1}$-maps extending $Df_t$.
Using Lemma 7.2 it is easy to see that lemma 2.6 has a
parameter-dependent version which gives the desired
conclusion.
\endpr

\demo Proof of Theorem 7.1: Set $M=M_0$ and $i=i_0\colon M\hra\C^n$.
We first apply Lemma 7.3 to get a \nbd\ $U\subset \C^n$
of $M$ and a continuous family of $\cC^k$-diffeomorphisms
$\wh{f}_t \colon U \to U_t \subset \C^n$ which are $\dibar$-flat
on $M$. The family of inverses $(\wh{f}_t)^{-1} \colon U_t \to U$
is then a continuous family of $\cC^k$-diffeomorphisms
on $\wt{U} = \bigcup_{t \in [0,1]} \{t\} \times U_t$
which are $\dibar$-flat on $M_t$ and which extend
$f_t^{-1}\colon M_t\to M$.

Let $\a_t = \sum_{|I|=p} \a_{I,t}dz^I$ be as in theorem 7.1,
with $\a_{I,t} \in \cC^k(M_t)$. Our assumption is that
$\a_{I,t} \circ f_t$ $(t\in [0,1])$ is a continuous family
of $\cC^k$-functions for each $I$. Applying lemma 7.3
we can extend it to a continuous family $\a'_{I,t}$
of $\cC^k$-functions on $[0,1]\times U$ which are $\dibar$-flat
on $M$. Set $\wt{\a}_{I,t}= \a'_{I,t}\circ (\wh{f}_t)^{-1}$
and $\wt{\a}_t = \sum_{|I| = p} \wt{\a}_{I,t}dz^I$; this is
a continuous family of $\cC^k$ $(p,0)$-forms on $\wt{U}$,
and $\wt{\a}_t$ is $\dibar$-flat on $M_t$.

The next step is to modify $\wt{\a}_t$ so as to make its
differential flat on $M_t$. We observe that both
$\wh{f}_t^*\wt{\a}_t$ and
$\b_t \colon = d\wh{f}_t^*\wt{\a}_t = \wh{f}_t^*(d\wh{\a}_t)$
are continuous families of $\cC^{k-1}$-forms on $U$.
By assumption, $di_t^*\a_t = 0$, hence $i^*\b_t = 0$.

It is clear that the proof of lemma 6.3 produces a $\cC^l$-family
of solutions $v_t$ for any $\cC^l$-family $u_t$ satisfying the assumptions
in that lemma. Applying this to the forms $u_t=\b_t$ constructed
above, we obtain a continuous family of $(p,0)$-forms
$\g'_t = \sum_{j=1}^N \z_j \g_{j,t}^{\prime}$ $(t\in [0,1])$
such that $d\g'_t - \b_t$ is $(k-1)$-flat on $M$,
where $\zeta_1,\cdots,\zeta_N$ are $\dibar$-flat $\cC^k$-functions
vanishing on $M$ and $\g_{j,t}^{\prime}$ are continuous families
of $\cC^{k-1}$ $(p,0)$-forms that are $\dibar$-flat on $M$.

Then $(\wh{f}_t^{-1})^* \g_{j,t}' = \sum_{|I|=p} a_{j,I,t}dz^I + \l_{j,t}$
where $a_{j,I,t}$ are continuous families of $\cC^{k-1}$-functions that are
$\dibar$-flat on $M_t$ and $\l_{j,t} = o(d_{M_t}^{k-1})$ uniformly in $t$.
Applying parameter dependent rough multiplication (lemma 7.2) to $\z_j$
and $a_{j,I,t}\circ \wh{f}_t$ gives continuous families $b_{j,I,t}$
of $\cC^k$-functions near $M$ which are $\dibar$-flat on $M$.
Setting $\g_t = \sum_{|I|=p} \sum_{j=1}^N (b_{j,I,t}\circ \wh{f}_t^{-1}) dz^I$
and $\wh{\a}_t = \wt{\a}_t - \g_t$, we get $\wh{\a}_t \mid_{M_t} = \a_t$,
$t \in [0,1]$, and $|d\wh{\a}_t| = o(d_{M_t}^{k-1})$ uniformly in $t$.

The next step is to approximate $\wh{f}_t$ well by biholomorphic maps in
tubular neighborhoods ${\cal T}_\d$ of $M$. $\wh{f}_t$ maps $M$ onto $M_t$ and is
a diffeomorphism from a neighborhood $U$ of $M$ on a neighborhood $U_t$
of $M_t$, with estimates on derivatives valid for all $t \in [0,1]$.
It follows that for some $\bar{a}>0$ and all sufficiently small $\d>0$
we have $\wh{f}_t({\cal T}_{\bar{a}\d}M) \subset {\cal T}_{\d}M_t$
and $\wh{f}_t^{-1}({\cal T}_{\bar{a}\d}M_t) \subset {\cal T}_{\d}M$ for
all $t\in [0,1]$.

If we apply the solution operator of theorem 3.1 to the equation
$\dibar R_t^\d = \dibar \wh{f}_t$ in ${\cal T}_\d={\cal T}_\d M$ and set
$h_t^\d = \wh{f}_t - R_t^\d$, we obtain a continuous family
of holomorphic maps $h_t^\d$ on ${\cal T}_\d$ satisfying
$\|h_t^\d - \wh{f}_t\|_{\cC^j({\cal T}_\d M)} = o(\d^{k-j})$ for
$j \le k$, where $k \ge 2$. It follows that for small $\d>0$
the map $h_t^\d$ is a biholomorphism of ${\cal T}_\d$ onto its image,
and $g_t^\d := \wh{f}_t^{-1}\circ h_t^\d$ is a
$\cC^k$-diffeomorphism of the tube ${\cal T}_\d$ onto a small
perturbation of ${\cal T}_\d$.

Since $h_t^\d$ is close to $\wh{f}_t$, it is not hard to see, using the
argument in the proof of theorem 1.2,  that if $0<a<\bar{a}$ and
$\e > 0$ are given then for $\d>0$ sufficiently small
(depending on $a$ and $\e$) we have the inclusions
$h_t^\d({\cal T}_{\d^{\prime}}M) \supset {\cal T}_{a\d^{\prime}}M_t$
for $\e\d\leq\d^{\prime}\leq\d$ and
$(h_t^\d)^{-1}({\cal T}_{\d^{\prime}}M_t) \supset {\cal T}_{a\d^{\prime}}M$
for $\e\d\leq\d^{\prime}\leq a\d$. We also have
${\cal T}_{d^{\prime}/2}M \subset g_t^\d(T_{d^{\prime}}M)
\subset {\cal T}_{2d^{\prime}}M$
for $\e\d\leq\d^{\prime}\leq\d$, for all $t$.

The next step is to approximate $\wh{\a}_t$ by a continuous family of
holomorphic $p$-forms $u_t^\prime (= u_t^{\prime \d})$ on tubes
${\cal T}_\d M_t$. Suppose that $\wh{\a}_t = \sum_{|I| = p} \wh{\a}_{t,I} dz^I$.
For small $\d>0$, $h_t^{\d/a}({\cal T}_{\d/a}M) \supset {\cal T}_\d M_t$ for
$t \in [0,1]$. Let $u_{t,I}^{\prime \prime}$ be holomorphic approximations
to $\wh{\a}_{t,I} \circ h_t^{\d/a}$, constructed as $F_\d$ in
section 4. Set
$u_{t,I}^\prime = u_{t,I}^{\prime \prime} \circ (h_t^{\d/a})^{-1}$.
Then the $p$-form $u_t^{\prime} = \sum_{|I| = p} u_{t,I}^\prime dz^I$
is holomorphic in ${\cal T}_\d M_t$ and satisfies
$\|u_t^\prime - \wh{\a}_t \|_{\cC^j({\cal T}_\d M_t)} = o(\d^{k-j})$,
uniformly in $t$. We also see that
$\|du_t^\prime \|_{L^{\infty}({\cal T}_\d M_t)} = o(\d^{k-1})$, and
if we set $v_{0,t} = i_t^* (u_t^\prime - \wh{\a}_t)$ then
$dv_{0,t} = i_t^*du_t^\prime$.

We wish to prove the existence of a continuous family of holomorphic
$(p-1)$-forms $v_t (=v_t^\d)$ on ${\cal T}_{b\d}M_t$ for some $b>0$, with
$\|v_t\|_{L^\infty ({\cal T}_{b\d}M_t}) = o(\d^k)$, uniformly in $t$, and solving
$dv_t = du_t^\prime$. Then $u_t^\d = u_t^\prime - v_t$ would be a continuous
family of closed holomorphic $p$-forms with
$\|u_t|_{M_t} - \a_t \|_{\cC^j(M_t)} = o(\d^{k-j})$,
uniformly in $t$, as required.

A parameter-dependent version of theorem 5.1 for the family $M_t$
would give that result. The following argument will give this for
a small $b>0$, but we shall restrict ourselves to the special case
we need. Choose $a < \bar{a}$ and $\e = a/2$. For $\d>0$
small, $w_t^\prime = \wh{f}_t^*(du_t^{\prime \d})$ are
$\cC^{k-1}$-forms on ${\cal T}_{\bar{a} \d}M$ with
$\|w_t^\prime \|_{L^\infty({\cal T}_{\bar{a}\d}M)} = o(\d^{k-1})$ and
$\|w_t^\prime \|_{\cC^s({\cal T}_{\bar{a}\d}M)} = o(\d^{k-1-s})$,
uniformly in $t$.

Furthermore, with $v_{0,t}^\prime = f_t^*v_{0,t}$, we have
$dv_{0,t}^\prime = i^*w_t^\prime$ on $M$, with
$\|v_{0,t}^\prime \|_{L^\infty} = o(\d^k)$ and
$\|v_{0,t}^\prime \|_{\cC^s} = o(\d^{k-s})$, uniformly in $t$.
Then the first part of the proof of theorem 5.1 and the remarks
on continuous $t$-dependence give a continuous family of
$\cC^{k-1}$-forms $\omega_t^\prime$ on ${\cal T}_{\bar{a}\d}M$
solving $d\omega_t^\prime = w_t^\prime$,
with $\|\omega_t^\prime\|_{L^\infty} = o(\d^k)$ and
$\|\omega_t^\prime \|_{\cC^s} = o(\d^{k-s})$, uniformly in $t$.
Then $\omega_t = (g_t^\d)^* \omega_t^\prime$ are defined
on ${\cal T}_{\bar{a}\d/2}M$ and satisfy the same kind of estimates, and
$d\omega_t = (h_t^\d)^*(\wh{f}_t^{-1})^*w_t^\prime = (h_t^\d)^*du_t^\prime$
is holomorphic. Since $a<\bar{a}$, the second part of the proof
of theorem 5.1 gives the existence of a continuous family of
holomorphic $p$-forms $v_t^\prime$ on ${\cal T}_{a\d/2}M$ satisfying
$dv_t^\prime = (h_t^\d)^*du_t^\prime$ and
$\|v_t^\prime \|_{L^\infty ({\cal T}_{a\d/2}M)} = o(\d^k)$, uniformly in $t$.
By assumption $h_t^\d({\cal T}_{a\d/2}M) \supset {\cal T}_{a^2\d/2}M_t$ for each t,
and $v_t^\d = (h_t^\d)^{-1*}v_t^\prime$ is a continuous family
of holomorphic $p$-forms on ${\cal T}_{a^2\d/2}M_t$ with $dv_t^\d = du_t^\prime$
on ${\cal T}_{a^2\d/2}M_t$ and $\|v_t^\d \|_{L^\infty} = o(\d^k)$ uniformly in $t$.

We now show that if $i_t^*\a_t$ is exact for every $t$, the holomorphic
forms $u_t^\d$ as above may be chosen to be exact. We recall that
the de Rham cohomology group $H^p(M,{\bf C})$ is finite dimensional
and $f_t^* \colon H^p(M_t,{\bf C}) \to H^p(M,{\bf C})$ is an isomorphism
for every $t$. We have that
$H^p(M,{\bf C}) \approx \{\a \in \cC_{(p)}(M) : d\a = 0\}/$(exact forms),
where derivatives are taken in the weak sense, and we may equip
$H^p(M,{\bf C})$ with the quotient norm.

For each $t_0 \in [0,1]$ there exist closed holomorphic $p$-forms
$\wh{u}_1,\cdots,\wh{u}_N$ on an open neighborhood $U$ of $M_{t_0}$
such that $[i_{t_0}^*\wh{u}_j]$, $1\leq j\leq N$,
is a basis for $H^p(M_{t_0},{\bf C})$. Then $t \to [f_t^*u_t^\d]$
is a continuous map $[0,1] \to H^p(M,{\bf C})$, and
$t \to [f_t^*\wh{u}_j]$ is continuous for $t$ neat $t_0$ and
$1\leq j\leq N$. It follows that $\{[f_t^*\wh{u}_j] \colon j\leq N\}$
is a basis for $H^p(M,{\bf C})$ for $t$ in a neighborhood $J\subset [0,1]$
of $t_0$, and that we may write
$[f_t^*u_t^\d] = \sum_{j=1}^N c_j^\d(t) [f_t^*\wh{u}_j]$ with $c_j^\d$
continuous on $J$. Each form $f_t^*\a_t$ is exact on $M$, so
$\|[f_t^*u_t^\d]\|\leq\|f_t^*(u_t^\d)-\wh{\a}_t)\|_{L^\infty(M)} = o(\d^k)$.
This means that for $J_1\ss J$ we have
$\max_{t\in J_1} |c_j^\d(t)| = o(\d^k)$ for all $j\leq N$.
For $\d>0$ small and $t\in J_1$ we have ${\cal T}_\d M_t \subset U$,
$u_t^{0\d}=u_t^\d - \sum_{j=1}^N c_j^\d(t)\wh{u}_j$ is exact on
${\cal T}_\d M_t$ (since $[i_t^*u_t^{0\d}] = 0$), and it approximates
$\a_t$ well enough. We can now patch these together with a partition
of unity in $t$ to obtain a solution $u^\d_t$ for $t\in [0,1]$
satisfying theorem 7.1.

Finally, assume $M_t$ is polynomially convex for all $t\in [0,1]$
and let $u_t^\d$ be the exact solution on $U_\d$. For $\d>0$ sufficiently
small we may also assume that ${\cal T}_\d M_t$ is Runge in $\C^n$ for all $t$.
Given $a < a^\prime < 1$ and $\e > 0$, there exist $t_j \in [0,1]$,
$j=1,\cdots,N$, and (relatively) open intervals $I_j \subset [0,1]$,
$t_j \in I_j$, such that
$U_{a\d} \subset \bigcup_{j=1}^N I_j \times
{\cal T}_{a^\prime \d}M_{t_j} \subset U_\d$,
and for all for $t \in I_j$ we have
$\|u_t^\d - u_{t_j}^\d \|_{\cC^k({\cal T}_{a^\prime \d}M_{t_j})} < \e$ and
$\|\wh{\a}_t^\d - \wh{\a}_{t_j}^\d \|_{\cC^k({\cal T}_{a^\prime \d}M_{t_j})} < \e$.
Let $\b_j$ be a holomorphic $(p-1)$-form on ${\cal T}_\d M_{t_j}$ such that
$d\b_j = u_{t_j}^\d$. By Oka's theorem there is an entire $(p-1)$-form
$v_j$ such that $\| \b_j-v_j \|_{L^\infty({\cal T}_\d M_{t_j})} < \e$.
The Cauchy estimates imply
$\|\b_j-v_j\|_{\cC^r({\cal T}_{a^\prime \d}M_{t_j})} = \e o(\d^{-r})$,
and hence
$\|\u_{t_j}^\d -dv_j\|_{\cC^r({\cal T}_{a^\prime \d}M_{t_j})} = \e o(\d^{-(r+1)})$.
Choosing $\e = o(\d^{k+1})$, we obtain
$\|dv_j - \wh{\a}_t^\d \|_{\cC^r({\cal T}_{a\d}M_t)} = o(\d^{k-r})$ whenever $t \in I_j$.
If $\chi_j(t)$ is a partition of unity on $[0,1]$ subordinate to the covering
$\{ I_j \}$ and we define  $v_t=\sum_{j=1}^N \chi_j(t)v_j(z)$, then
$u_t = dv_t$ is an entire form for each $t$ which satisfies theorem 7.1.
\endpr

\demo Proof of Theorem 1.7:
By assumption $f_t \colon M \to M_t$ is $\cC^1$-family of
$\cC^k$-diffeomorphisms and $\omega$ is one of the forms
(1.4), (1.5). Let $X_t$ be the infinitesimal generator of
$f_t$, i.e., $\di_t f_t(z) = X_t(f_t(z))$ for $z\in M$ and
$t\in [0,1]$. Then $\a_t = X_t\rfloor \omega$ is a continuous family of
$(p,0)$-forms on $M_t$, with $p=n-1$ when $\omega$ is the volume
form (1.4)) and $p=1$ when $\omega$ is the symplectic form (1.5).
Since $f_t$ is an $\omega$-flow, $i_t^*\a_t$ is closed on $M_t$
for each $t$, by the remark after definition 2.

By theorem 7.1 there exists an extension of $\a_t$ to a continuous
family $\wh{\a}_t$ of $(p,0)$-forms of class $\cC^k$ on a neighborhood of
$\wt{M} = \bigcup_{t \in [0,1]} \{t\} \times M_t$ such that for all
sufficiently small $\d > 0$ there exists a continuous family of closed
holomorphic $p$-forms $u_t^\d$ on
$U_{\d} = \bigcup_{t \in [0,1]} \{t\} \times {\cal T}_{\d}M_t$
with $\|u^{\d}_t - \wh{a}_t \|_{\cC^r({\cal T}_\d M_t)} = o(\d^{k-r})$,
uniformly in $t$, for $0 \le r \le k$.

The equation $u_t^\d = Y_t^\d \rfloor \omega$ uniquely defines a
time-dependent holomorphic vector field $Y_t^\d$ on $U_\d$.
Since $u_t^\d$ is closed, the flow $F_t^\d$ of $Y_t^\d$ is
a holomorphic $\omega$-flow wherever it is defined
(see definition 2). If we let $X_t$ denote the extension
of $X_t$ to $U_\d$ defined by
$\wh{\a}_t = X_t \rfloor \omega$, then
$\|Y_t^\d - X_t \|_{\cC^r({\cal T}_\d M_t)} = o(\d^{k-r})$,
uniformly in $t$. We may apply lemma 4.1 of [FL] to see that
for small $\d>0$ the flow $F_t^\d (z)$ exists for all $t \in [0,1]$
and $z \in {\cal T}_\d M_0$, and
$\| F_t^\d - f_t \|_{\cC^r({\cal T}_\d M_0)} = o(\d^{k-r})$,
uniformly in $t$. In fact, it follows from the proof of this lemma
(see section 4 and [FL]) that the same approximation also holds
for the flow from time $t$ to time $s$; if we let
$f_{t,s}=f_s\circ f_t^{-1} \colon {\cal T}_\d M_t \to \cn$ denote the flow of
$X_t$ from $t$ to $s$ and $F^\d_{t,s}=F^\d_s\circ (F^\d_t)^{-1}$ the flow
of $Y_t^\d$ from $t$ to $s$, then for small $\d>0$ the flow
$F^\d_{t,s}$ exists for all $s,t \in [0,1]$, and we have
$\|F^\d_{t,s}-f_{t,s}\|_{\cC^r({\cal T}_\d M_t)} = o(\d^{k-r})$, uniformly in
$s$ and $t$. Since $f_t^{-1} = f_{t,0}$, the second estimate
in theorem 1.7 follows.

Finally, if $f_t$ is an exact $\omega$-flow, i.e., $i_t^*\a_t$ is
exact on $M_t$ for each $t$, and if each $M_t$ is also polynomially convex,
then by (the proof of) theorem 7.1 above we may choose
$u_t^\d (z) = \sum_{j=1}^N \chi_j^\d (t)dv_j(z)$, where $v_j(z)$
are entire $(p-1)$-forms on $\C^n$ and $\chi_j^\d$
$(1\le j\le N)$ are $\cC^\infty$ functions with compact support
in $\R$ which form a partition of unity on $[0,1]$. We may even
assume that $v_j$ are $(p-1)$-forms with polynomial coefficients.
This means that the polynomial vector fields $X_j$ on $\C^n$,
uniquely defined by the equation $dv_j = X_j \rfloor \omega$,
are divergence free (resp.\ Hamiltonian). By proposition 4.1
in [F4] these can be written as finite sums
$X_j(z) = \sum_{k=1}^{N_j} X_{j,k}(z)$, where
$X_{j,k}$ are complete divergence free (resp.\ Hamiltonian)
polynomial vector fields on $\C^n$ (in fact they are shear fields).
Completeness means that the fields $X_{jk}$ may be
integrated in time for all $t \in \C$ (and initial points $z\in \C^n$).
Then $Y_{jk}(t,z) \colon = \chi_j^\d (t) X_{jk}(z)$ is also a complete
vector field whose integral curves are reparametrizations of the
integral curves of $X_{jk}$. Hence we may write
$Y_t^\d = \sum_{j,k} Y_{jk}(t,z)$, i.e., $Y_t^\d$ is the sum of
complete, divergence free (resp.\ Hamiltonian), time-dependent,
polynomial (in $z\in \C^n$) vector fields. For the rest of this
proof it is more convenient to write this sum as
$\sum_{l=1}^N Y_l (t,z)$, where each $Y_l$ is one
of the $Y_{jk}$ above.

Let $G_{t,t+s}^l$ be the flow of $Y_l (t,z)$ from time $t$ to time
$t+s$. This means that $G_{t,t}^l(z) = z$ and ${d \over ds}G_{t,t+s}^l(z) =
Y_l(t+s,G_{t,t+s}^l(z))$. Define
$G_{t,t+s}(z) = (G_{t,t+s}^N \circ \cdots \circ G_{t,t+s}^1)(z)$. We
can regard this as the flow of a time-dependent vector field
$X^{t,t+s}_{t^\prime}$, defined for times $t^\prime$ between $t$
and $t+s$; for $t+{j-1 \over N}s\leq t^\prime\leq t+{j \over N}s$ we define
$X^{t,t+s}_{t^\prime}(z)= {1\over N}Y_j(t+N(t^\prime-{j-1 \over N}s),z)$.
If we reparametrize time such that the joints are passed at zero speed,
we may even assume that $X^{t,t+s}_{t^\prime}$ is smooth and vanishes near the
endpoints. We denote this smooth flow by $G^{t,t+s}_{t'}(z)$. By definition,
$G_{t,t+s}(z)=G^{t,t+s}_{t+s}(z)$. Since the vector fields $Y_j$ are complete
divergence free (resp.\ Hamiltonian) entire vector fields, it follows that
$G^{t,t+s}_{t^\prime}$ is a \holo\ $\omega$-flow, i.e.,
$(G^{t,t+s}_{t^\prime})^*\omega = \omega$ when $t\le t'\le t+s$.

For each $m \in \N$ we define the concatenations
$F_1^m(z) = (G_{1- {1\over m},1} \circ \cdots \circ G_{0,{1 \over m}})(z)$.
Then by supplement 4.1.A of [AMR] we have
$\lim_{m\to\infty} F_1^m(z)= F_1^\d(z)$, uniformly for $z\in {\cal T}_\d M_0$.
As above, we can view $F^m_1(z)$ as the time-one map of the flow
of the vector field $X_t$ defined by
$X_t=X_t^{{j-1\over m},{j\over m}}$ for
$t \in [{j-1\over m},{j\over m}]$, $1\leq j\leq m$.
Let $F^m_t(z)$ be the flow of this vector field. It is easy
to see that we can arrange that $\lim_{m\to\infty} F^m_t =F_t^\d$,
uniformly in $[0,1] \times {\cal T}_\d M_0$, and the Cauchy estimates
imply $\|F_t^m - F_t^\d\|_{\cC^k(M_0)} < \e$ for all $t\in [0,1]$
and all sufficiently large $m\in \N$.
Similarly, $(F_1^m)^{-1}$ is a concatenation and hence
$\lim_{m\to\infty} (F_1^m)^{-1}= (F_1^\d)^{-1}$ uniformly on ${\cal T}_\d M_1$;
it follows that $\lim_{m\to\infty} (F^m_t)^{-1} =(F^\d_t)^{-1}$ on $M_t$,
hence the result follows by the Cauchy estimates.
\endpr

\demo Proof of Theorem 1.8:
We shall see that in all cases except (iii) and (vi) the
pull-back $i_t^*\a_t$ of the form $\a_t=X_t\rfloor\omega$ to $M_t$
is exact for each $t$; hence $f_t$ is an exact $\omega$-flow
and the result follows from the second part of theorem 1.7.

In case (i) we have $i_t^* \a_t=0$ by degree reason.
In cases (ii), (iv), (v) and (vii) we first see that the form
$i_t^*\a_t$ is closed on $M_t$, either by degree reasons or
by the comment after definition 2 in sect.\ 1; hence
the cohomological assumptions imply in each of these cases
that $i_t^*\a_t$ is exact on $M_t$.

For the two remaining cases (iii) and (vi) it is shown on pages 439
and 441 of [F3] that the initial family $f_t$ may be altered to an
exact, \tr\ and \pc\ $\omega$-flow, without changing the maps
$f_0=Id$ and $f_1$; hence the result again follows from theorem 1.7.
\endpr

\leftline{\bf References.}
\medskip

\ii{[AMR]} A.\ Abraham, R.\ E.\ Marsden, T.\ Ratiu:
Manifolds, Tensor Analysis, and Applications.
(Second ed.) Applied Math.\ Sciences {\bf 75},
Springer, New York--Berlin, 1988.

\ii{[AB]} M.\ Andersson, B.\ Berndtsson:
Henkin--Ramirez formulas with weight factors.
Ann.\ Inst.\ Fourier (Grenoble), {\bf 32}, 91--110 (1982).

\ii{[Ati]} M.\ F.\ Atiyah: K-Theory.
(Lecture notes by D.\ W.\ Anderson.)
Benjamin, New York--Amsterdam, 1967.
Second ed.: Addison-Wesley, Redwood City, CA, 1989.

\ii{[BT1]} M.\ S.\ Baouendi, F.\ Tr\`eves:
A property of the functions and distributions annihilated by a locally
integrable system of complex vector fields.
Ann.\ of Math.\ (2) {\bf 113}, 387--421 (1981).

\ii{[BT2]}  M.\ S.\ Baouendi, F.\ Tr\`eves:
Approximation of solutions of linear PDE with analytic coefficients.
Duke Math.\ J.\ {\bf  50}, 285--301 (1983).

\ii{[BB]}  J.\ Bruna, J.\ Burgu\'es:
Holomorphic approximation in $\cC^m$-norms on totally real
compact sets in $\C^n$.
Math.\ Ann.\ {\bf 269}, 103--117 (1984).

\ii{[Ca]} H.\ Cartan: Espaces fibr\'es analytiques.
Symposium Internat.\ de topologia algebraica, Mexico, 97--121 (1958).
(Also in Oeuvres, vol.\ 2, Springer, New York, 1979.)

\ii{[DF]} K.\ Diederich, J.\ E.\ Forn\ae ss:
A smooth curve in ${\C}^2$ which is not a pluripolar set.
Duke Math.\ J.\ {\bf 49}, 931--936 (1982).

\ii{[F1]} F.\ Forstneri\v c:
Approximation by automorphisms on smooth submanifolds of $\C^n$.
Math.\ Ann. {\bf 300}, 719--738 (1994).

\ii{[F2]} F.\ Forstneri\v c: A theorem in complex symplectic geometry.
J.\ Geom.\ Anal.\ {\bf 5}, 379--393 (1995).

\ii{[F3]} F.\ Forstneri\v c: Equivalence of real submanifolds under
volume preserving holomorphic automorphisms of $\cn$.
Duke Math.\ J.\ {\bf 77}, 431--445 (1995).

\ii{[F4]} F.\ Forstneri\v c:
Holomorphic automorphisms of $\cn$: A survey.
Complex Analysis and Geometry (Trento, 1993), 173--200,
Lecture Notes in Pure and Appl.\ Math., \bf 173\rm,
Marcel Dekker, New York, 1996.

\ii{[FL]} F.\ Forstneri\v c, E.\ L\o w:
Global Holomorphic Equivalence of Smooth Submanifolds in $\cn$.
Indiana Univ.\ Math. J.\ {\bf 46}, 133--153 (1997).

\ii{[FR]} F.\ Forstneri\v c, J.-P.\ Rosay:
Approximation of biholomorphic mappings by automorphisms of ${\C}^n$.
Invent.\ Math.\ {\bf 112}, 323--349 (1993).

\ii{[GLR]} I.\ Gohberg, P.\ Lancaster, L.\ Rodman:
Invariant subspaces of matrices with applications.
John Wiley and Sons, New York, 1986.

\ii{[HaW]} F.\ R.\ Harvey, R.\ O.\ Wells, Jr.:
Holomorphic approximation and hyperfunction theory on a
$\cC^1$ totally real submanifold of a complex manifold.
Math.\ Ann.\ {\bf 197}, 287--318 (1972).

\ii{[HL]} G.\ M.\ Henkin, J.\ Leiterer:
Theory of Functions on Complex Manifolds.
Monographs in Mathematics {\bf 79}, Birkhäuser, Boston, 1984.

\ii{[H\"o]} L.\ H\"ormander:
An Introduction to Complex Analysis in Several Variables, 3rd ed.
North Holland, Amsterdam, 1990.

\ii{[H\"oW]} L.\ H\"ormander, J.\ Wermer:
Uniform approximations on compact sets in $\cn$.
Math.\ Scand.\ {\bf 23}, 5--21 (1968).

\ii{[M]} J.\ Moser: On the volume elements on a manifold.
Trans.\ Amer.\ Math.\ Soc.\ {\bf 120}, 286--294 (1965).

\ii{[\O1]}  N.\ \O vrelid:
Integral representation formulas and $L\sp{p}$-estimates for
the $\bar \partial $-equation.
Math.\ Scand.\ {\bf 29}, 137--160 (1971).

\ii{[\O2]} N.\ \O vrelid:
Integral representation formulae for differential forms,
and solutions of the $\overline \partial $-equation.
Fonctions analytiques de plusieurs variables et analyse complexe
(Colloq.\ Internat.\ CNRS, No.\ 208, Paris, 1972), pp. 180--198.
Gauthier-Villars, Paris, 1974.

\ii{[RS]}  R.\ M.\ Range, Y.\ T.\ Siu:
$\cC^k$ approximation by holomorphic functions and
$\bar \partial $-closed forms on $\cC^k$ submanifolds of
a complex manifold.
Math.\ Ann.\ {\bf 210}, 105--122 (1974).

\ii{[S]} E.\ M.\ Stein:
Harmonic analysis: real-variable methods, orthogonality, and
oscillatory integrals. (With the assistance of T.\ S.\ Murphy.)
Princeton Mathematical Series {\bf 43}, Princeton
University Press, Princeton, 1993.

\ii{[T]} J.\ C.\ Tougeron: Id\'eaux de fonctions diff\'erentiables.
Ergebnisse der Mathematik und ihrer Grenzgebiete {\bf 71},
Springer, Berlin--New York, 1972.

\ii{[We]} R.\ O.\ Wells: Differential analysis on complex manifolds.
(Second ed.) Graduate Texts in Math.\ {\bf 65},
Springer, New York--Berlin, 1980.

\ii{[Wh1]} H.\ Whitney: Analytic extensions of differentiable
functions defined on closed sets.
Trans.\ Amer.\ Math.\ Soc.\ {\bf 36}, 69--89 (1934).

\ii{[Wh2]} H.\ Whitney: Differentiable Manifolds.
Ann.\ of Math.\ {\bf 37}, 645--680 (1936).

\vfill\eject
\settabs 5\columns
\+\ \ Franc Forstneri\v c            &&& Erik L\o w and Nils \O vrelid \cr
\+\ \ Department of Mathematics      &&& Department of Mathematics \cr
\+\ \ University of Wisconsin        &&& University of Oslo \cr
\+\ \ Madison, WI 53706, USA         &&& P.O.Box 1053, Blindern \cr
\+                                   &&& N-0315 Oslo, Norway \cr

\+\ \ {\it Current address:}    \cr
\+\ \ IMFM, University of Ljubljana \cr
\+\ \ Jadranska 19   \cr
\+\ \ 1000 Ljubljana, Slovenia    \cr

\bye